\documentclass[11pt,psfig]{report}
\usepackage{amssymb,amsmath,amscd,amsthm}
\newtheorem{conj}{Conjecture}

\newtheorem{thm}{Theorem}[chapter]
\newtheorem{lem}[thm]{Lemma}
\newtheorem{cor}[thm]{Corollary}
\newtheorem{defn}{Definition}[chapter]
\newtheorem{no}{Notation}[chapter]
\newtheorem{pro}{Problem}
\newtheorem{thm*}{Theorem}
\newtheorem{lem*}[thm*]{Lemma}
\newcommand{\fb}{\framebox(10,10){\;}}
\newcommand{\fp}{\framebox(10,10){+}}
\newcommand{\fm}{\framebox(10,10){-}}
\newcommand{\mbo}{\makebox(10,10){\;}}
 
\newcommand{\f}[1]{\mathfrak{#1}}

\newcommand{\p}{\prime}
\newcommand{\Ind}{{\rm Ind}}
\newcommand{\Mat}{{\rm Mat}}
\newcommand{\spec}{{\rm{Spec}}}
\newcommand{\mb}{\mathbb}
\newcommand{\commentout}[1]{}

\newcommand{\mc}{\mathcal}
\newcommand{\arr}[1]{\left( \begin{array}{clcr} #1 \end{array} \right)}

\newcommand{\diag}{{\,\rm diag}}
\newcommand{\Hom}{\rm Hom\,}
\newcommand{\sgn}{\, {\rm sgn}}
\newcommand{\trivial}{{\rm trivial}}

\begin{document}
\title{Unipotent Representations and Quantum Induction}
\author{Hongyu He \\
Department of Mathematics \\
Louisiana State University \\
email: hongyu@math.lsu.edu \\
}
\date{Abstract: In this paper, we construct unipotent representations for the real orthogonal 
groups and the metaplectic groups in the sense of Vogan. Our construction is based on quantum induction which involves the  compositions of even number of theta correspondences. In particular, our results imply that
there are irreducible unitary representations attached to each special nilpotent orbit. }

\maketitle

\chapter*{Introduction}
\addcontentsline{toc}{chapter}{\bf{INTRODUCTION}  }

The existence of certain irreducible unitary representations, often called unipotent representations, was conjectured by Arthur, Barbasch and Vogan. Precise predictions about the representation-theoretic invariants associated with these representations were given in \!~\cite{arthur83}, \!~\cite{bv}, \!~\cite{vogan89} and \!~\cite{arthur}. The goal of this article is to verify the existence of unipotent representations for $O(p,q)$
and $Mp_{2n}(\mathbb R)$ in the sense of Vogan. The main idea is to study the theta correspondence and its compositions with respect to the dual pairs $(O(p,q), Sp_{2n}(\mathbb R))$ along the line of \!~\cite{howe} \!~\cite{li2} \!~\cite{theta} \!~\cite{unit} \!~\cite{basic}. Other classical orthogonal and symplectic groups can all be realized in dual pairs that are of type I in the sense of Howe \!~\cite{howe}, and unipotent representations for these groups  can be constructed in the same spirit. \\
\\ 
Unipotent representations play a crucial role in the representation theory of real reductive groups and in the theory of automorphic forms. 
Constructing unipotent representations is one core part of Vogan's program to classify the unitary duals of real reductive groups (\!~\cite{vogan86}).  For general linear groups over Archimedean fields, unipotent representations can be produced by parabolic induction and the classification of unitary duals was carried out by Vogan himself in \!~\cite{vo86}. For complex semisimple groups, unipotent representations were constructed and studied by Barbasch-Vogan (\!~\cite{bv}). Later, Barbasch proved the unitarity of these representations which led to his classification of the unitary duals of complex classical groups (\!~\cite{b89}). For other types of real reductive groups, Adams, Barbasch and Vogan gave certain descriptions about special unipotent representations in Arthur packets (\!~\cite{abv} and \!~\cite{arthur}). The unitarity and construction of these representations remain open problems. \\
\\
In this paper, we propose a different way of obtaining unipotent representations that will automatically be unitary. The motivation came from the work of Segal-Shale-Weil on the oscillator representations and the work of Howe, Li, Przebinda and many others on theta correspondences. In \!~\cite{basic}, we defined the concept of quantum induction as a composition of theta correspondences (see \!~\cite{howe}) within a certain range. Moreover we proved that quantum induction, if nonvanishing, preserves unitarity.
In this paper, we further study quantum induction in the framework of the orbit
method (\!~\cite{vogan86}, \!~\cite{vo94}). Using quantum induction, we attach to certain
real rigid nilpotent orbits $\mc O$ 
a packet of irreducible {\bf unitary} representations, $\mc N(\mc O)$. Each irreducible representation in $\mc N(\mc O)$ satisfies the characterizations given in \!~\cite{vogan89}. This construction gives rises to  many unipotent representations beyond what was known before.\\
\\
In this introduction, we will give a brief review of unitary representations, nilpotent orbits and how they are related for semisimple Lie groups. Then we will discuss theta correspondences, quantum induction before we state our main results.  
\section{Unitary dual, conjugacy classes and coadjoint orbits}
Let $G$ be a Lie group with finitely many connected components. Let $\f g$ be its Lie algebra and $\f g^*$ be the dual linear space of $\f g$ over $\mathbb R$. The group $G$ acts on the Lie algebra $\f g$ adjointly.  The dual action of $G$ on $\f g^*$ is often called the coadjoint action. \\
\\
Let $\Pi_u(G)$ be the equivalence classes of irreducible unitary representations of $G$. $\Pi_u(G)$ is  often called the unitary dual of $G$, also denoted by $\hat{G}$. $\Pi_u(G)$ can be made into a topological space using Fell topology (\cite{wallach}). For many groups, including most  classical groups discussed in this paper, $\Pi_u(G)$ is not completely classified. Historically, unitary dual was often studied by linking it to conjugacy classes of $G$, or coadjoint orbits in $\f g^*$. This link is already evident for finite  groups. For instance the number of elements in $\Pi_u(G)$ is exactly the number of conjugacy classes of $G$ (\cite{se}). For simply connected nilpotent groups, Kirillov proved that $\Pi_u(G)$ is in one-to-one correspondence with the coadjoint orbits $\f g^*// G$ (\cite{kir0}). Later Auslander and Kostant extended this correspondence to type I solvable groups (\cite{ak}). Since then, the link between coadjoint orbits and unitary representations has been explored by many authors and proven to be a powerful way to study unitary representations (\cite{kir}). This method is generally known as the orbit method.\\
\\
 Recall that every Lie group is a semidirect product of a reductive group and  a nilpotent group. By the virtue of Mackey's theory (\cite{mac}), the core of classifying $\Pi_u(G)$ for general Lie groups, is to classify  $\Pi_u(G)$ for $G$ reductive.  Furthermore a reductive Lie group is a direct product of its center and  a semisimple Lie group.  Hence we shall focus on semisimple Lie groups. \\
 \\
Assume that {\bf $G$ is connected and semisimple}. 
Let $\Pi_{2}(G)$ be the reduced dual of $G$, consisting of those irreducible unitary representations appearing in the Plancherel formula of $L^2(G)$ ( \cite{wallach}).  For $G$ compact, $\Pi_2(G)=\Pi_u(G)$ and $\Pi_u(G)$ is well-known,  by the theory of Cartan and Weyl. For $G$ noncompact,  due to the work of Harish-Chandra, Rossmann and others, there is a somewhat satisfactory way to link $\Pi_2(G)$ to a subset of  coadjoint orbits in $\f g^*/G$ (\cite{ro0}).  One special feature of a representation in $\Pi_2(G)$ is that its matrix coefficients are \lq\lq tempered \rq\rq. Thus $\pi \in \Pi_2(G)$ is also called tempered. The main trouble here is that
$\Pi_2(G)$ is a proper subset of $\Pi_u(G)$ when $G$ is noncompact. There are plenty of irreducible unitary representation that are not tempered. Roughly speaking, nontempered representations tend to be more "singular" and exhibit more symmetries.  Many nontempered representations arise naturally in harmonic analysis, quantum mechanics and  number theory. 
 \\
\\
{\it Fix a connected  noncompact semisimple Lie group $G$}. The easiest nontempered unitary representation is the trivial representation where the Lie algebra $\f g$ acts trivially. 
For $G$ the metaplectic covering of $Sp_{2n}(\mb R)$, the most well-known nontempered representation is perhaps the oscillator representation,  also known as the Segal-Shale-Weil representation.
 It is a unitary representation on $L^2(\mathbb R^n)$ and 
 the Lie algebra acts as algebraic differential operators on $\mathbb R^n$. Glancing at these two examples, one can see that the study of nontempered unitary representations requires deeper understanding of the action of the Lie algebra, consequently the action of universal enveloping algebra $U(\f g)$. To this end, one powerful tool  is developed to understand the  action of $U(\f g)$, namely,   $Ann(\pi) \subset U(\f g)$,  the annihilator of the representation $\pi$ (\cite{dix}). Using the natural filtration on $U(\f g)$ and some standard noncommutative algebra technique, one obtains the associated variety $\mc V(Ann(\pi)) \subseteq \f g_{\mathbb C}^*$. Here $\f g_{\mathbb C}$ is the complex Lie algebra and $\f g^*_{\mathbb C}$ is the complex dual space. Since $Ann(\pi)$ is stable under the coadjoint action of $G$, $\mc V(Ann(\pi))$ will be stable under the action of the adjoint group $G^{ad}=Ad(G)$, consequently stable under the adjoint action of the complexified  group $G^{ad}_{\mathbb C}$. For the trivial representation, the associated variety consists of one single point $\{ 0 \}$. For the oscillator representation, the associated variety will be the closure of the minimal nilpotent orbit of $Sp_{2n}(\mb C)$. Now we attach a geometric object $\mc V(Ann(\pi))$ to a unitary representation $\pi$. \\
 \\
 In 1980's, Borho-Brylinski and Joseph proved that for $\pi$ irreducible,
$\mc V(Ann(\pi))$ is the closure of a single (coadjoint) nilpotent orbit in $\f g^*_{\mathbb C}$  (\cite{bb} \cite{jo}). More details about nilpotent orbits will be given in the next section and it is enough to remind the reader that there are only finite number of nilpotent orbits in a semisimple Lie algebra. 
 Therefore, we can now attach to any $\pi \in \Pi_u(G)$  a unique (complex) nilpotent orbit. Conversely, we may  divide $\Pi_u(G)$ into packets of representations $\Pi_u(\mc O)$ that share the same associated variety and attach unitary representations to nilpotent orbits. For example, the trivial representation will be attached to the trivial nilpotent orbit $\{0\}$. More details about associated variety can be  found in Ch 1.5.  \\
\\
For many nilpotent orbits, there are infinitely many irreducible unitary representations attached to it. For other nilpotent orbits, it is expected that there are only finite number of (nonequivalent) irreducible unitary representations attach to it.  For some nilpotent orbits,  there may be no unitary representation attached to it (\cite{hl}). The focus of our paper is to construct a finite number of unitary representations in $\Pi_u(\mc O)$, for a class of  nilpotent orbits to which $\Pi_u(\mc O)$ is conjectured to be finite and nonempty (\cite{abv} \cite{vogan89}).

\section{Nilpotent orbits and unitary representations}\label{0.2} 
The nilpotent orbits we discussed in last section, are coadjoint nilpotent orbits. Recall that for semsimple Lie algebra, the Killing form provides an identification of $\f g_{\mathbb C}$ with its dual $\f g_{\mathbb C}^*$. In addition, this identification intertwins the adjoint action and the coadjoint action of the Lie group. Hence nilpotent coadjoint orbits can be identified with nilpotent (adjoint) orbits canonically. \\
\\
Recall that nilpotent orbits of $\f{sl}(n, \mathbb C)$ (or equivalently of $\f{gl}(n, \mathbb C)$) are in one-to-one correspondence with Jordan canonical forms. In particular, nilpotent orbits of $\f{gl}(n, \mathbb C)$ can be parametrized by partitions of $n$:
$$p_1+p_2 + \ldots +p_l=n$$
where $p_i \geq 1 \, (i \in [1, l])$ are the dimensions of the Jordan blocks.  Similarly, Nilpotent orbits  of $\f{sp}_{2n}(\mathbb C)$ can be  parametrized by
partitions of $2n$ such that odd parts occur with even multiplicities; nilpotent orbits of  $\f{o}(m, \mathbb C) $ can be  parametrized by partitions of $m$ such that even parts occur with even multiplicities. Throughout this paper, nilpotent orbits will always mean nilpotent coadjoint orbits and partitions will often be represented by a Young diagram with row lengths $(p_1, p_2, \ldots, p_l)$ in descending order. For example, the minimal nilpotent orbit of $Sp_{4}(\mb C)$ corresponds to the following Young diagram:
\begin{center}
 \fb \fb  \\
\fb  \mbo \\
\fb  \mbo \\
\end{center}
For a comprehensive account of nilpotent orbits, See \cite{cm}  and references therein.\\
\\
Generally speaking, a nilpotent orbit of a complex semisimple Lie algebra $\f g$  is said to be induced if it intersects
$\mc O_{\mathfrak l}+ \mathfrak n$ densely. Here $\mathfrak l \oplus \mathfrak n$ is a parabolic subalgebra of $\f g$, with $\mathfrak l$ the Levi factor, $\mathfrak n $ the nilradical, and $\mc O_{\mathfrak l}$ a nilpotent orbit of $\mathfrak l$. Such a nilpotent orbit is said to be induced  from $\mc O_{\mathfrak l}$. If a nilpotent orbit is induced from the zero orbit, it is often called a Richardson orbit. It is very easy to see that all nilpotent orbits 
of $\f{gl}(n, \mathbb C)$ are induced orbits, in fact Richardson orbits. However, there are  nilpotent orbits of $\f{sp}_{2n}( \mathbb C) ( n>1)$ and 
 of $\f o(m, \mathbb C) (m > 4)$ that are not induced. Nilpotent orbits that are not induced are called rigid nilpotent orbits, or simply rigid orbits. For example, the minimal orbit of $\f {sp}_{2n}(\mathbb C)$, consisting of rank 1 elements in $\f {sp}_{2n}(\mathbb C)$, is rigid.  \\
 \\
 There is a good characterization of rigid orbits for type B, C, D Lie algebras, in terms of the corresponding Young diagrams. Let $\mathfrak g$ be a complex simple Lie algebra of type B, C or D. Let $\mc O_{\bf d}$ be the nilpotent orbit parametrized by the Young diagram $\bf d$. Let $\bf d^t$ be the transpose. If $\bf d^t$ has any parts with multiplicity more than $1$, or equivalently, the Young diagram $\bf d$ has two columns of the same size, then $\mc O_{\bf d}$ is induced. For example,  the following  Young diagram
 \begin{center}
 \fb \fb \fb \fb \fb \fb  \\
\fb \fb \fb \fb \fb \fb  \\
\fb \fb \fb \mbo \mbo \mbo  \\
\fb \mbo \mbo \mbo \mbo \mbo  \\
\fb \mbo \mbo \mbo \mbo \mbo  \\
\end{center}
corresponds to a nilpotent orbit $\mc O_{\bf d}$ of $O(17, \mathbb C)$. Then $\mc O_{\bf d}=\Ind_{\f{gl}(3)\f{gl}(2) \f{o}(7)}^{\f{o}(17)} \mc O^{\prime} $ with $\mc O^{\prime}$ the nilpotent orbit of $GL(3, \mathbb C) \times GL(2, \mathbb C) \times O(7, \mathbb C)$ corresponding to the following three Young diagrams
\begin{center}
 \fb  \hspace{1 in}  \fb \hspace{1 in}  \fb \fb  \\
 \fb   \hspace{1 in}  \fb \hspace{1 in} \fb \fb  \\
 \fb  \hspace{1 in}  \mbo  \hspace{1 in} \fb \mbo  \\
 \mbo  \hspace{1 in}  \mbo \hspace{1 in} \fb \mbo  \\
 \mbo  \hspace{1 in}  \mbo \hspace{1 in} \fb \mbo  \\
\end{center}
We can see that if $\mc O_{\bf Y}$ is rigid then $\bf Y^t$ must be multiplicity free. \\
\\
Let us get back to our discussion about the relations between unitary dual and nilpotent orbits.  We fix a real semisimple Lie group $G$. For each irreducible unitary representation of $G$, its associated variety will be the closure of a complex nilpotent orbit in $\f{g}_{\mathbb C}^*$. If a complex nilpotent orbit $\mc O$ is induced from $\mc O_{\f l_{\mathbb C}}$ of a (complex) Levi algebra $\f l_{\mathbb C}$, there are techniques of parabolic induction and cohomological induction to construct an irreducible unitary representation $\pi$  such that $\mc V(Ann(\pi))$ is the closure of $\mc O$ (\cite{bv83} \cite{vogan86book} \cite{vogan86}). Of course these constructions will depend on whether there are unitary representations of a smaller group $L$ with associated variety $\mc O_{\f l_{\mathbb C}}$. We shall refer the readers to Vogan's works on orbital methods for the details of these constructions (\cite{vogan86book} \cite{vo94}). In any case, for rigid nilpotent orbits (which are not induced), these two constructions fail to produce the desired representations. The main purpose of this book is to produce unitary representations associated with rigid nilpotent orbits. Most of them shall be called unipotent representations since conjecturally, there are only finitely many such representations (\cite{vogan89}).\\
\\
There is yet one more subtlety. Not every complex rigid orbit can be the associated variety of an irreducible unitary representations of $G$. For example, the minimal nilpotent orbit of $\f{o}(2n+1, \mathbb C)$ is not the associated variety of any $SO(p, 2n+1-p)$-representation for $n \geq 3$ (\cite{hl}). 
\section{Theta Correspondence and Unipotent Representations}
One motivation for this book comes from the works of R. Howe, Jian-Shu Li, Przebinda and many others on theta correspondences. In 1980's, Howe and Li succeeded  in constructing a certain class of unitary representations, called the lower rank representations (\cite{howesmall} \cite{li891} \cite{li2}). Some of these lower rank representations are in fact unipotent (\cite{hl}). Roughly,   lower rank representations are those "close" to the trivial representation in the unitary dual, and their matrix coefficients tend to \lq\lq decay \rq\rq slowly.  However, in relationship to nilpotent orbits,  the associated varieties of lower rank representations correspond to certain small nilpotent orbits, consisting of nilpotent elements whose matrix rank is less than half of the size of the matrix (\cite{hl} \cite{he08}). So there are plenty of unitary representations that are not of lower rank. In any case, the main instrument behind the construction of lower rank representations is the local theta correspondence, also known as Howe's correspondence. Roughly, local theta correspondence, denoted by $\theta$, is a one-to-one correspondence between certain  infinitesimal equivalence classes of irreducible representations of the dual pair $(O(p,q), Sp_{2n}(\mathbb R))$ (\cite{howe}). For example, low rank representations of even $O(p,q)$ can be obtained from irreducible unitary representations of $Sp_{2n}(\mathbb R)$ with $n$ smaller than both $p$ and $q$. In addition, Przebinda  established that under certain conditions, the asymptotic cycle, a refinement of associated variety, behaves reasonably well under theta correspondence (\cite{pr1}). Przebinda's work suggested that theta correspondence yields a correspondence at the associated variety level under certain restrictions.  Therefore, one may be able to iterate theta correspondence to obtain unipotent representations inductively.  \\
\\
\commentout{Now we take the natural steps to generalize the constructions of Howe and Li to obtain a larger class of unipotent representations. For this purpose, we study the local theta correspondence $\theta$. For the sake of simplicity, we define $\theta(\pi)=0$ if $\pi$ does not appear in Howe's duality correspondence. Generalizing Li's definition of $(\, , \,)_{\pi}$ for stable range dual pairs, we define the semistable range for all dual pairs. We prove that if $\pi$ is in the semistable range, $\theta(\pi)$ can be constructed by an integration process (\cite{theta}). This paves the way to apply analytic techniques to study $\theta$. For simplicity we impose the condition that $\pi(\epsilon)=-1$  where $\{1, \epsilon\}$  is the preimage of the identity under the metaplectic covering. Within the semistable range, I obtain the following results.
\noindent}
In this paper, we carry out of the construction of unipotent representation using compositions of theta correspondences.  The analytic foundation of such process is laid in \cite{basic}, which   focused on the composition of two theta correspondences. Notice that when we compose two ortho-symplectic theta correspondences correctly, we obtain a correspondence of representations of two orthogonal groups, or two symplectic groups. we call this composition quantum induction. The reason behind this terminology is that the limit cases of quantum induction is closely related to the usual parabolic induction. The main focus of this book is  the quantum induction with respect to an alternating  chain of orthogonal groups and symplectic groups. We illustrate this process by the following example. \\
\\
Recall the (real) nilpotent orbits of the symplectic group is in one-to-one correspondences with certain signed Young diagrams (Ch 9. \cite{cm}). See Chapter \ref{nilorbits} for a review.
Let $\mc O_{\mathbf D}$ be a real nilpotent coadjoint orbit of $Sp_{36}(\mathbb R)$, with $\bold D=$
\begin{center}
\fm \fp \fm \fp \fm \fp  \\
\fm \fp \fm \fp \fm \mbo \\
\fp \fm \fp \fm \fp \mbo \\
\fm \fp \fm \fp \mbo \mbo  \\
\fm \fp \fm \fp \mbo \mbo  \\
\fm \fp \fm \mbo \mbo \mbo \\
\fp \fm \fp \mbo \mbo \mbo \\
\fm \fp \mbo \mbo \mbo \mbo \\
\fm \fp \mbo \mbo \mbo \mbo \\
\fm \mbo \mbo \mbo \mbo \mbo \\
\fp \mbo \mbo \mbo \mbo \mbo 
\end{center}
Hence the generic element, as a $36 \times 36$ matrix, will have the Jordan blocks of dimension 
$$(6,5,5,4,4,3,3,2,2,1,1).$$
Let $j$ be a positive integer. Define $\bold D-\bold j$ to be the Young diagram obtained by deleting the first $j$ columns from the left. Let $p_j$ be the number of $+$ in $\bold D-\bold j$ and $q_j$ be the number of $-$ in $\bold D-\bold j$. Construct the sequence
$$[(p_0,q_0),(p_1,q_1),(p_2,q_2), \ldots ],$$
obtaining
$$[(18,18), (15, 10), (8, 8), (6,3), (2,2), (1,0)].$$
Let $\chi$ be a unitary character of $O(1)$. Let $Mp_{2n}(\mb R)$ be the unique double cover of $Sp_{2n}(\mb R)$.
Define 
\begin{equation}
\begin{split}
& Q(36;15,10;16;6,3;4;1,0)(\chi) \\
= & \theta_s(O(15,10), Mp_{36}(\mathbb R)) \theta_s(Mp_{16}(\mathbb R), O(15,10)) \theta_s(O(6,3), Mp_{16}(\mathbb R)) \\
& \theta_s(Mp_{4}(\mathbb R), O(6,3)) \theta_s(O(1,0), Mp_{4}(\mathbb R)) (\chi). \\
\end{split}
\end{equation}
where $\theta_s$ is the theta correspondence in the semistable range (see Definitions \!~\ref{thetas} and \!~\ref{quantum}). As proved in \!~\cite{basic}, $Q(36;15,10;16;6,3;4;1,0)(\chi)$  is a genuine irreducible unitary representation of $Mp_{36}(\mathbb R)$ if it does not vanish 
(see Theorem \!~\ref{quantuminduction1} and Theorem \!~\ref{quantuminduction2}).  The nonvanishing of $Q(*)(\chi)$ is proved in  Chapter \ref{nonvan} (Theorems \!~\ref{non2}, also \!~\ref{exi}, \!~\ref{non3}). The details of our construction are then given in Chapter \ref{construct}.  
In Chapter 5 an 6, we compute the infinitesimal character and the associated variety of $Q(*)(\chi)$ based on the results obtained by Przebinda (\!~\cite{pr} \!~\cite{pr1}). We further prove that
the associated variety of $Q(*)(\chi)$ is indeed the closure of the complexified orbit of $\mc O_{\mathbf D}$. The resulting representation
$Q(*)(\chi)$ have the right associated variety and infinitesimal character, conjectured for  unipotent representations (\cite{arthur} \cite{abv} \cite{vogan89}).
\section{Main Results}
We shall now discuss how general our approach is in terms of constructing unipotent representations.
As we mentioned earlier, not every complex nilpotent orbit can be the associated variety of an irreducible unitary representation. As far as we know, there is no complete characterization of associated varieties of irreducible unitary representations for orthogonal groups or symplectic groups except for some special cases. Nevertheless, there is a fairly good characterization of associated varieties of irreducible representations with integral infinitesimal characters, namely, those infinitesimal character coming  from  finite dimensional representations.  The main result, following the line of studies on primitive ideals, states that, the associated variety of an irreducible representation with integral infinitesimal character, is a special orbit in the sense of Lusztig (\cite{lus}). We shall refer the readers to \cite{du} \cite{jo80}  \cite{bv82} \cite{jan}  and references therein for historical account of primitive ideals and associated varieties. We shall also mention that special orbits for symplectic and orthogonal groups have a simpler characterization in terms of the transpose of Young diagram associated with the nilpotent orbits (Prop 6.3.7 \cite{cm}). For Lie algebras of type $C$ and $D$, $\mc O_{\bold d}$ is special if any odd number in the transpose $\bold d^t$ occurs with even multiplicities. \\
\\
We shall remind the reader that the description above pertains to Lie algebra representations and {\it complex} nilpotent orbits.  As to representations of  a fixed real form $G$, there are special complex nilpotent orbits  not associated to any representations. A necessary condition for this to occur, is that there is no real nilpotent orbit of $G$, that complexifies to such special nilpotent orbit. Indeed, by the results in this paper, it is sufficient. Here is one highlight of our paper. \commentout{ stated not in full generality, but sufficient for the purpose of an introduction.}
\commentout{In order to address this issue, we call a real nilpotent orbit, special if the corresponding complex orbit is special; rigid, if the corresponding complex orbit is rigid.}  
\begin{thm}~\label{main1}
Let $G=Sp_{2n}(\mathbb R)$ or $G=O(p,q)$ with $p+q$ even.
Let $\mc O_{\bold d}$ be a special rigid nilpotent  orbit of $G_{\mathbb C}$ parametrized by the partition $\bold d$ (see Ch 6.3,7.3 \!~\cite{cm}). Let $\bold d^t=(m_1> m_2 > m_3 > \ldots > \ldots)$ be the transpose of $\bold d$ (considered as a Young diagram). Then all $m_i$ will be even. Let $\mc O$ be a real nilpotent orbit of $G$ in $\mc O_{\bold d}$.
Then there exists a nonempty set $\mc N(\mc O)$ of irreducible unitary representations of $G$  such that for every $\pi \in \mc N(\mc O)$
the associated variety of $\pi$,
$\mc V(Ann\,\pi)=cl(\mc O_{\bold d})$.
 the infinitesimal character
$\mc I(\pi)$ only depends on $\bold d$.
\begin{itemize}
\item If $G=Sp_{2n}(\mathbb R)$, $\mc I(\pi)=(\rho(\f{sp}_{m_1}(\mathbb C)), \rho(\f o(m_2, \mathbb C)), \rho(\f{sp}_{m_3}(\mathbb C)), \rho(\f o(m_4, \mathbb C)), \ldots)$;
\item If $G=O(p,q)$, $\mc I(\pi)=(\rho(\f o(m_1, \mathbb C)), \rho(\f{sp}_{m_2}(\mathbb C)), \rho(\f o(m_3, \mathbb C)), \rho(\f{sp}_{m_4}(\mathbb C)), \ldots)$.
\end{itemize}
Here $\rho(\f g)$ is the half sum of the positive roots of $\f g$.
 \end{thm}
 
 In short, we obtain unitary representations attached to special rigid orbits. 
 Let me make a few remarks here.\\
 \\
{\bf Remark.} 
\begin{enumerate} 
\item From Ch. \ref{0.2}, a nilpotent orbit  $\mc O_{\bold d}$ is rigid implies that $\bold d^t$ is multiplicity-free.  For Lie algebras of type $C$ and $D$, $\mc O_{\bold d}$ is special if every odd number in $\bold d^t$ occurs with even multiplicities (see 
\cite{lus} also Proposition 6.3.7 \cite{cm}).  Therefore, $\mathcal O_{\bold d}$ being rigid and special implies that no odd number can appear in $\bold d^t$. Hence  $m_i$ must be all even and multiplicity free. Similarly, for type B simple Lie algebras, a nilpotent orbit $\mc O_{\bf d}$ being rigid and special implies that $\bf d^t$ can only have odd parts and must be multiplicity free.  We can construct unitary representations of odd orthogonal groups attached to special rigid orbits in the same fashion. Same statements hold for odd $O(p,q)$. See Theorems \!~\ref{exi}, \!~\ref{un}, \!~\ref{un}, \!~\ref{infi}, \!~\ref{as} and Definition \!~\ref{u} for details. 

\item Our construction of $\mc N(\mc O)$ is based on the real nilpotent orbit $\mc O$ in $\mc O_{\bold d} \subseteq  \mathfrak g_{\mathbb C}$. A complex nilpotent orbit $\mc O_{\bold d}$ can contain many different real nilpotent orbits $\bold O_{\bf D}$ of the same real group $G$. Essentially, this amounts to putting signs $\pm$ in the Young diagram $\bold d$ and making signed Young diagram $\bold D$. Then all
$\mc N(\mc O_{\bold D})$ will be attached to $\mc O_{\bold d}$. 
\item Some of the representations constructed in this paper has been studied before. For example, a  class of representations,  called  minimal representations, the representations attached to minimal orbits, were studied by Brylinski-Kostant 
 in a series of papers (\!~\cite{bk}) and by Binegar-Zierau for $O(p,q)$ (\!~\cite{bz}) from a different angle. Beyond that,
representations attached to a perhaps wider class of small rigid orbits have also been studied by 
Kashiwara-Vergne, Howe, Li, Sahi, Tan, Huang, Zhu and others (see \!~\cite{kave}, \!~\cite{howesmall}, \!~\cite{sahi},  \!~\cite{ht}, \!~\cite{hz}, \!~\cite{hl} and the references within them). 
\commentout{\item The intrinsic connection between coadjoint orbits and the unitary dual of a Lie group was first explored by Kirillov and Kostant. It is now known as the {\it orbit method}.
One core part of the orbit method for real reductive Lie groups is to construct unipotent representations attached to rigid nilpotent orbits. Once we know how to attach representations to rigid nilpotent orbits, there are various ways to attach unitary representations to induced orbits (see \!~\cite{vo94}, \!~\cite{vogan86book}). What we have accomplished in this paper is the construction of
some unipotent representations attached to special rigid orbits and some other nonspecial rigid nilpotent orbits.
It is not clear whether our list of $\mc N(\mc O)$ exhausts all unipotent representations attached to these rigid orbits. At least for some special rigid orbits, 
it does (\!~\cite{hl}). We conjecture that $\mc N(\mc O)$ is exhaustive for special rigid orbits.}
\item For the metaplectic group $Mp_{2n}(\mathbb R)$, similar statements hold. See Theorems \!~\ref{exi}, \!~\ref{un}, \!~\ref{un}, \!~\ref{infi}, \!~\ref{as} and Definition \!~\ref{u} for details. The representations constructed for $Mp_{2n}(\mathbb R)$ will have half integral infinitesimal characters and no longer correspond to special orbits. As in the earlier example for $Mp_{36}(\mathbb R)$, the nilpotent orbit $\mc O_{\bf d}$ is not special, but rigid.  The  infinitesimal character of the constructed representation  
$$[\frac{3}{2}, \frac{1}{2}, \frac{3}{2}, \frac{1}{2}, \frac{7}{2}, \frac{5}{2}, \frac{3}{2},\frac{1}{2}, \frac{7}{2}, \frac{5}{2}, \frac{3}{2},\frac{1}{2}, \frac{11}{2}, \frac{9}{2}, \frac{7}{2}, \frac{5}{2}, \frac{3}{2},\frac{1}{2}].$$
 is not integral. This can be proved using a Theorem of Przebinda (\cite{pr}). Hence
our construction of $\mc N(\mc O)$ is effective for a large number of nonspecial rigid nilpotent orbits as well.
\item One of the most well-known nonspecial rigid orbit is the minimal orbit of $Sp_{2n}(\mb C)$, namely, the unique rank $1$ nilpotent orbit in $\f{sp}_{2n}(\mb C)$. The corresponding irreducible representations are  the two irreducible constituents of the oscillator representation (also called the Segal-Shale-Weil representation) and their two dual representations. They are perhaps the most well-known  unipotent representations attached  to a non-special rigid orbit. 
 The representations constructed in this paper can be regarded as derivations of the oscillator representation.

\end{enumerate}

 Finally, recall that every special nilpotent orbit is induced from a special rigid orbit. This amounts to reducing the multiplicities in $\bf d^t$ to $1$.  Unitary representations associated with induced orbits can be constructed by unitary parabolic induction (see \!~\cite{bv}, \!~\cite{vogan01}). Hence we can construct irreducible unitary representations of $Sp_{2n}(\mb R)$ or even orthogonal groups $O(p,q)$ attached to any special nilpotent orbits.
 \begin{thm}
 Let $\mc O_{\bold d}$ be a special nilpotent orbit of $Sp_{2n}(\mathbb C)$ or $O(p+q, \mathbb C)$. Suppose that  $\mc O_{\bold D}$ is a nilpotent orbit of a real group $G=Sp_{2n}(\mb R)$ or $O(p,q)$  contained in $\mc O_{\bold d}$. Then  there exists an irreducible unitary representation $\pi$ of $G$ such that
$\mc V(Ann \pi)= cl(\mc O_{\bold d})$.
\end{thm}
Also included in this paper is a theorem concerning the relationship between quantum induction and parabolic induction. This theorem will provide us a way to prove the nonvanishing theorem of quantum induction $Q(p,q)(\pi)$. 
Based on a theorem of Kudla-Rallis \!~\cite{kr} and of Lee-Zhu \!~\cite{lz}, we decompose a certain parabolically induced representation $I^{\alpha}(\pi)$  of $Mp_{2n}(\mathbb R)$ into a direct sum of quantum induced representations $Q(p,q)(\pi)$ with $p+q=n+1$ and a fixed parity on $p$ 
(see Definition \!~\ref{ipi}, Theorem \!~\ref{induction} and Theorem \!~\ref{induction1}). Each $Q(p,q)(\pi)$ is either irreducible or vanishes. This is the limit case for which quantum induction can be constructed from parabolic induction.
Correspondingly, there is a decomposition theorem for
the wave front set of $I^{\alpha}(\pi)$, namely
$$WF(I^{\alpha}(\pi))= \bigcup_{p+q=n+1, p-q \equiv \alpha \pmod 4} WF(Q(p,q)(\pi)) .$$
On the one hand, $WF(I^{\alpha}(\pi))$ is computable and consists of a finite
number of irreducible components. On the other hand, the wave front set of each $Q(p,q)(\pi)$ possesses a certain distinctive characteristic that can be derived from \!~\cite{pr1} and \!~\cite{pan}. Under certain hypotheses on $\pi$, we sort out the occurrence of $WF(Q(p,q)(\pi))$ in $WF(I^{\alpha}(\pi))$ completely. This provides us with one nonvanishing theorem  for
$Q(p,q)(\pi)$ and for $\theta_s(p,q)(\pi)$.\\
\\
Let me point out one advantage of our construction related to the theory of automorphic forms. Theta correspondence originated from the studies of theta series (see \!~\cite{si}, \!~\cite{weil}, \!~\cite{weilsiegel}, \!~\cite{howe79}). As pointed out to me by Jian-Shu Li, theta correspondence should map automorphic representations to automorphic representations (see \!~\cite{howe79}, \!~\cite{ra87}, \!~\cite{li94}). Thus representations in $\mc N(\mc O)$ should all be automorphic. In the final part of this paper, along the lines of \!~\cite{vogan86}, we make some conjectures regarding the automorphic dual in the sense of Burger-Li-Sarnak \!~\cite{bls}. \\
\\
The first draft of this paper was finished in early 2002. Since then, there have been many interesting developments related to this paper, notably \cite{pt} \cite{bp} \cite{ar13} \cite{he15} \cite{bmsz}.
  I wish to acknowledge my gratitude to Monica Nevins, Shu-Yen Pan, Tomasz Przebinda and 
David Vogan for some very helpful e-mail communications and to Ray Fabec, Bill Graham,  Gestur Olafsson and Chen-Bo Zhu for reading the manuscript and for their comments. I also want to thank the referees for their very valuable suggestions and comments, in particular their suggestions concerning the structure of the introduction, Definition \!~\ref{unip}, the proof of Theorem\!~\ref{induction1} and the structure of Chapter \ref{ch7}. Special thanks go to the former AMS Memoirs editor Daniel Bump and the current AMS Memoirs editor Henri Darmon for their efforts to get this paper published.
In the midst of revising this paper, my father Decai He,  passed away. My father had always encouraged and inspired me in my mathematical endeavors. I would like to dedicate this paper to him. 

\tableofcontents

\chapter{Invariants}
In mathematics, classification problems are often approached by constructing invariants.
In this chapter, we attach invariants to the equivalence classes of irreducible admissible Hilbert representations of a reductive group $G$. One hopes that these invariants could shed some light on the classification and construction of irreducible unitary representations. Our purpose here is not to give a historical account of these invariants, but rather, to review some basic facts we need concerning these invariants. We adopt the standard notations from \!~\cite{knapp} and \!~\cite{wallach}.

\section{Notation}
Reductive groups and semisimple groups in this paper are assumed to have at most a finite number of components. Let $G$ be a Lie group. We adopt the following notation:
\begin{enumerate}
\item $G_0$---the identity component of $G$;
\item $\f g$---the real Lie algebra of $G_0$;
\item $\f g^*$---the dual space of $\f g$;
\item $\f g_{\mathbb C}$---the complex Lie algebra of $G_0$;
\item $U(\f g)$---the complex universal enveloping algebra;
\item $Z(\f g)$---the center of $U(\f g)$;
\end{enumerate} 
Let $G$ be a real reductive group. We adopt the following notation:
\begin{enumerate}
\item $K$---a maximal compact subgroup of $G$;
\item $\Pi(G)$---the set of equivalence classes of irreducible $(\f g, K)$-modules;
\item $\Pi_u(G)$---the set of equivalence classes of unitarizable $(\f g, K)$- modules, or equivalently, the set of equivalence classes of irreducible unitary representations of $G$.
\item $(\ , \ )$---an invariant quadratic form on $\f g$ such that $(\ , \ )|_{\f k}$ is positive definite;
\item $\f p$---the orthogonal complement of $\f k$ with respect to $(\ , \ )$;
\item $\f a$---a maximal Abelian Lie subalgebra of $\f p$;
\item $A$---the connected Abelian group generated infinitesimally by $\f a$;
\item $KAK$---a Cartan decomposition.
\end{enumerate}
Let $V$ be a finite dimensional vector space over $\mathbb F$. We use $V^*$ to denote 
$\Hom_{\mathbb F}(V, \mathbb F)$. 
Let $G$ be a real reductive group. 
\begin{no}
Let $(\pi, \mc H)$ be a Hilbert representation of $G$. Fix a maximal compact subgroup $K$. Let $V( \pi)$ be the space of smooth $K$-finite vectors in $\mc H$. To emphasize $\pi$, $\mc H$ is often denoted by $\mc H_{\pi}$.
\end{no}
The space $V(\pi)$ is a $(\f g, K)$-module. $\pi$ is said to be admissible if each $K$-type occurs in $V(\pi)$ with finite multiplicity.  If, in addition, $V(\pi)$ is finitely generated as a $U(\f g)$-module, then $V(\pi)$ is often called a Harish-Chandra module.  \\
\\
If $\pi$ is unitary and irreducible, then
the Hilbert norm is unique up to a scalar multiplication. So, equivalence classes of irreducible unitary representations are in one-to-one correspondence with irreducible unitarizable  
Harish-Chandra modules.\\
\\ 
{\bf All unitary representations in this paper, unless otherwise stated, are taken as suitable
unitarized Harish-Chandra modules}. \\
\\
This convention is consistent with the way we construct irreducible unitary representations. First, we will construct a
$(\f g, K)$-module $(\pi, V)$. Then we will show that $V$ is an irreducible $(\f g, K)$-module. Hence $V$ is a Harish-Chandra module. Finally, we will prove the existence of an invariant inner product on $V$. Then $\pi$ is an irreducible unitary representation of $G$.
For simplicity, we use $\pi$ to denote the group action and the
Lie algebra action. We define three involutions in the category of Harish-Chandra modules:
\begin{enumerate}
\item $\pi^*$, the contragredient representation;
\item $\pi^c$, the representation $\pi$ equipped with the conjugate complex linear structure;
\item $\pi^h$, the Hermitian dual representation of $\pi$. 
\end{enumerate}
We have $(\pi^*)^c=\pi^h$.
If $\pi$ is unitary, then $\pi^h \cong \pi$. 
\begin{no}~\label{preceq}
Suppose $a, b \in \mathbb R^n$.
Write $a \preceq b$ if and only if
for every $1 \leq k \leq n$,
$$\sum_{j=1}^k a_j \leq \sum_{j=1}^k b_j;$$
write $a \prec b $ if and only if for every $1 \leq k \leq n$,
$$\sum_{j=1}^k a_j < \sum_{j=1}^k b_j.$$
\end{no}
The ordering $\preceq$ is a partial ordering.

\section{Infinitesimal Character and Harish-Chandra Homomorphism}
Let $\f g$ be a complex reductive Lie algebra. 
Let $\f h$ be a Cartan subalgebra of $\f g$. Let $W(\f g, \f h)$ be the Weyl group generated by the root system $\Sigma(\f g, \f h)$. Let $U(\f h)^{W(\f g , \f h)}$ be the space of $W(\f g, \f h)$-invariant vectors in $U(\f h)$.
Then the Harish-Chandra homomorphism 
$${\mc{HC}}: Z(\f g) \rightarrow U(\f h)^{W(\f g, \f h)}$$
is an algebra isomorphism (see \!~\cite{knapp}, Chapter VIII. 5 or \!~\cite{wallach}).
Identify $ U(\f h)$ with the symmetric algebra of 
$\f h$. For each vector $\Lambda$ in the complex dual space of $\f h$, define a character $\chi_{\Lambda}$ of $Z(\f g)$ by
$$\chi_{\Lambda}(x)=\Lambda({\mc{HC}}(x))={\mc{HC}}(x)(\Lambda) \qquad (x \in Z(\f g))$$
It is well-known that every character of $Z(\f g)$ can be obtained this way and that $\Lambda$ is unique up to the action of $W(\f g, \f h)$. In short, $Spec(Z(\f g)) \cong \f h^*// W(\f g, \f h)$. Here the categorical quotient ${\f h}^*//W(\f g, \f h)$ coincides with the geometric quotient since $W(\f g, \f h)$ is finite.\\
\\
Let $G$ be a connected real reductive Lie group with Lie algebra $\f g$.
Let $K$ be a maximal compact subgroup of $G$. 
Let $(\pi, \mc H)$ be an irreducible admissible Hilbert representation of $G$. Let $V({\pi})$ be the Harish-Chandra module of $(\pi, \mathcal H)$, consisting of all the $K$-finite vectors in the Hilbert space $\mc H$. We retain $\pi$ for the infinitesimal action of $U(\f g)$ on the smooth vectors of $(\pi, \mc H)$.  Since $\pi$ is irreducible, $V({\pi})$ is an irreducible $(U(\f g),K)$-module. Since $G$ is connected, $Z(\f g)=U(\f g)^G$.
By Schur's lemma, $Z(\f g)$ acts on $V({\pi})$ by a character
$$\chi: Z(\f g) \rightarrow \mathbb C.$$
Thus there exists a $\Lambda$ such that $Z(\f g)$ acts on $V(\pi)$ by $\chi_{\Lambda}$. For simplicity, we call
$\Lambda$ an {\bf infinitesimal character} of $\pi$.\\
\\
Let $G$ be a real reductive group with a finite number of components. Let $\f h$ be a complex Cartan subalgebra in $\f g_{\mathbb C}$.  Then $U(\f g)^{G_0}=Z(\f g)$. Let $U(\f g)^G$ be the $G$-invariant vectors in $U(\f g)$. Consider the adjoint action of $G/G_0$ on $Z(\f g)$. Each element in $G/G_0$ acts on $Z(\f g)$ as an algebra automorphism.  Furthermore, $U(\f g)^G$ is precisely the subalgebra of $Z(\f g)$ invariant under the action of $G/G_0$. By the Harish-Chandra homomorphism,  $G/G_0$ acts on 
$$ U(\f h)^{W(\f g_{\mathbb C}, \f h)}$$
as algebra automorphisms. This action induces an action of $G/G_0$ on 
$$\spec(U(\f h)^{W(\f g_{\mathbb C}, \f h)}).$$
Since $G/G_0$ is finite, by invariant theory
$$\spec({\mc{HC}}(U(\f g)^G))$$ is precisely in one-to-one correspondence with the $G/G_0$-orbits
in
$$ \spec(U(\f h)^{W(\f g_{\mathbb C}, \f h)}).$$ 
\\
Let $\pi$ be an irreducible $(\f g, K)$-module. By Schur's lemma, $U(\f g)^G$ must act by a character $\xi$. 
\begin{defn}
Let $\f h$ be a Cartan subalgebra of $\f g_{\mathbb C}$. We say that $\Lambda \in \f h^*$ is an infinitesimal character of $\pi$ if $\chi_{\Lambda}|_{U(\f g)^G} = \xi$. 
\end{defn}
$\Lambda$ is unique up to the action of $G/G_0$ and $W(\f g_{\mathbb C}, \f h)$. 
For semisimple Lie group $G$, the action of $G/G_0$ on $\f h^*// W(\f g_{\mathbb C}, \f h)$ is not difficult to understand.
Let $\iota$ be the action of $G/G_0$ on $\f h^*// W(\f g_{\mathbb C}, \f h)$, induced from the action of $G/G_0$ on $\f g^* // G_0$. Then $\iota$ induces an action of $G/G_0$
on a closed Weyl chamber, consequently on the positive root system. Hence $\iota(G/G_0)$ can be regarded as
 automorphisms of the Dynkin diagram of $\f g_{\mathbb C}$. The automorphisms of the Dynkin diagram are not difficult to classify. \\
\\
The infinitesimal character is one of the main tools used in the literature to study irreducible admissible representations. In fact, there is only a finite number of infinitesimal equivalence classes of irreducible representations with a fixed infinitesimal character. 
\begin{no}
Let $\Pi_{\Lambda}(G)$ be the set of infinitesimal equivalence classes of irreducible admissible Hilbert representations of $G$ with infinitesimal character
$\Lambda$.
\end{no}

\section{Leading Exponents of Irreducible Representations}
Let $G$ be a real reductive Lie group. Fix a maximal compact subgroup $K$.
Unless otherwise stated, matrix coefficients in this paper are assumed to be $K$-finite. 
Fix a nondegenerate real invariant bilinear form $(\ , \ )$ on $\f g$ and a maximal
Abelian Lie subalgebra $\f a$ of $\f p$ as in (1.1). Let $r$ be the real dimension of $\f a$. We call $r$ the {\it real rank } of $G$. Let $\Sigma(\f g, \f a)$ be the restricted root system. \\
\\
Fix a positive root system $\Sigma^+(\f g, \f a)$. Let $\rho(G)$ be the half sum of all positive roots in $\Sigma(\f g, \f a)$. Let $M$ be the centralizer of $\f a$ in $K$. Let $W(G, \f a)$ be the normalizer of $\f a$ in $K$ modulo $M$.  We call $W(G, \f a)$ the real Weyl group. Let $W(\f g, \f a)$ be the Weyl group generated by the root system. Clearly, $W(\f g, \f a) \subseteq W(G, \f a)$.  \\
\\
For $Sp_{2n}(\mathbb R)$, $\f a$ is isomorphic to $\mathbb R^n$.
The real Weyl group $W(Sp_{2n}(\mathbb R), \f a)$ is generated by permutations and sign changes on $n$ variables. For $O(p,q)$, $\f a$ is isomorphic to
$\mathbb R^{\min\{p,q\}}$. The real Weyl group $W(O(p,q), \f a)$ is also generated
by permutations and sign changes on $\min\{p,q\}$ variables. See 1.1.1 and 1.3.2. \\
\\
Attached to each irreducible admissible representation $\pi$ of a connected $G$, is a finite
number of vectors in $\f a_{\mathbb C}^*$, called leading exponents. Leading exponents are
 the main data used to produce the Langlands classification (see \!~\cite{langlands}, \!~\cite{knapp}, \!~\cite{wallach}). For connected $G$, leading exponents depend on the choice of $\Sigma^+(\f g, \f a)$. For $G$ with finite number of connected components, we define the leading 
exponents of $\pi \in \Pi(G)$ to be the leading exponents of the irreducible subrepresentations of $\pi|_{G_0}$. \\
\\
Leading exponents are closely related to the infinitesimal character. For an irreducible finite dimensional representation $\pi$ of a connected real reductive group $G$ with real rank equal to complex rank, leading exponent is just the highest weight of $\pi$ with respect to a maximally split Cartan subgroup. In this situation,
finite dimensional representation theory says that the  highest weight $v$ is related to the infinitesimal character
$\Lambda$ by the following equation
$$v+\rho= w \Lambda$$
for some $w$ in $W(G, \f h)$. A similar statement holds for leading exponents for any irreducible admissible representation of a real reductive group.
By Theorem 8.33 from \!~\cite{knapp}, we have
\begin{thm}~\label{mat}
Let $\f b$ be a Cartan subalgebra of $\f m$. Take 
$\f h= \f a \oplus \f b.$
Suppose $\pi$ is an irreducible admissible representation of a real reductive group $G$. Let $v$ be a leading exponent of $\pi$. Then there exists an infinitesimal character $\Lambda$ of $\pi$ in $\f h_{\mathbb C}^*$ such that
$$v+ \rho(G)=\Lambda|_{\f a}.$$
\end{thm}
Notice that for non-split groups, $v$ is in $\f a_{\mathbb C}^*$ and $\Lambda$ is in
$(\f a \oplus \f b)_{\mathbb C}^*$. By Theorem \!~\ref{mat}, for any $\pi \in \Pi_{\Lambda}$, the set of leading exponents is finite.
\begin{no}
For $v \in \f a_{\mathbb C}^*$, the complex dual of $\f a$, we denote the real part of $v$ by $\Re (v)$.
\end{no}
Leading exponents are extracted from a certain asymptotic expansion of the matrix coefficients at
$\infty$. Therefore, they control the growth of matrix coefficients.
We cite the following estimate (Theorem 8.47 \!~\cite{knapp}).
\begin{thm}~\label{kn}
Let $\pi$ be an irreducible admissible representation of a real reductive group $G$. 
Let $v_0$ be in the real dual of $\f a$. Let $\f a^+$ be the closed Weyl Chamber associated with 
$ \Sigma^+(\f g, \f a)$ (see Prop. 5.14 \!~\cite{knapp}). Let $A^+$ be the closed Weyl Chamber in $A$ obtained by exponentiating $\f a^+$. Let $a(g)$ be the middle term of the $KA^+K$ decomposition of $g$. 
If every leading exponent $v$ of $\pi$ satisfies
\begin{equation}~\label{v0}
(v_0-\Re(v))(H) \geq 0 \qquad (\forall \ \ H \in \f a^+),
\end{equation}
then there is an integer $q \geq 0$ such that each $K$-finite matrix coefficient of $\pi$ is dominated  by a multiple of
$$\exp (v_0 (\log a(g))) (1 + (\log a(g), \log a(g))^{\frac{1}{2}})^q.$$
The converse also holds.
\end{thm}
When $G=O(p,p)$, due to the action of $G/G_0$, $A^+$ and $\f a^+$ can be made smaller (see 1.3.1). 
Theorem \!~\ref{kn} remains valid for this generalized Weyl chamber $A^+$. \\
\\
Theorem \!~\ref{kn} provides a uniform bound for matrix coefficients for all
$\pi \in \Pi_{\Lambda}(G)$.
Let $\pi$ be a unitary representation in $\Pi_{\Lambda}(G)$.
If $\Lambda$ is \lq\lq small \!\rq\rq, we can easily find a $v_0$ to dominate the set
$$\{ w \Lambda|_{\f a}- \rho(G) \mid w \in W(G, \f h_{\mathbb C}) \}.$$ 
By Theorem \!~\ref{kn}, we gain fairly good control over the growth of the matrix coefficients for all $\pi \in \Pi_{\Lambda}$. If $\Lambda$ is large, then the convex cone spanned by the set 
$$\{ w \Lambda|_{\f a}- \rho(G) \mid w \in W(G, \f h_{\mathbb C}) \}$$ 
is widely spread. The uniform bound from Theorem \!~\ref{kn} does not yield any useful  information. \\
\\
For unitary representations in $\Pi(G)$, the matrix coefficients are all bounded by constant functions. Theorem \!~\ref{kn} then implies
\begin{cor}
Let $G$ be a real reductive group.
If $\pi \in \Pi_u(G)$, then every leading exponent $v$ of $\pi$ satisfies
$$ \Re(v) (H) \leq 0 \qquad (\forall \ H \in \f a^+).$$
\end{cor}
In this paper, we will mostly deal with small
$\Lambda$. To give a broader picture of the importance of $\Pi_u(G)$ with small infinitesimal character, we recall one definition from \!~\cite{vo00} and \!~\cite{sv}. 
\begin{defn} A unitary representation $\pi \in \Pi_{\Lambda}(G)$
is  called \lq\lq unitarily small \rq\rq if $\Lambda$ is in the convex hull spanned by
$$\{ w (\rho(\f g_{\mathbb C})) \mid w \in W(\f g, \f h) \}.$$
\end{defn}
Salamanca-Vogan conjectured that any irreducible unitary representation of $G$ not small may be constructed by parabolic or cohomological induction from
an irreducible unitarily small representation of a reductive subgroup of $G$ (\!~\cite{sv}, \!~\cite{vo00}).\\
\\
It is precisely for \lq\lq unitarily small \!\rq\rq representations that Theorem \!~\ref{kn} will
produce useful information about matrix coefficients beyond what is stated in
the corollary. One goal of this paper is to construct and study some unitarily small representations
far away from the tempered unitary representations.
\subsection{Example I: The groups $Mp_{2n}(\mathbb R)$ and $Sp_{2n}(\mathbb R)$}
Let $G$ be either $Mp_{2n}(\mathbb R)$ or $Sp_{2n}(\mathbb R)$. Fix
$$\f a=\{\diag(a_1,a_2, \ldots a_n, -a_1, -a_2, \ldots -a_n) \}.$$
Then the Weyl group $W(G, \f a)$ is generated by permutations and sign changes on $\{a_i\}_{i=1}^n$.
Fix positive roots $\Sigma^+=\{ e_i \pm e_j \mid i \leq j; i,j \in [1,n] \}$. Then
$$\rho(G)=(n,n-1,\ldots, 1)$$
and 
$$\f a^+=\{ \diag(a_1, a_2, \ldots, a_n, -a_1, -a_2, \ldots, -a_n) \mid a_1 \geq a_2 \geq \ldots \geq a_n \geq 0 \}.$$
Clearly,
$$A^+=\{\diag(A_1, A_2, \ldots A_n, A_1^{-1}, A_2^{-1}, \ldots A_n^{-1}) \mid
A_1 \geq A_2 \geq \ldots \geq A_n \geq 1 \}.$$
Finally, the Inequality (\!~\ref{v0}) is equivalent  to
$$ v_0 \succeq \Re (v)$$
(see Notation \!~\ref{preceq} and 4).
\subsection{Example II: The Groups $O(p,q)$}
Let $G=O(p,q)$ with $q \geq p$. Here $O(p,q)$ is the isometry group of
$$(x, y)=\sum_{i=1}^p ( x_i y_{p+i}+ x_{p+i} y_i) + \sum_{j=1}^{q-p} x_{2p+i} y_{2p+i} .$$
Fix
$$\f a=\{\diag(a_1, a_2, \ldots a_p, -a_1, -a_2, \ldots - a_p, 0, 0, \ldots 0 ) \mid a_i \in \mathbb R \}.$$
Then the Weyl group $W(G, \f a)$ is again generated by permutations and sign changes on
$\{a_i\}_{i=1}^p$. Notice for $p=q$, $W(G, \f a)$ is bigger than $W(\f g, \f a)$. \\
\\
Fix restricted positive roots
$$\Sigma^+=\{ e_i \pm e_j \mid i < j; i,j \in [1,p] \} \cup \{ e_i \mid i \in [1,p]\} \ \  \mbox{if }  p < q$$
and $$\Sigma^+=\{ e_i \pm e_j \mid i < j; i,j \in [1,p] \} \ \  \mbox{if } p=q .$$
Then
$$\rho(G)=(\overbrace{\frac{p+q-2}{2}, \frac{p+q-4}{2}, \ldots \frac{q-p}{2}}^p).$$
Fix a (generalized) Weyl chamber 
$$\f a^+=\{\diag(a_1, a_2, \ldots a_p, -a_1, -a_2, \ldots - a_p, 0, 0, \ldots 0 ) \mid a_1 \geq a_2 \geq \ldots a_n \geq 0 \}.$$
Then
$$A^+=\{ \diag(A_1, A_2, \ldots A_p, A_1^{-1}, A_2^{-1}, \ldots, A_p^{-1}) \mid A_1 \geq A_2 \geq \ldots A_p \geq 1 \}.$$
The Inequality (\!~\ref{v0}) is again equivalent to $v_0 \succeq \Re(v)$. Notice that
for $O(p,p)$, due to the action of
$O(p,p)/SO_0(p,p)$, $A^+$ is half of the Weyl Chamber determined by $\Sigma^+$.
\section{Global Characters}
The global character is also known as the Harish-Chandra character.
For each irreducible admissible Hilbert representation $\pi$ of $G$ and a compactly supported smooth function $f$ on $G$, define $\Theta(\pi)(f)$ to be the trace of the operator
$$\pi(f)=\int f(g) \pi(g) d g.$$
 A theorem of Harish-Chandra states that 
$\Theta_{\pi}$ is a distribution on $G$ which can be identified with a locally integrable function (still denoted by $\Theta_{\pi}$) and $\Theta_{\pi}$ is real analytic on the set of regular semisimple elements of $G$ (see \!~\cite{hc}).\\
\\
For $\pi$ unitary, Mili\v{c}i\'{c} defined a notion of the rate of growth of $\Theta_{\pi}$ which we will not recall here. We will state one theorem relating 
$\gamma$ to the matrix coefficients of $\pi$ (see \!~\cite{mi}). 
\begin{thm}[Thm. 1 \!~\cite{mi}]~\label{mi}
Let $\pi \in \Pi_u(G)$. Let $\Xi(g)$ be Harish-Chandra's $\Xi$ function. 
Then the following are equivalent:
\begin{enumerate}
\item $\Theta_{\pi}$ has the rate of growth $\gamma \in \mathbb R$;
\item
every $K$-finite
 matrix coefficient of $\pi$ is bounded by
$$C \Xi^{1-\gamma}(g) (1+\|\log(a(g))\|)^s \qquad (C, s \geq 0, a(g) \in A^+);$$
\item every leading exponent $v$ of $\pi$ satisfies
$$ \Re(v)(H) \leq \ (\gamma-1) \rho(G)(H) \qquad (\forall \ H \in \f a^+).$$
\end{enumerate}
\end{thm}
So if $\gamma=0$ then $\pi$ is tempered and $\gamma=1$ if $\pi$ is trivial. For $G=Mp_{2n}(\mathbb R)$ or $G=O(p,q)$, the last statement is equivalent to
$$\Re(v) \preceq (\gamma-1) \rho(G).$$

\section{Associated Variety, Asymptotic Cycle and Wave Front Set}
Recall that $U(\f g)$ has a natural filtration
$$ \mathbb C= U_0(\f g) \subseteq  U_1(\f g) \subseteq  U_2(\f g) \subseteq   U_3(\f g) \subseteq \ldots \subseteq  U_n(\f g) \subseteq \ldots.$$
\subsection{Annihilator and Associated Variety}
Let $\pi$ be an admissible irreducible representation of a reductive group $G$. Let $V(\pi)$ be the Harish-Chandra module of $\pi$. Then
$V(\pi)$ is a $U(\f g)$ module.  Consider the annihilator of $V(\pi)$,
$$ Ann(V(\pi)) =
\{ D \in U(\f g) \mid \pi(D) V(\pi)=0  \}.$$
$Ann(V(\pi))$ is an ideal of $U(\f g)$ and inherits a filtration from the standard filtration of $U(\f g)$.
It follows that the induced graded algebra $gr(Ann(V(\pi)))$  is an ideal of $gr(U(g))=S(g)$ and necessarily commutative. 
\begin{no}
Let $\mc V(Ann\;\pi)$ be the associated variety of $gr(Ann(V(\pi)))$.
\end{no} 
\begin{thm}[Borho-Brylinski \!~\cite{bb}, Joseph \!~\cite{jo}] Suppose $G$ is a connected semisimple group and $\pi$ is an irreducible representation of $G$. Let $\f g_{\mathbb c}^*$ be the dual space of $\f g_{\mathbb C}$. Then $\mc V(Ann\;\pi)$ is the closure of a single
nilpotent orbit in $\f g_{\mathbb C}^*$.
\end{thm}
For $G$ having a finite number of components, $\mc V(Ann\;\pi)$ will be the closure of a finite number of nilpotent coadjoint
orbits of equal dimension. In fact, $G/G_0$ acts on the set of nilpotent coadjoint orbits. The associated variety $\mc V(Ann\;\pi)$ will be the closure of the union of one nilpotent orbit with its translations under $Ad(G/G_0)$. For a detailed account of associated varieties of Harish-Chandra modules, see \!~\cite{vogan89}.
\subsection{Asymptotic Cycle and Wave Front Set}
Let $G$ be a semisimple Lie group with finitely many components. Let us consider the global character $\Theta_{\pi}$. One can lift $\Theta_{\pi}$ to an invariant distribution $D$ on $\f g$. There exists a Taylor expansion of $D$ near $0$,
$$D ( f( tx))\cong \sum_{i=-r}^{\infty} t^i D_i(f (x)).$$
In \!~\cite{bv0}, Barbasch-Vogan proved that the Fourier transform $\widehat{D_i}$ is supported on the nilcone of $\f g^*$. 
\begin{no}  Let $supp(AS(\pi))$  be the closure of the union of
supports of $\widehat{D_i}$. We call $supp(AS(\pi))$ the support of the asymptotic cycle of $\pi$.
\end{no}
\begin{thm}[Barbasch-Vogan]
Suppose $G$ is a connected semisimple Lie group. Then the real dimension of $supp(AS(\pi))$ is equal to the complex dimension of $\mc V(Ann\;\pi)$. $supp(AS(\pi))$ is a union of nilpotent orbits contained in $\f g \cap \mc V(Ann\;\pi)$. Let 
$r=\frac{\dim(AS(\pi))}{2}$. Then $D_{-r}$ is the lowest nonzero term of the Taylor expansion of $D$ and $supp(\widehat{D_{-r}})$ is of maximal dimension in $supp(AS(\pi))$.
\end{thm}
This theorem holds for $G$ with a finite number of components. See \!~\cite{bv0}. \\
\\
Another notion similar to $supp(AS(\pi))$ is the wave front set of $\pi$ defined by Howe (\!~\cite{howewave}). Originally, $WF(\pi)$ was defined as a closed subset of the cotangent bundle $T^* G$. Because of the $G$-action, $WF(\pi)$ 
can be regarded as a closed subset of $\f g^*$. For $G$ semisimple and $\pi$ irreducible,
$WF(\pi)$ is in the nilpotent cone of $\f g^*$. In \!~\cite{howewave}, Howe studied the behavior of wave front sets $WF(\pi)$ under
restrictions to certain subgroups.\\
\\
Suppose that $\pi$ is irreducible and $G$ is semisimple. Rossmann proved that $WF(\pi)$ and 
$supp(AS(\pi))$ are identical (see Theorem C, \!~\cite{ro}). In what follows, we will not distinguish
between $WF(\pi)$ and $supp(AS(\pi))$ as long as $\pi$ is irreducible. There are two basic facts the reader should keep in mind. The first fact is that $WF(\pi)$ lies in the real space $\f g^*$. The second fact is that the algebraic closure of $WF(\pi)$ is exactly the associated variety $\mathcal V(Ann\; \pi) \subset \f g^*_{\mathbb C}$.

\begin{no}
From now on, we will identify $\f g^*$ with $\f g$ using a fixed invariant bilinear form on $\f g$. We will also identify $\f g^*_{\mathbb C}$ with $\f g_{\mathbb C}$ as complexifications of $\f g^*$ and $\f g$.
\end{no}

\chapter{Nilpotent Orbits }\label{nilorbits}
An element $x$ in a Lie algebra $ \f g$ is called nilpotent if $ad(x)$ is nilpotent. Let $G$ be a Lie group with Lie algebra $\f g$. If $x$ is  nilpotent and $g \in G$, then $Ad(g) x$ is also nilpotent.  It follows that under the adjoint action of $G$, nilpotent elements are grouped into $G$-orbits.  Each $G$-orbit is called a nilpotent adjoint orbit. For semisimple Lie groups, there are finitely many nilpotent orbits
and the classification of nilpotent adjoint orbits is completely known.
We cite Collingwood-McGovern's book \!~\cite{cm} as the main reference for this chapter. Unlike in \!~\cite{cm}, nilpotent orbits in this paper depend on the Lie groups, not just on the Lie algebra. In particular, for the orthogonal groups, a nilpotent orbit may not be connected. For our convenience, we state the Springer-Steinberg theorem slightly differently than in \!~\cite{cm}. The reason for doing this is given in Theorem \!~\ref{reasonss}. The main results are Theorems \!~\ref{-1}, \!~\ref{reasonss}, \!~\ref{induced1}.
\section{Young Diagrams and Complex Nilpotent Orbits}
A sequence of positive integers
$$ \bold d= (d_1 \geq d_{2} \geq \ldots \geq d_{r-1} \geq d_{r} > 0=d_{r+1})$$
is said to be a partition of $n$ if
$n=\sum_{j=1}^r d_j$. Write
$\|\bold d \|= n$. 
\begin{no}
In this paper, a partition of $n$ will be represented by a Young diagram of $n$ boxes, arranged as follows: 
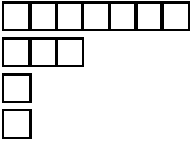
\begin{center}
\fb \fb \fb \fb \fb \fb \fb  \\
\fb \fb \fb \mbo \mbo \mbo \mbo \\
\fb \mbo \mbo \mbo \mbo \mbo \mbo \\
\fb \mbo \mbo \mbo \mbo \mbo \mbo \\
\end{center}
with the $i-th$ row of length $d_i$. The transpose of $\bold d$ is denoted by $\bold d^t$.
\end{no}
In our notation, the Young diagram above can be written as  $\bold d=(7,3,1,1)$ and its transpose
$\bold d^t =(4,2,2,1,1,1,1).$
\begin{defn}
If $d_j$ is odd for every $j \leq r$, we say $\bold d$ is {\it very odd}. If $d_j$ is even for every $ j \leq r$, we say $\bold d$ is {\it very even}.
If the nonzero $d_j$ are all distinct, we say $\bold d$ is {\it multiplicity free}. A row of even length is called an even row. A row of odd length is called an odd row.
\end{defn}

\subsection{Complex Nilpotent Orbits of $Sp_{2n}(\mathbb C)$}
\begin{defn}
A Young diagram $\bold d$ is said to be a symplectic Young diagram of size $2n$ if odd rows
occur with even
multiplicity and $\| \bold d \|=2n$. We use ${\mc YD}_{-}(2n)$ to denote the 
set of symplectic Young diagrams of size $2n$. 
\end{defn}
A symplectic group is the linear group that preserves a nondegenerate skew-symmetric
form.
\begin{thm}  Nilpotent adjoint orbits  of $Sp_{2n}(\mathbb C)$ are 
in one to one correspondence
with ${\mc YD}_{-}(2n)$.
\end{thm}
We denote the nilpotent adjoint orbit corresponding to $\bold d$ by $\mc O_{\bold d}(-)$.
The subscript $-$ is used to indicate that the group $G$ is defined with respect to a {\it skew}-symmetric form. If $\mathcal O_{\bold d}(-)$ is known to be symplectic, we  will simply write $\mathcal O_{\bold d}$.
For our convenience, we denote
the symplectic group $Sp_{2n}(\mathbb C)$ by $G(\mc O_{\bold d}(-))$ or $G(\mathcal O_{\bold d})$. The group $Sp_{2n}(\mathbb C)$ is thus attached
to the orbit implicitly. 

\subsection{Complex Nilpotent Orbits of $O(n, \mathbb C)$}
\begin{defn}
A Young diagram $\bold d$ is said to be an orthogonal Young diagram of size $n$ if even rows occur
with even multiplicity and $\|\bold d\|=n$. We use ${\mc YD}_{+}(n)$ to denote 
the set of all orthogonal Young diagrams of size $n$.
\end{defn}
An orthogonal group is a linear group preserving a nondegenerate symmetric form.
\begin{thm}  Nilpotent adjoint orbits  of $O(n,\mathbb C)$ are in one to one correspondence
with ${\mc YD}_{+}(n)$.
\end{thm}
We denote the nilpotent adjoint orbit corresponding to $\bold d$ by $\mc O_{\bold d}(+)$. The script $+$ is used to indicate that the group $G$ preserves a {\it symmetric} form. If the orbit $\mc O_{\bold d}(+)$ is known to be orthogonal, we will simply denote it by $\mc O_{\bold d}$.
In this context, the orthogonal group $O(n, \mathbb C)$ is  denoted by $G(\mc O_{\bold d}(+))$ or
$G(\mc O_{\bold d})$. \\
\\
Depending on the context,
$\mc O_{\bold d}$ can refer to either $\mc O_{\bold d}(+)$ or $\mc O_{\bold d}(-)$. 
\subsection{Operator $\bold -1$}
\begin{defn}
Let $\bold d$ be a Young diagram of size  $n$. Define $\bold d -\bold 1$ to be the new Young diagram
obtained by deleting the first column of $\bold d$. Define $\bold d- \bold i$ to be the new Young diagram
obtained by deleting the first $i$ columns of $\bold d$.
\end{defn}
Consider $\bold d =[3^2, 2^1, 1^2]$. The Young diagrams $\bold d$, $\bold d-\bold 1$ and $\bold d-\bold 2$ are listed as follows
\clearpage
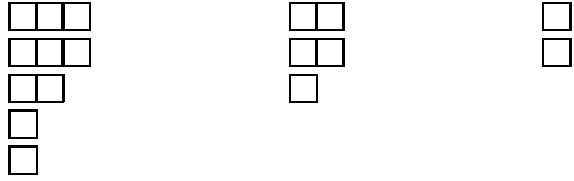
\begin{center}
\fb \fb \fb \hspace{1 in} \fb \fb \hspace{1 in} \fb \\
\fb \fb \fb \hspace{1 in} \fb \fb \hspace{1 in} \fb \\
\fb \fb \mbo \hspace{1 in} \fb \mbo \hspace{1 in} \mbo \\
\fb \mbo \mbo \hspace{1 in} \mbo \mbo \hspace{1 in} \mbo \\
\fb \mbo \mbo \hspace{1 in} \mbo \mbo \hspace{1 in} \mbo \\
\end{center}
Notice that $\bold d$ is a symplectic Young diagram,  
$\bold d - \bold 1$
is an orthogonal Young diagram, and  
$\bold d -\bold 2$
is a symplectic
Young diagram.
\begin{thm}~\label{-1}
If $\bold d$ is a symplectic Young diagram, then $\bold d - \bold 1$ is an orthogonal 
Young diagram. If $\bold d$ is an orthogonal Young diagram, then $\bold d -\bold 1$ is a
symplectic Young diagram.
\end{thm}
\section{Signed Young Diagrams and Real Nilpotent Orbits}
 Type I classical groups are subgroups of the general linear groups that preserve certain sesquilinear forms (see \!~\cite{li2} for the definition). Let $G$ be a real classical group of type I. Real nilpotent orbits of $G$ are simply nilpotent $G$-orbits in the real Lie algebra $\f g$. For type I classical groups, the real nilpotent orbits for $G$ are parameterized
 by equivalence classes of signed
Young diagrams. 
\begin{defn}[ Ch. 9.3 \!~\cite{cm}]
Signed Young diagrams are Young diagrams with $+$ or $-$
labeling the boxes in such a way that signs alternate across rows. Two signed Young diagrams are considered to be equivalent if one signed diagram can be obtained from the other by interchanging the rows of same lengths.
\end{defn}
Let $\bold D$ be a signed Young diagram. 
We will use $ D^+$ to denote the number of positive boxes in $\bold D$ and $D^-$ to denote the number of negative boxes in $\bold D$. We call $(D^+,D^-)$ the signature of $\bold D$. \\
\\
In this chapter, we are only interested in the nilpotent orbits of $Sp_{2n}(\mathbb R)$ and $O(p,q)$. 
\begin{no}
For the sake of our discussion, we fix a matrix realization for each $G$. So the group $O(p,q)$ and $Sp_{2n}(\mathbb R)$ in this paper, will come with a fixed sesquilinear form. 
\end{no}
The correspondence between signed Young diagrams and real nilpotent orbits depends on the sesquilinear forms. I would like to thank the referee for pointing this to me.

\subsection{Real Orbits of $Sp_{2n}(\mathbb R)$}
\begin{thm}[Springer and Steinberg]~\label{ss}
Nilpotent adjoint orbits of $Sp_{2n}(\mathbb R)$ are parametrized
by the equivalence classes of signed Young diagrams with the following properties:
\begin{itemize}
\item the signature of $\bold D$ is $(n, n)$;
\item for every $s$, rows of length $2s+1$ must occur with even multiplicity and must have their leftmost boxes labeled
$$-,+,-,+,\ldots,-,+$$
from the highest row to the lowest row.
\end{itemize}
\end{thm}
Clearly, the signed Young diagrams that satisfy our second condition must have signature $(n,n)$. So, the first condition is redundant. In the second condition, the choice of  the sign pattern of rows of odd lengths is artificial. In fact, it suffices that there are same number of boxes labeled $+$ and $-$ in the rows of length $2s+1$. Nevertheless, I have included this in
Theorem \!~\ref{ss} for two purposes. The first is to maximize the analogy with Theorem \!~\ref{ss2}. The second is to make the signed young diagrams easy to manipulate. \\
\\
We denote the set of signed Young diagrams in Theorem \!~\ref{ss} by ${\mc YD}_{-}(n,n)$. The following signed Young diagrams are in ${\mc YD}_{-}(n,n)$:
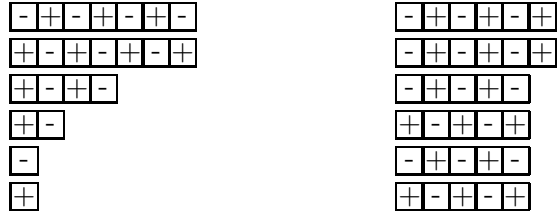
\begin{center}
\fm \fp \fm \fp \fm \fp \fm \hspace{1.in} \fm \fp \fm \fp \fm \fp \\
\fp \fm \fp \fm \fp \fm \fp \hspace{1 in} \fm \fp \fm \fp \fm \fp \\
\fp \fm \fp \fm \mbo \mbo \mbo \hspace{1. in} \fm \fp \fm \fp \fm \mbo \\
\fp \fm \mbo \mbo \mbo \mbo \mbo \hspace{1. in} \fp \fm \fp \fm \fp \mbo \\
\fm \mbo \mbo \mbo \mbo \mbo \mbo \hspace{1.in} \fm \fp \fm \fp \fm \mbo \\
\fp \mbo \mbo \mbo \mbo \mbo \mbo \hspace{1. in} \fp \fm \fp \fm \fp \mbo
\end{center}
Notice that $Mp_{2n}(\mb R)$-adjoint orbits coincide with $Sp_{2n}(\mathbb R)$- adjoint orbits. 
\begin{no}
We denote the nilpotent orbit corresponding 
to $\bold D$ by $\mc O_{\bold D}(-)$ or simply $\mc O_{\bold D}$ if $\bold D$ is specified to be symplectic. We denote the group $Mp_{2n}(\mathbb R)$ by $G(\mc O_{\bold D})$. 
\end{no}
The reader should notice that $G(\mathcal O_{\bold d})$ is a complex group of which $G(\mathcal O_{\bold D})$ is a real form.
\begin{defn}~\label{tauauto} Define
$$\tau(x)=\arr{I & 0 \\ 0 & -I } x \arr{I & 0 \\ 0 & -I } \qquad (x \in {\f{sp}}_{2n}(\mathbb R)).$$
\end{defn}
Notice that 
$$\arr{I & 0 \\ 0 & -I } \arr{0 & I_n \\ -I_n & 0} \arr{I & 0 \\ 0 & -I } =
-\arr{0 & I_n \\ -I_n & 0 }$$
and changing the symplectic form $\langle \ ,\  \rangle$ to $-\langle \ , \ \rangle$ results in the same group $Sp_{2n}(\mathbb R)$. Therefore $\tau$ defines an involution on $\f {sp}_{2n}(\mathbb R)$ and on $Sp_{2n}(\mathbb R)$. By sorting out the signs in the proof of Lemma 9.3.1 in \!~\cite{cm}, we obtain
\begin{lem}
\begin{enumerate}
\item
$\tau$ defines an involution on the real Lie algebra ${{\f{sp}}}_{2n}(\mathbb R)$.
\item  $\tau$ induces an involution on the set of nilpotent orbits, namely
$$\tau(\mc O_{\bold D}) = \mc O_{\tau(\bold D)}.$$
Here
$\tau(\bold D)$ is the signed Young diagram obtained by switching the signs ($+ \leftrightarrow -$) for the even rows of $\bold D$. 
\item The involution $\tau$ defines a diffeomorphism from $\mc O_{\bold D}$ onto $\mc O_{\tau(\bold D)}$.
\end{enumerate}
\end{lem}
The proof is left to the reader. \\
\\
Notice that the odd rows of $\tau(\bold D)$ remain the same as those of $\bold D$. For example, $\tau$:
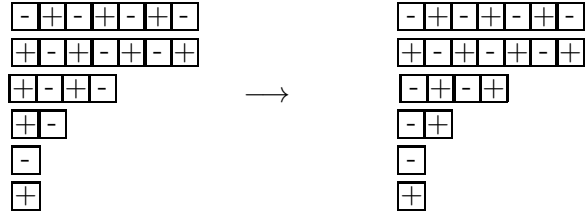
\begin{center}
\fm \fp \fm \fp \fm \fp \fm \hspace{1.in} \fm \fp \fm \fp \fm \fp \fm\\
\fp \fm \fp \fm \fp \fm \fp \hspace{1 in} \fp \fm \fp \fm \fp \fm \fp\\
\fp \fm \fp \fm \mbo \hspace{.48 in} $\longrightarrow$  \hspace{.48 in} \fm \fp \fm \fp \mbo \mbo \mbo \\
\fp \fm \mbo \mbo \mbo \mbo \mbo \hspace{1. in} \fm \fp \mbo \mbo \mbo \mbo \mbo  \\
\fm \mbo \mbo \mbo \mbo \mbo \mbo \hspace{1.in} \fm \mbo \mbo \mbo \mbo \mbo  \mbo \\
\fp \mbo \mbo \mbo \mbo \mbo \mbo \hspace{1. in} \fp \mbo \mbo \mbo \mbo \mbo  \mbo
\end{center}
One can further explore the involution $\tau$ on the Lie group level.
\begin{defn}
\begin{enumerate}
\item Define $\tau$ on $Sp_{2n}(\mb R)$ by
$$\tau(g)=\arr{I_n & 0 \\ 0 & -I_n } g \arr{I_n & 0 \\ 0 & -I_n}.$$
Clearly $\tau$ defines an involution on $Sp_{2n}(\mb R)$. $\tau$ is a topological homeomorphism and an isometry with respect to a certain left invariant
Riemannian metric.
\item Lift $\tau$ to an involution on $Mp_{2n}(\mb R)$. The lift $\tau$ exists and is unique. By abusing notation,
denote this involution by $\tau$.
\item Let $\pi \in \Pi(Mp_{2n}(\mb R))$. Define a new representation $(\pi^{\tau}, \mc H_{\pi})$ by
$$\pi^{\tau}(g)=\pi(\tau(g)).$$
\end{enumerate}
\end{defn}
\begin{lem}
$\pi \rightarrow \pi^{\tau}$ defines an involution on $\Pi(Mp_{2n}(\mb R))$ and on $\Pi_u(Mp_{2n}(\mb R))$.
\end{lem}

\subsection{Real Orbits of $O(p,q)$}
\begin{thm}[Springer \& Steinberg]~\label{ss2}
Nilpotent adjoint orbits of $O(p,q)$ are parametrized by equivalence classes of signed Young diagrams
with the following properties
\begin{itemize}
\item the signature of $\bold D$ is $(p,q)$;
\item for every $s$, rows of even length $2s$ must occur with even multiplicity and must
have their leftmost boxes labeled
$$+,-,+,-, \ldots, +,-$$
from the highest row to the lowest row.
\end{itemize}
\end{thm}
Again, the choice of sign patterns for rows of even length is artificial.
We denote the set of signed Young diagrams in Theorem \!~\ref{ss2} by ${\mc YD}_{+}(p,q)$.
The following signed Young diagrams are in ${\mc YD}_{+}(7,9)$:
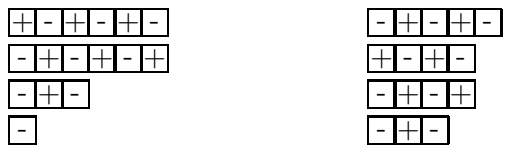
\begin{center}
\fp \fm \fp \fm \fp \fm \hspace{1.in} \fm \fp \fm \fp \fm \\
\fm \fp \fm \fp \fm \fp \hspace{1 in} \fp \fm \fp \fm  \mbo \\
\fm \fp \fm \mbo \mbo \mbo \hspace{1. in} \fm \fp \fm \fp \mbo \\
\fm \mbo \mbo \mbo \mbo \mbo \hspace{1. in} \fm \fp \fm \mbo \mbo 
\end{center}
\begin{no}
We use $\mc O_{\bold D}$ to denote the nilpotent orbit corresponding 
to $\bold D \in {\mc YD}_{+}(p,q)$. We use $G(\mc O_{\bold D})$ to denote the group $O(p,q)$. 
Any  $\bold D$ in this paper is interpreted as a signed Young diagram in a previously chosen set ${\mc YD}_{-}$ or
$\mc YD_{+}$. The orbit $\mc O_{\bold D}$ refers to either $\mc O_{\bold D}(+)$ or $\mc O_{\bold D}(-)$ with the understanding that $+$ or $-$ is implicitly known once $\bold D$ is given. 
\end{no}
Of course, if $\bold D$ is specified to be in ${\mc YD}_{+}$
or ${\mc YD}_{-}$, the notation $\mc O_{\bold D}$ causes no confusion.

\section{Nilpotent Orbits of Class $\mc U$}
\begin{no}
Let $\mc O_{\bold D}$ be a real nilpotent orbit. We use  $\bold d$ to denote the Young diagram obtained from $\bold D$ by removing the signs. 
\end{no}
Every real nilpotent orbit $\mc O$ induces 
a complex
nilpotent orbit $\mc O_{\mb C}$ by considering the complex group $G_{ad}$ acting
on $\mc O$ in $\f g$. This map is simply
$$\mc O_{\bold D} \rightarrow \mc O_{\bold d}$$
for orthogonal groups and symplectic groups. 
\begin{defn}
Let $\bold D$ be a signed Young diagram.
Define $\bold D-\bold 1$ to be the signed Young diagram obtained from $\bold D$ by deleting
the first column.
\end{defn}
\begin{thm}~\label{reasonss}
$-\bold 1$ defines an operation from ${\mc YD}_{+}(p,q)$ to the disjoint union
$$ \cup_{n \leq \min(p,q)} {\mc YD}_{-}(n,n)$$
$-\bold 1$ also defines an operation from ${\mc YD}_{-}(n,n)$ to the disjoint
union
$$ \cup_{\max(p,q) \leq n} {\mc YD}_{+}(p,q)$$
\end{thm}
By Theorem \!~\ref{ss} and Theorem \!~\ref{ss2},
$-\bold 1$ induces an orbital correspondence from real nilpotent orbits of
symplectic groups to real nilpotent orbits of orthogonal groups, and conversely.  
\begin{defn}~\label{prerigid}
We say that $\mc O_{\bold D}$ (or $\mc O_{\bold d}$) is pre-rigid if $\bold d^t$ is multiplicity free.
\end{defn}
This amounts to saying that every integer between $1$ and $d_1$ appears in the partition $\bold d$. In other words,
the partition $\bold d$ is of the following form
$$([d_1]^{*}, [d_1-1]^{*}, [d_1-2]^{*}, \ldots, [1]^{*})$$
with each multiplicity greater than zero. The reason why we call these orbit pre-rigid will become clearly in the next section. Roughly, they are not properly induced with a few exceptions.
\begin{defn}[Nilpotent orbits of Class $\mc U$]~\label{u}
\begin{enumerate}
\item $\mc U$ consists of a class of real nilpotent orbits $\mc O_{\bold D}$. For technical reasons, we {\bf exclude} from $\mc U$ those $\bold D$ whose last 2 columns are of the same length with at least two rows and are of the following forms
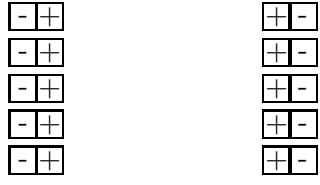
\begin{center}
\rm{ \fm \fp \hspace{1.in}   \fp \fm \\
  \fm \fp \hspace{1.in}   \fp \fm \\
\fm \fp \hspace{1.in}   \fp \fm \\
\fm \fp \hspace{1.in} \fp \fm  \\
\fm \fp  \hspace{1.in}  \fp \fm  }
\end{center}
(See the remark after this definition).
\item
Let $\bold D \in {\mc YD}_{-}(n,n)$ and $\bold d^t=(m_1 \geq m_2 \geq \ldots \geq m_{d_1})$. 
The nilpotent orbit $\mc O_{\bold D}$ of $Mp_{2n}(\mb R)$ is said to be in $\mc U(Mp_{2n}(\mb R))$ if $\bold d^t$ satisfies the following three conditions:
\begin{enumerate}

\item For every $i$, $ m_{2i} > m_{2i+1}$ and $ m_{2i+1} \geq m_{2i+2}$;
\item (1);
\item $\bold d^t$ is either very even or very odd.
\end{enumerate}
\item 
Let $\bold D \in {\mc YD}_{+}(p,q)$ and
$$\bold d^t=(m_1 \geq m_2 \geq \ldots \geq m_{d_1}).$$ 
The nilpotent orbit $\mc O_{\bold D}$ is said to be in $\mc U(O(p,q))$ if
$\bold d^t$ satisfies the following three conditions:
\begin{enumerate}

\item For every $i$, 
$m_{2i} \geq m_{2i+1}$ and $m_{2i+1} > m_{2i+2}$;
\item (1);
\item $\bold d^t$ is either very even or very odd.
\end{enumerate}
\end{enumerate}
\end{defn}
Here are two $\bold D$, one in $\mc U(Mp_{30}(\mb R))$, the other in $\mc U(O(8,14))$. 
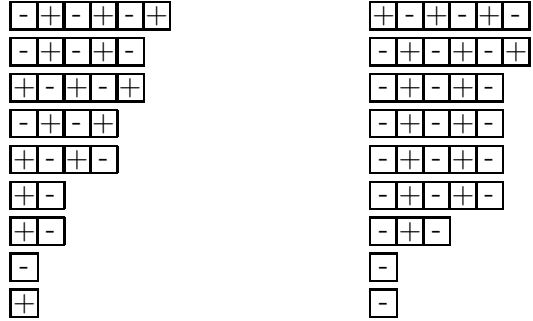
\begin{center}
\fm \fp \fm \fp \fm \fp  \hspace{1.in}   \fp \fm \fp \fm \\
\fm \fp \fm \fp \fm \mbo \hspace{1 in}   \fm \fp \fm  \fp \\
\fp \fm \fp \fm \fp \mbo \hspace{1. in}  \fm \fp \fm \mbo \\
\fm \fp \fm \fp \mbo \mbo  \hspace{1. in}   \fm  \fp \fm \mbo \\
\fp \fm \fp \fm \mbo \mbo \hspace{1. in}  \fm  \fp \fm \mbo \\
\fp \fm \mbo \mbo \mbo \mbo \hspace{1. in}  \fm  \fp \fm \mbo \\
\fp \fm \mbo \mbo \mbo \mbo \hspace{ 1. in}  \fm \mbo \mbo \mbo \\
\fm \mbo \mbo \mbo \mbo \mbo \hspace{ 1. in}   \fm \mbo \mbo \mbo \\
\fp \mbo \mbo \mbo \mbo \mbo \hspace {1. in}  \mbo \mbo \mbo \mbo
\end{center}
The condition $(a)$ is 
weaker than saying $\bold d^t$ is multiplicity free.
The condition $(c)$ amounts to saying that all multiplicities of $\bold d$ except perhaps the first are even. \\
\\
{\bf Remark.} Condition $(1)$ in our definition essentially excludes the following kind of real nilpotent orbits of $Mp_{2n}(\mathbb R)$ from the class of $\mc U$
\begin{center}
\rm{ \fm \fp  \fm   \fp \fm \fp \\
  \fm \fp  \fm   \fp \fm \fp \\
\fm \fp    \fm \fp \mbo \mbo \\
\fm \fp \fm \fp \mbo \mbo   \\
\fm \fp    \mbo \mbo \mbo \mbo \\
\fm \fp    \mbo \mbo \mbo \mbo
 }
\end{center}
There are a couple of reasons to exclude these orbits.
First of all, these orbits are not rigid. Hence excluding these orbits do not have any impact on our main results. Secondly, these orbits are induced at the complex orbit level. One can obtain infinitely many unitary representation attached to them by applying the cohomological induction process (\cite{k-v}). Indeed, when only two columns are involved, the corresponding representations can include all unitary highest or 
lowest weight modules.
Thirdly, our inductive construction will encounter  some analytic technicality involving the estimation of matrix coefficients and   cohomological induction is more straight forward for these orbits.   \\
\\
Recall
\begin{lem} [Thm. 7.3.5 and Prop. 6.3.7 \!~\cite{cm}] If $\mc O_{\bold d}$ is rigid, then $\mc O_{\bold d}$ is pre-rigid. If $\mc O_{\bold d}$ is special and rigid, then ${\bold d^t}$ must be either very even or very odd and must be multiplicity free.
\end{lem}
Based on the remark after Definition \ref{u}, we have 
\begin{cor} Special rigid orbits of $O(p,q)$ and $Sp_{2n}(\mb R)$ are contained in
$\mc U$.
\end{cor}
By Theorem \!~\ref{reasonss}, we have
\begin{cor}~\label{inu}
If $\mc O_{\bold D}$ is in $\mc U(Mp)$, then $\mc O_{\bold D-\bold 1}$ is in $\mc U(O)$.
If $\mc O_{\bold D}$ is in $\mc U(O)$, then $\mc O_{\bold D-\bold 1}$ is in $\mc U(Mp)$.
\end{cor}
\section{Induced Orbits}
\subsection{Complex Induced Orbits}
Let $\f g$ be either $\f o(n, \mb C)$ or ${\f{sp}}_{2n}(\mb C)$. Let $\f p=\f l \oplus \f n$ be a parabolic subalgebra. Let $\mc O_{\f l}$ be a nilpotent orbit in $\f l$. 
\begin{defn}
Define $\Ind_{\f l}^{\f g} \mc O_{\f l}$ to be the nilpotent $G$-orbit that intersects $ \mc O_{\f l} + \f n$ in an open dense set.
\end{defn}
Our definition differs slightly from the definition
in Ch. 7.1 of  \!~\cite{cm} as we require $\Ind_{\f l}^{\f g} \mc O_{\f l}$ be a $G$-orbit and $G$ may not be connected.  In any case, $\mc O_{\f g}$ is the $G$-orbit of the "generic" elements in $\mc O_{\f l}+ \f n$. In addition, $\Ind_{\f l}^{\f g} O_{\f l}$ only depends on $\f l$, not on the choice of $\f p$ (see Ch. 7.1\!~\cite{cm}).
\begin{lem}[Proposition 7.1.4 \!~\cite{cm}]
Let $\f l_1$ and $\f l_2$ be two Levi subalgebras of $\f g$ and $\f l_1 \subset \f l_2$. Then
$$\Ind_{\f l_2}^{\f g}(\Ind_{\f l_1}^{\f l_2} \mc O_{\f l_1})=\Ind_{\f l_1}^{\f g} \mc O_{\f l_1}.$$
\end{lem}
Thus orbital induction is \lq\lq associative\rq\rq. \\
\\
If $\mathfrak g=\f{sp}_{2n}(\mb C)$, then any Levi subalgebra $\f l$ is of the form
$$\f l=\f{sp}_{2n_0}(\mb C) \times \f{gl}(n_1, \mb C) \times \ldots \times \f{gl}(n_r, \mb C)$$
where $n_0, n_1, \ldots, n_r$ are nonnegative integers summing to $n$. Similarly, if $\mathfrak g=\f {o}(m, \mb C)$, then any Levi subalgebra $\f l$ is of the form
$$\f l= \f{o}(m_0, \mb C) \times \f{gl}(m_1, \mb C) \times \ldots \times \f{gl}(m_s, \mb C)$$
where $m_0, m_1, \ldots m_s$ are nonnegative integers and $m_0+2m_1+ \ldots + 2m_s=m$.\\
\\
 Recall that
nilpotent orbits of $SL(n, \mb C)$ are parametrized by partitions of $n$ in terms of the Jordan form. 
\begin{lem}Let $n=n_1+n_2$. Let $\f g= {\f{sl}}(n, \mb C)$.
Let $\f l$ be the block-diagonal matrices of size $(n_1, n_2)$ in $\f g$. Let $\mc O_{\bold s}$ and $ \mc O_{\bold t}$ be nilpotent orbits in ${\f{sl}}(n_1, \mb C)$ and ${\f{sl}}(n_2, \mb C)$ respectively. Then
$$\Ind_{\f l}^{\f g} (\mc O_{\bold s} \times \mc O_{\bold t})=\mc O_{\bold d}$$
with
$$d_j=s_j + t_j \qquad (\forall \ j ).$$
\end{lem}
\begin{defn}
We call $\bold d$ the merging of $\bold s$ and $\bold t$.
\end{defn}
For $G=O(n, \mb C)$ or $G=Sp_{2n}(\mb C)$, computation of induced orbits  is slightly more complicated than simply a \lq\lq merging\rq\rq. One needs the concept of \lq\lq collapse\rq\rq (see Ch. 6 \!~\cite{cm}). For the sake of simplicity, we will not introduce the concept of \lq\lq collapse\rq\rq. We only state one special result concerning $Sp_{2n}(\mb C)$. The general result can be found in \!~\cite{cm} ( Lemma 6.3.3 ) and is due to Gerstenhaber.  
\begin{lem}~\label{in} Let $G= Sp_{2n}(\mb C)$ and $L=GL(n-m, \mb C) \times Sp_{2m}(\mb C)$. Let $\mc O_{\bold s}$ be a nilpotent orbit of $Sp_{2m}(\mb C)$
such that the number of rows in $\bold s$ is less than or equal to $n-m$.
Then
$$\Ind_{\f l}^{\f g} \{0\} \times \mc O_{\bold s} =\mc O_{\bold d}$$
with
$$d_j=s_j+ 2  \qquad (\forall \, 1 \leq  j \leq n-m ).$$
\end{lem}
So under the assumption that the number of the rows in $\bold s$ is less or equal to $n-m$, $\bold d$ is the merging of $\bold s$ with two copies of $1^m $.
$\bold d$ remains symplectic since $\bold s$ is symplectic.
\subsection{Induced Real Orbits}
Let $G$ be a reductive Lie group and $L$ be a Levi subgroup. Let $P=LN$ be a parabolic subgroup of $G$. Let $\mc O_{\f l}$ be a nilpotent orbit in $\f l$. 
\begin{defn}
Define 
$\Ind_{\f l}^{\f g} \mc O_{\f l}$ to be the union of the nilpotent $G$-orbits $\mc O_{\f g}^{(i)}$ that contain an open subset of  $\mc O_{\f l} + \f n$.
\end{defn} 
Implicit in the notation is the fact that $\Ind_{\f l}^{\f g}$ does not depend on the choice of $\f n$. The proof of this fact is essentially the same as the proof for the complex orbits (see Theorem 7.1.3 of \!~\cite{cm}).
$\Ind_{\f l}^{\f g} \mc O_{\f l}$ for a classical group of type I may no longer be a single
nilpotent orbit. Nevertheless, it is contained in a single complex nilpotent orbit.
\begin{lem}
We have the following commutative diagram:
\begin{equation}~\label{commu}
\begin{CD}
\mc O_{\f l} @>{complex \ orbit} >> \mc O_{\f l_{\mathbb C}} \\
@VV{\Ind}V @VV{\Ind}V \\
\Ind_{\f l}^{\f g} \mc O_{\f l} @>{complex \ orbit}>> \Ind_{\f l_{\mathbb C}}^{\f g_{\mathbb C}} \mc O_{\f l_{\mathbb C}}
\end{CD}
\end{equation}
\end{lem} 
Let $\mc O_{\bold S}$ be a nilpotent orbit in $\f {sp_{2m}(\mb R)}$.
By the lemma above,
$$\Ind_{\f{gl}(n-m, \mb R) \oplus {\f{sp}}_{2m}(\mb R)}^{{\f{sp}}_{2n}(\mb R)} \{0\} \times \mc O_{\bold S}$$
must be contained in the complex induced orbit
$$\Ind_{\f{gl}(n-m, \mb C) \oplus {\f{sp}}_{2m}(\mb C)}^{{\f{sp}}_{2n}(\mb C)} \{0\} \times \mc O_{\bold s}.$$
Under the assumption in Lemma \!~\ref{in}, this complex induced orbit is parametrized by $\bold d$ with
$$ d_j = s_j + 2 \ \ (\forall \, 1 \leq j \leq n-m).$$
So 
$$\Ind_{\f{gl}(n-m, \mb R) \oplus {\f{sp}}_{2m}(\mb R)}^{{\f{sp}}_{2n}(\mb R)} \{0\} \times \mc O_{\bold S}$$
is contained in the intersection of $\mc O_{\bold d}$ with $\f{sp}_{2n}(\mb R)$. We have
\begin{thm}~\label{induced1} 
Let  $\bold S \in \mc{YD}_{-}(m,m)$ and $\bold s=(s_1 \geq s_2 \geq \ldots \geq s_r > 0)$. Suppose that $n-m \geq r$. Then
$$\Ind_{\f{gl}(n-m, \mb R) \oplus {\f{sp}}_{2m}(\mb R)}^{{\f{sp}}_{2n}(\mb R)} \{0\} \times \mc O_{\bold S}= \cup_{j=0}^{n-m-r} \mc O_{\bold D^{(j)}}.$$
Each $\bold D^{(j)} \in \mc{YD}_{-}(n,n)$
is uniquely defined as follows:
\begin{enumerate}
\item merge one column of length ${n-m}$ to ${\bold S}$ from the left;
\item merge one column of length ${n-m}$ to ${\bold S}$ from the right;
\item extend the signs of $\bold S$ for rows of even lengths;
\item fill in the $n-m-r$ rows of length $2$ in $\bold D^{(j)}$  with
$j$  {\bf{ \fm \fp}} \, s and $n-m-r-j$  {\bf{ \fp \fm}} \, s;
\item the signs of $\bold D^{(j)}$ for rows of odd length are uniquely dictated by the rules we set for ${\mc YD}_{-}(n,n)$, thus not subject to 
change.
\end{enumerate}
 \end{thm}
We give an example first and a sketch of the proof will follow.\\
\\
{\bf Example}  Let $\bold t= 1^{n-m}$ be the partition of $n-m$ which corresponds to the zero orbit for $GL(n-m, \mb R)$. Then we can write the conclusion of Theorem \!~\ref{induced1} as 
$$\Ind_{\f{gl}(n-m, \mb R) \oplus {\f{sp}}_{2m}(\mb R)}^{{\f{sp}}_{2n}(\mb R)} \mc O_{\bold t} \times \mc O_{\bold S}= \cup_{j=0}^{n-m-r} \mc O_{\bold D^{(j)}}.$$
Suppose $\bold S$ is the signed partition in the middle of the diagram below and $n-m=7$. 
\begin{center}
\fb \hspace{1. in} \fm \fp \fm \fp \fm \hspace{1.in} \fb \\
\fb \hspace{1. in}  \fp \fm \fp \fm \fp \hspace{1.in} \fb \\
\fb \hspace{1. in} \fm \fp \fm \fp \mbo \hspace{1.in} \fb \\
\fb \hspace{1. in} \fm \fp \mbo \mbo \mbo \hspace{1.in} \fb \\
\fb \hspace{1. in} \fm \fp \mbo \mbo \mbo \hspace{1.in} \fb \\
\fb \hspace{1. in} \mbo \mbo \mbo \mbo \mbo \hspace{1.in} \fb \\
\fb \hspace{1. in} \mbo \mbo \mbo \mbo \mbo \hspace{1.in} \fb \\
\end{center}
The procedure explained
in Theorem \!~\ref{induced1} will result in the following three signed symplectic Young diagrams:
\begin{center}
\fm \fp \fm \fp \fm \fp \fm \hspace{1. in} \fm \fp \fm \fp \fm \fp \fm \hspace{1. in} \fm \fp \fm \fp \fm \fp \fm \\
\fp \fm \fp \fm \fp \fm \fp \hspace{1.in} \fp \fm \fp \fm \fp \fm \fp \hspace{1. in} \fp \fm \fp \fm \fp \fm \fp \\
\fp \fm \fp \fm \fp \fm \mbo \hspace{1. in} \fp \fm \fp \fm \fp \fm \mbo \hspace{1. in} \fp \fm \fp \fm \fp \fm \mbo \\
\fp \fm \fp \fm \mbo \mbo \mbo \hspace{1.in} \fp \fm \fp \fm \mbo \mbo \mbo  \hspace{1.in} \fp \fm \fp \fm \mbo \mbo \mbo \\
\fp \fm \fp \fm \mbo \mbo \mbo \hspace{1.in} \fp \fm \fp \fm \mbo \mbo \mbo \hspace{1.in} \fp \fm \fp \fm \mbo \mbo \mbo \\
\fp \fm \mbo \mbo \mbo \mbo \mbo \hspace{1.in} \fp \fm \mbo \mbo \mbo \mbo \mbo \hspace{1.in} \fm \fp \mbo \mbo \mbo \mbo \mbo \\ 
\fp \fm \mbo \mbo \mbo \mbo \mbo \hspace{1.in} \fm \fp \mbo \mbo \mbo \mbo \mbo \hspace{1. in} \fm \fp \mbo \mbo \mbo \mbo \mbo
\end{center}
 These three signed Young diagrams are $\bold D^{(0)}$, $\bold D^{(1)}$ and $\bold D^{(2)}$.\\
 \\
Sketch of the proof:  $(1), (2)$ are evident by the commutative diagram \!~\ref{commu} and $(5)$ follows from our convention (see Theorem \!~\ref{ss}). Let
$$\bold d=(s_1+2, s_2+2, \ldots, s_r+2, \overbrace{2, \ldots, 2}^{n-m-r}).$$
An odd number $2l+1$ must occur even times. {\bf We group $2l+1$'s into pairs}.\\
\\
We adopt the notation from Ch 5. 2 of \!~\cite{cm} and the main reference is Ch 5.2 and 9.3 of \!~\cite{cm}. A signed Young diagram can be expressed as 
$$a_1^{\pm} \oplus a_2^{\pm} \oplus \ldots a_s^{\pm}$$
with $a_i^{\pm}$ symbolizing a row of length $a_i$ beginning with the sign $\pm$. Notice that odd number $2l+1$ always comes in pairs as $((2l+1)^+, (2l+1)^-)$ according to our convention.\\
\\
Fix a symplectic space $(\mb R^{2n}, \Omega)$ and a nondegenerate symplectic subspace $\mb R^{2m}$.  Decompose $\mb R^{2n}$ into
$$\mb R^{2m} \oplus V_0 \oplus V_0^{\prime}$$ such that 
 $V_0 \oplus V_0^{\prime}$ is the orthogonal complement of $\mb R^{2m}$ with respect to $\Omega$ and $V_0$ and $V_0^{\prime}$ are both isotropic subspaces. Then $\Omega|_{V_0 \oplus V_0^{\prime}}$ is nondegenerate. Let
 $$\f n = \{ u \in \f{sp}_{2n}(\mb R) \mid u V_0=0, u ( \mb R^{2m} \oplus V_0) \subseteq V_0 \}.$$
 It is easy to see that for any $u \in \mathfrak n$ we have $ker(u) \supseteq V_0$ and $ u \mb R^{2n} \subseteq V_0 \oplus \mb R^{2m}$.\\
 \\
Let $X_{\bold S}$ be an element in $\mc O_{\bold S} \subset \f{sp}_{2m}(\mb R)$. 
Fix an arbitrary $u \in \f n$ such that $X_{\bold S}+ u \in \mc O_{\bold d}$. Then $X_{\bold S}+u \in \mc O_{\bold D}$ for some $\bold D$. Our main task is to determine all possible $\bold D$. Since the signs for odd rows of $\bold D$ are prefixed, we only need to determine the signs for even rows of $\bold D$.  \\
\\
Let $\{H, X_{\bold S}+u, Y \}$ be a standard triple.  Then
$\mb R^{2n}$ decomposes into a direct sum of subspaces
$$V_1 \oplus V_2 \oplus \ldots \oplus V_{r} \oplus V_{r+1} \oplus \ldots \oplus V_{n-m},$$
such that
\begin{itemize}
\item $\dim V_i= d_i=s_i+2$ for $i \in [1,r]$ and  $\dim V_{i}=2$ for $i > r$.
\item For even $d_i$, $\Omega|_{V_i}$ is nondegenerate.
\item  For an odd pair $d_i=d_{i+1}$, $\Omega|_{V_i \oplus V_{i+1}}$ is nondegenerate, $V_i$ and $V_{i+1}$ are both isotropic.
\item All $ V_i $ for even $d_i$   and $ V_i \oplus V_{i+1} $ for $(d_i, d_{i+1})$  an  odd  pair, are perpendicular to each other with respect to the symplectic form $\Omega$.  In other words,
$$\Omega= [\oplus_{i \ even} \Omega|_{V_i}] \oplus [\oplus_{(d_i, d_{i+1}) \ an \ odd \ pair} \Omega|_{V_i \oplus V_{i+1}}].$$
\item  $X_{\bold S}+u$ decomposes accordingly
as 
$$X_1 \oplus X_2 \oplus \ldots \oplus X_r \oplus X_{r+1} \ldots \oplus X_{n-m},$$
where
$X_i$ is a principal nilpotent element in $\f{sp}(V_i, \Omega|_{V_i})$ for even dimensional $V_i$ and $X_i \oplus X_{i+1}$ is a nilpotent element in $\mc O_{d_i^+, d_{i+1}^-} \subset \f{sp}(V_i \oplus V_{i+1}, \Omega|_{V_i \oplus V_{i+1}})$ for the odd pair $(d_i= d_{i+1})$.
\end{itemize}
This is essentially the abstract version of Jordan decomposition for a nilpotent element in a real symplectic Lie algebra. 
For even dimensional $V_i$, there are two principal nilpotent classes parametrized by $d_i^+$ and $d_i^-$ respectively. We now need to determine the sign attached to $X_i$ for each even $d_i$. \\
\\
Let $d_i$ be an even integer. Generally speaking, the sign attached to $X_i$  can be determined as follows.
There exists a basis of $V_i$
 $$ \{ e_i^{(1)}, e_i^{(2)}, \ldots, e_i^{(d_i)} \}$$
such that
$$X e_i^{(j)}= e_i^{(j+1)} \qquad (\forall \ j \in [1, d_i-1]) ,$$
and $X_i e_i^{(d_i)}=0$. Then the signature of
$$\Omega(X^{d_i-1} e_i^{(1)}, e_i^{(1)})$$
gives the sign attached to $d_i$ for $X_i$ (see Page 139 \!~\cite{cm}). For the standard triple $\{h, x, y \} \subset \mathfrak{sl}_2(\mb R)$, $\mc O_{x}$ will correspond to $2^+$ and $\mc O_{-x}$ will correspond to $2^-$.\\
\\
Suppose that $d_i$ is an even integer greater than $2$.
Now we need to connect the signature attached to $X_i$ with its \lq\lq projection \!\rq\rq to $\mb R^{2m}$.
Notice that $X_{\bold S} V_0=0$ and $u V_0=0$. So $V_0 \subseteq ker(X_{\bold S}+u)$. Since $X_{\bold S}+u \in \mc O_{\bold d}$, the kernel of $X_{\bold S}+ u$ is of dimension $n-m$ which is equal to the dimension of $V_0$. We see that 
$ker(X_{\bold S}+u)=V_0$. Hence $$(X_{\bold S}+u) (\mb R^{2n})= V_0^{\perp}=V_0 \oplus \mb R^{2m}.$$
In particular,
 $${\rm{span}} \{e_i^{(2)}, \ldots, e_i^{(d_i)} \}=X_i (V_i)=V_i \cap (V_0  \oplus \mb R^{2m}).$$
 For $j \in [2, d_i]$, put $e_i^{(j)}=f_i^{(j)}+ v_i^{(j)}$ where $f_i^{(j)} \in \mb R^{2m}$ and $v_i^{(j)} \in V_0$.
 We have for $j \in [2, d_i-1]$,
  $$(X_{\bold S}+u)(f_i^{(j)}+v_i^{(j)})=(X_{\bold S}+u) e_i^{(j)}= X_i(e_i^{(j)})=e_i^{(j+1)}=f_i^{(j+1)}+v_i^{(j+1)}.$$
 Since $u(\mb R^{2m} \oplus V_0) \subseteq V_0$ and $X_{\bold S} V_0=0$, we have
 $$X_{\bold S} f_i^{(j)}=f_i^{(j+1)}, \qquad u e_i^{(j)}=v_i^{(j+1)}.$$
 Let $W_i={\rm span}\{f_i^{(j)} \mid j \in [2,d_i-1] \}$. 
Notice that $\Omega$ restricted to ${\rm span}\{e_i^{(2)}, \ldots, e_i^{(d_i-1)} \}$ is nondegenerate. Since $V_0$ is isotropic and $W_i \perp V_0$, we see that $\Omega$ restricted onto $W_i$ is nondegenerate. Furthermore, $W_i$ are mutually perpendicular with respect to $\Omega$ since $V_i$'s are perpendicular to each other. So $X_{\bold S}|_{W_i}$ is a Jordan block of size $d_i-2=s_i$. The sign attached to it, is the signature of
$$\Omega(X_{\bold S}^{s_i-1} f_i^{(2)}, f_i^{(2)})=\Omega((X_{\bold S}+u)^{s_i-1} e_i^{(2)}, e_i^{(2)}).$$
Observe that
$$\Omega((X_{\bold S}+u)^{s_i-1} e_i^{(2)}, e_i^{(2)})=\Omega((X_{\bold S}+u)^{s_i} e_i^{(1)}, (X_{\bold S}+u) e_i^{(1)})=-\Omega((X_{\bold S}+u)^{s_i+1} e_i^{(1)}, e_i^{(1)}).$$
Hence the sign attached to $d_i$ for $X_{\bold S}+u$ differs from the sign attached to $s_i$ for $X_{\bold S}$. $(3)$ is proved. We have shown that
$$\Ind_{\f{gl}(n-m, \mb R) \oplus {\f{sp}}_{2m}(\mb R)}^{{\f{sp}}_{2n}(\mb R)} \{0\} \times \mc O_{\bold S} \subseteq \cup_{j=0}^{n-m-r} \mc O_{\bold D^{(j)}}.$$
\\
Lastly, we need to show that each $\mc O_{\bold D^{(j)}}$ appears in 
$$\Ind_{\f{gl}(n-m, \mb R) \oplus {\f{sp}}_{2m}(\mb R)}^{{\f{sp}}_{2n}(\mb R)} \{0\} \times \mc O_{\bold S}.$$
It suffices to construct an element $X_{\bold D^{(j)}} \in Y_{\bold S} + \f n$ for an element $Y_{\bold S} \in \mc O_{\bold S}$. Let $\{ e_1, e_2, \ldots, e_{n-m} \}$ be a basis for $V_0$ and $\{ e_{n-m+1}, e_{n-m+2}, \ldots, e_{2n-2m} \}$ be the dual basis for $V_0^{\prime}$ such that the natural bilinear form coincides with the restriction of $\Omega$ onto $V_0 \otimes V_{0}^{\prime}$. Fix a standard basis $\{e_{2n-2m+1}, \ldots, e_{2n} \}$ for
$(\mb R^{2m}, \Omega|_{\mb R^{2m}})$. Let $E_{i,j}$ be the $2n \times 2n$ matrix with $(i,j)$-th entry equal to $1$ and the rest equal to zero.  Without loss of generality, suppose that 
$$Y_{\bold S}=Y_1 \oplus Y_2 \oplus \ldots Y_r$$
is a Jordan decomposition such that 
\begin{itemize}
\item For even $s_i$, $$Y_i=\pm E_{j_i+\frac{s_i}{2}, m+j_i+\frac{s_i}{2}}+ \sum_{k=1}^{\frac{s_i}{2}-1} E_{j_i+k, j_i+k+1}-E_{m+j_i+k+1, m+j_i+k}$$ 
for some $j_i \geq 2n-2m$;
\item For a pair of odd integer $(s_i=s_{i+1})$, $$Y_i=\sum_{k=1}^{s_i-1} E_{j_i+k, j_i+k+1} \qquad  Y_{i+1}=-\sum_{k=1}^{s_i-1} E_{m+j_i+k+1, m+j_i+k},$$ 
for some $j_i \geq 2n-2m$.
\end{itemize}
For even $s_i$, define 
$$X_i=E_{i,j_{i}+1}+ Y_i-E_{m+j_{i}+1, n-m+i}.$$ 
 For odd pair $(s_i=s_{i+1})$, define 
 $$X_i=E_{i, j_{i}+1}+Y_i+E_{j_i+s_i, i+1}, \qquad X_{i+1}=-E_{m+j_{i}+1,n-m+i}+Y_{i+1}-E_{n-m+i+1, j_i+s_i+m}.$$
 For $i \in [r+1, n-m]$, define $X_i= \pm E_{i, n-m+i}$ such that the number of $-$ is $j$ and number of $+$ is
 $n-m-r-j$. Let $X=\sum_{i=1}^{n-m} X_i$. It is then easy to check that $X \in \mc O_{\bold D^{(j)}}$. $\Box$.\\
 \\
Our method can be generalized to compute induced real nilpotent orbits. We shall do this in a subsequent paper.
\begin{cor}~\label{induced0}  Let $\bold t=[1^{n}]$.
Then 
$$\Ind_{{\f{gl}}(n, \mb R)}^{{\f{sp}}_{2n}(\mb R)} \mc O_{\bold t}$$
is the union of $n+1$ orbits consisting of all $\mc O_{\bold D^{(j)}}$ with 
$\bold d^{(j)}=[2^n]$.
\end{cor}
\begin{thm}~\label{tauinduction} 
Let $\mc O_{\bold S}$ be a nilpotent orbit of $Sp_{2m}(\mb R)$. Let $\bold s=(s_1 \geq s_2 \geq \ldots \geq s_r > 0)$. Suppose $n-m \geq r$.
Then $$\Ind_{{\f{sp}}_{2m}(\mb R) \oplus \f{gl}(n-m, \mb R)}^{{\f{sp}}_{2n}(\mb R)} \tau(\mc O_{\bold S}) \times \{0\}$$
is a union of nilpotent orbits $\mc O_{\bold D^{(j)}}$ and each $\bold D^{(j)}$
can be obtained by 
\begin{enumerate}
\item merging two copies of $[1^{n-m}]$ to ${\bold S}$ from the left;
\item extending the signs of $\bold S$;
\item filling the $n-m-r$ rows of length $2$ in $\bold D^{(j)}$  with
$j$ {\rm \fm \fp} \ and $n-m-r-j$ {\rm \fp \fm} .
\end{enumerate}
\end{thm}
Sketch of the proof: Recall that $\tau(\mc O_{\bold S})=\mc O_{\tau(\bold S)}$ and $\tau(\bold S)$ is obtained by switching the signs of even rows in $\bold S$. Suppose that $s_i^{\pm}$ represents an even row in $\bold S$. Then
$s_i^{\mp}$ represents the corresponding even row in $\tau(\bold S)$. By Theorem \!~\ref{induced1}, $d_i^{\pm}$ will represent the corresponding even row in $\bold D^{(j)}$.
 This is the same as adding two boxes to the even rows of $\bold S$ and extend the signs of $\bold S$. $(3)$ follows directly from Thm. \!~\ref{induced1}. $\Box$

\chapter{Theta Correspondences and Quantum Induction }
\begin{no}~\label{m}
Let $Mp_{2n}(\mb R)$ be the metaplectic covering of $Sp_{2n}(\mb R)$.
Let $\{1, \epsilon\}$ be the preimage of $1 \in Sp_{2n}(\mb R)$.  
For any subgroup $G$ of $Sp_{2n}(\mb R)$, let $MG$ be the preimage of $G$ under the metaplectic covering.
\end{no}
Fix a dual pair $(O(p,q), Sp_{2n}(\mb R))$ in $Sp_{2n(p+q)}(\mb R)$ (see \!~\cite{howe79}).
 Fix a maximal compact subgroup $K$ for $Mp_{2n(p+q)}(\mb R)$ such that $K \cap MO(p,q)$ is a maximal compact subgroup of $MO(p,q)$ and $K \cap MSp_{2n}(\mb R)$ is a maximal compact subgroup of $MSp_{2n}(\mb R)$.
Let $\omega(p,q;2n)$ be the oscillator representation of the group $Mp_{2n(p+q)}(\mb R)$. If the parameters $(p,q;2n)$ are apparent, 
we will write $\omega$. Let $V(\omega)$ be the Harish-Chandra module.
Unless otherwise stated, any representation $\pi$ of $MG$ in this paper satisfies
$\pi(\epsilon)=-1$.
Very often, this assumption reduces our discussion and computation to representations of the linear group. 
\begin{no}Let $G$ be any reductive subgroup of $Sp$.
Let $\mc R(MG, \omega)$ be the equivalence classes of irreducible Harish-Chandra modules of $MG$ that occur as quotients of $V(\omega)$.
\end{no}
\begin{thm}[Howe \!~\cite{howe}]
Let $(G_1, G_2)$ be a dual pair in $Sp$. Then $\mc R(MG_1 G_2, \omega)$ yields a one-to-one correspondence between $\mc R(MG_1, \omega)$ and $\mc R(MG_2, \omega)$.
\end{thm}
This correspondence is often called (local) theta correspondence, Howe's correspondence or duality correspondence. Let 
$\theta(p,q;2n)$ be the theta correspondence from $\mc R(MO(p,q),\omega)$ to
$\mc R(MSp_{2n}(\mb R), \omega)$. Let $\theta(2n;p,q)$ be its inverse. If the parameters $(p,q;2n)$ are apparent, we will just write $\theta$. The description of the sets $\mc R(MO(p,q),\omega)$ and
$\mc R(MSp_{2n}(\mb R), \omega)$ is  not known in general. \\
\\
If $G$ is a subgroup of $Sp$, we say that $MG$ splits if $MG \simeq G \times \{\pm 1\}$. For the dual pair $(O(p,q), Sp_{2n}(\mathbb R))$, $MO(p,q)$ always splits and $MSp_{2n}(\mathbb R)$ splits if and only if $p+q$ is even. If $p+q$ is odd, $MSp_{2n}(\mb R)$ is isomorphic to the metaplectic group $Mp_{2n}(\mb R)$.
In all cases, $\theta$ can be regarded as a one-to-one correspondence between a
certain subset of $\Pi(Mp_{2n}(\mb R))$ and a certain subset of $\Pi(O(p,q))$. Since this is the viewpoint in some literature (for example \!~\cite{kr}), we will also take this viewpoint if necessary.

\section{ Theta Correspondence in Semistable Range}
We adopt the notation from (1.3.1) and (1.3.2). Let $\bold n$ be the constant vector
$(n,n,n, \ldots, n)$ of a fixed dimension.
\begin{defn}[Def. 5.1 \!~\cite{theta}, Def. 3.2.1 \!~\cite{basic}, Thm. 8.17, 8.18, Def. 8.3, 8.4 \!~\cite{thesis}]
An irreducible representation $\pi$ of $MO(p,q)$ is said to be
in the {\it semistable range} of
$\theta(p,q;2n)$
if and only if every leading exponent $v$ of $\pi$ satisfies
$$ \Re (v)-\bold n+ 2 \rho(O(p,q)) \prec 0.$$
An irreducible representation $\pi$ of $MSp_{2n}(\mb R)$ is said to be
in the semistable range of $\theta(2n;p,q)$ if and only if every leading exponent
$v$ of $\pi$ satisfies
$$\Re(v)-\bold {\frac{p+q}{2}}+ 2 \rho(Sp_{2n}(\mb R)) \prec 0.$$
\end{defn}
We denote the semistable ranges by $R_{s}(MO(p,q), \omega)$ and $R_s(MSp_{2n}(\mb R), \omega)$ respectively. The reader should note that a representation $\pi$ in $R_{s}(G, \omega)$ may not occur in $\mc R(G, \omega)$. This is pointed out to me by Chen-Bo Zhu.
\begin{defn}[\!~\cite{li2}, \!~\cite{theta}]~\label{thetas}
Consider $(G_1, G_2)=(O(p,q), Sp_{2n}(\mb R))$ or
$(G_1,G_2)=(Sp_{2n}(\mb R), O(p,q)).$
Let $\pi \in R_s(MG_1, \omega)$. Define a bilinear form $(\ , \ )_{\pi}$ on $V(\omega) \otimes V({\pi}^c)$ 
$$(\phi \otimes u, \psi \otimes v)_{\pi}=\int_{MG_1} (\omega(g) \phi, \psi)( v, \pi(g) u) d g \qquad (\phi, \psi \in V(\omega), u, v \in V({\pi})).$$
 Let $\mc R_{\pi}$ be the radical of $(\ , \ )_{\pi}$. Define
$$\theta_s(MG_1, MG_2)(\pi)= V(\omega) \otimes V(\pi^c) / \mc R_{\pi}.$$
\end{defn}
$\theta_s(MG_1, MG_2)(\pi)$ inherits an infinitesimal $MG_2$ action from
$\omega(MG_1, MG_2)$. It is a $(\f g, K)$-module of $MG_2$.
\begin{thm}[Theorem 1.1 \!~\cite{theta}, \!~\cite{unit}]~\label{irreducibility}
Suppose $\pi$ is a unitarizable Harish-Chandra module in the semistable range of $\theta(MG_1,MG_2)$.
Then $(\ , \ )_{\pi}$ is well-defined.
If $(\ , \ )_{\pi} \neq 0$, then $\theta_s(MG_1, MG_2)(\pi)$ is an irreducible Hermitian Harish-Chandra module of $MG_2$ and $\theta_s(MG_1, MG_2)(\pi)$ is equivalent to $\theta(MG_1, MG_2)(\pi)$.
\end{thm}
This theorem basically says that if $(\ , \ )_{\pi}$ is well-defined and nonvanishing, then $\theta_s$ is Howe's correspondence (\!~\cite{howe}) on the Harish-Chandra module level. By Howe's Theorem, $\theta_s(MG_1, MG_2)(\pi)$ is an irreducible Harish-Chandra module of $MG_2$.
The notation $(\ , \ )_{\pi}$ is essentially due to Jian-Shu Li (\!~\cite{li2}).\\
\\
In \!~\cite{non}, we proved the following nonvanishing theorem.
\begin{thm}[Corollary 1.1 \!~\cite{non}] ~\label{nonvanishing} 
Suppose $p+q \leq 2n+1$. Suppose $\pi \in R_{s}(MO(p,q), \omega)$. Let $\det$ be the lift of the determinant of $O(p,q)$ to $MO(p,q)$, i.e.,
$$ \ker(\det) = MSO(p,q).$$
Then either $\theta_s(p,q;2n)({\pi}) \neq 0 $ or $\theta_s(p,q;2n)({\pi \otimes \det}) \neq 0$.
\end{thm}
\section{Unitarity and Strongly Semistable Range}
In \!~\cite{unit}, we proved the unitarity of $\theta_s(\pi)$ for $\pi$ in the strongly semistable range. 
\begin{defn}~\label{qi}
An irreducible admissible representation $\pi$ is  
in $R_{ss}(MO(p,q), \omega)$ if every leading exponent of $\pi$ satisfies
\begin{equation}~\label{ss1}
\Re(v)-(\bold{n-\frac{p+q}{2}})+\rho(O(p,q)) \preceq 0.
\end{equation}
An irreducible admissible representation $\pi$ of $MSp_{2n}(\mb R)$ is  in $R_{ss}(MSp_{2n}(\mb R), \omega)$ if every leading exponent of $\pi$ satisfies
\begin{equation}
\Re(v)-(\bold{\frac{p+q}{2}-n-1})+\rho(Sp_{2n}(\mb R)) \preceq 0.
\end{equation}
We call $R_{ss}$ the strongly semistable range.
\end{defn}
Notice that $R_s(MSp_{2n}(\mb R), \omega)$ and $R_{ss}(MO(p,q), \omega)$ only depend on $p+q$, not on a particular
pair $(p,q)$. Since
$$ \Re (v)-\bold n+ 2 \rho(O(p,q)) \prec \Re(v)-(\bold{n-\frac{p+q}{2}})+\rho(O(p,q)) $$
and
$$ \Re(v)-\bold {\frac{p+q}{2}}+ 2 \rho(Sp_{2n}(\mb R)) \prec \Re(v)-(\bold{\frac{p+q}{2}-n-1})+\rho(Sp_{2n}(\mb R)) ,$$
by the definition of semistable range, $R_{ss} \subseteq R_s$. Thus $\theta_s$ in the strongly semistable range is the same as the original $\theta$.
\begin{thm}[Corollary 5.2 \!~\cite{unit}]~\label{unitary1}
Suppose 
\begin{itemize}
\item $p+q \leq 2n+1$;
\item $\pi \in R_{ss}(MO(p,q), \omega)$;
\item $\pi$ is unitary;
\end{itemize}
Then $(\ , \ )_{\pi}$ is positive semidefinite.
If $\theta_s(p,q;2n)(\pi) \neq 0$, then $\theta_s(p,q;2n)(\pi)$ is an irreducible unitary representation of $MSp_{2n}(\mb R)$.
\end{thm}
\begin{thm}[Corollary 5.1 \!~\cite{unit}]~\label{unitary2}
Suppose
\begin{itemize}
\item $n < p \leq q$
\item $\pi \in R_{ss}(MSp_{2n}(\mb R), \omega)$;
\item $\pi$ is unitary;
\end{itemize}
Then $(\ , \ )_{\pi}$ is positive semidefinite.
If $\theta_s(2n;p,q)(\pi) \neq 0$, then $\theta_s(2n;p,q)(\pi)$ is an irreducible unitary representation of $MO(p,q)$.
\end{thm}
Strictly speaking, $\theta_s(*)(\pi)$ is an irreducible unitarizable Harish-Chandra module. Notice that $(\ , \ )_{\pi}$ can be regarded as an inner product on $\theta_s(*)(\pi)$. We can simply complete $(\ , \ )_{\pi}$ to obtain the Hilbert space of $\theta_s(*)(\pi)$. In general, invariant inner product on an irreducible Harish-Chandra module is unique up to the multiplication of a constant. \\
\\
Theorem \!~\ref{unitary1} and Theorem \!~\ref{unitary2} provide us the basis to construct unitary representation through the
theta correspondence in the strongly semistable range. 

\section{Quantum Induction $\mc Q$}
\begin{no}
Suppose that $p+q + p^{\p}+q^{\p} \equiv 0 \pmod 2$. Suppose that $p^{\p}+q \geq 2n+1$ and $q^{\p}+p \geq 2n+1$.
Let $\mathcal E(p^{\p}+q, q^{\p}+p;2n)$ be the representation $\pi_n$ of $O(p^{\p}+q,q^{\p}+p)$ studied by Zhu-Huang in \!~\cite{hz}. 
The representation $\mathcal E(p^{\p}+q, q^{\p}+p;2n)$ is essentially $\theta(2n;p^{\p}+q,q^{\p}+p)(\trivial)$ restricted to the $O(p^{\p}+q, q^{\p}+p)$ component in $MO(p^{\p}+q,q^{\p}+p)$. 
\end{no}
Notice here that 
$MO(p^{\p}+q,q^{\p}+p)$ splits into $O(p^{\p}+q,q^{\p}+p)$ and $O(p^{\p}+q,q^{\p}+p) \epsilon$.
\begin{defn} 
Suppose 
$$p^{\p}+q \geq 2n+1, \qquad p+q^{\p} \geq 2n+1, \qquad p+q \equiv p^{\p}+q^{\p} \ \pmod 2.$$
Let $O(p,q) \times O(p^{\prime}, q^{\prime})$ be  a subgroup diagonally embedded in $O(p^{\p}+q,q^{\p}+p)$. Fix a maximal compact subgroup $K=O(p^{\p}+q) \times O(q^{\prime}+p)$ such that
$$ K \cap O(p,q), \qquad K \cap O(p^{\prime}, q^{\prime})$$
are maximal compact subgroups of $O(p,q)$ and $O(p^{\prime}, q^{\prime})$ respectively.
 Let $\pi \in \Pi(O(p,q))$ be such that 
$$(u_1 \otimes v_1, u_2 \otimes v_2)=\int_{O(p,q)}(\mc E(p^{\p}+q,q^{\p}+p;2n)(g) u_1, u_2)
(\pi(g) v_1, v_2) d g  $$
converges absolutely for every $u_1, u_2 \in V(\mc E(p^{\p}+q,q^{\p}+p;2n))$ and every $v_1, v_2 \in V({\pi})$. Let $\mc R$ be the 
radical of $(\ , \ )$ as a Hermitian form on $ V(\mc E(p^{\p}+q,q^{\p}+p;2n)) \otimes V(\pi)$. Then 
$$(V(\mc E(p^{\p}+q,q^{\p}+p;2n)) \otimes V(\pi))/\mc R$$ inherits an infinitesimal $O(p^{\p}, q^{\p})$-action.
Define 
$$\mc Q(p,q;2n;p^{\p},q^{\p})(\pi)=V(\mc E(p^{\p}+q,q^{\p}+p;2n)) \otimes V({\pi})/ \mc R.$$
\end{defn}
$\mc Q(p,q;2n;p^{\p},q^{\p})(\pi)$ is a $(\f o(p^{\p}, q^{\p}), K \cap O(p^{\p}, q^{\p}))$-module. 
\begin{conj}\label{qconj1}
Suppose 
$$p^{\p}+q \geq 2n+1, \qquad p+q^{\p} \geq 2n+1, \qquad p+q \equiv p^{\p}+q^{\p} \ \pmod 2, \qquad p^{\prime}+q^{\prime} \geq p+q.$$
If $\mc Q(p,q;2n;p^{\p},q^{\p})(\pi) \neq 0$, then $\mc Q(p,q;2n;p^{\p},q^{\p})(\pi)$ is an irreducible admissible representation of
$O(p^{\p},q^{\p})$. If $\pi$ is unitary, then $\mc Q(p,q;2n;p^{\p},q^{\p})(\pi)$ is also unitary.
\end{conj}
\begin{no}
For $p+q \leq n+n^{\p}+1$, put $\mc E(2n+2n^{\p};p,q)=\theta(p,q;2n+2n^{\p})(\trivial)$.
\end{no}
If $p+q \leq n+n^{\p}$, $\mc E(2n+2n^{\p};p,q)$ is unitary according to 
Howe-Li's theory on stable range dual pairs (\!~\cite{howesmall}, \!~\cite{li2}). If $p+q =n+n^{\p}+1$, $\mc E(2n+2n^{\p};p,q)$ is unitary
according to the result of Przebinda on almost stable range dual pairs (see Lemma 8.6  \!~\cite{pr1}) or a Theorem of Kudla and Rallis (\!~\cite{kr}).
$\mc E(2n+2n^{\p};p,q)$ is a genuine representation of $Mp_{2n+2n^{\p}}(\mb R)$ 
if $p+q$ is odd.
\begin{defn} Suppose $n+n^{\p}+1 \geq p+q$. 
Let $Sp_{2n}(\mb R) \times  Sp_{2n^{\p}}(\mb R)$ be  diagonally  embedded into $Sp_{2n+2n^{\prime}}(\mb R)$. Let $K$ be a maximal compact subgroup of $Sp_{2n+2n^{\prime}}(\mb R)$ such that
$$Sp_{2n}(\mb R) \cap K , \qquad Sp_{2n^{\p}}(\mb R) \cap K$$
are maximal compact subgroups of $Sp_{2n}(\mb R)$ and $Sp_{2n^{\p}}(\mb R)$ respectively. 
Let $\pi$ be an irreducible representation of $Mp_{2n}(\mb R)$ such that the following 
Hermitian form $(\ , \ )$ on
$V(\mc E(2n+2n^{\p};p,q)) \otimes V({\pi})$ converges:
$$(\varphi \otimes u, \varsigma \otimes v)=
\int_{MSp_{2n}(\mb R)} (\mc E(2n+2n^{\p};p,q)(g) \varphi, \varsigma)( \pi(g) u,v) dg $$
$$\forall \ \  \varphi, \varsigma \in V(\mc E(2n+2n^{\p};p,q)); u,v \in V({\pi}).$$
Define $\mc Q(2n;p,q;2n^{\p})(\pi)$ to be $V(\mc E(2n;p,q;2n^{\p})) \otimes V({\pi})$ modulo the radical of $(\ , \ )$. $\mc Q(2n;p,q;2n^{\p})(\pi)$ inherits an
infinitesimal $\f{sp}_{2n^{\p}}(\mb R)$-action from $\mc E(2n+2n^{\p};p,q)$. It is a $({\f{sp}}_{2n^{\p}}(\mb R), MK \cap Mp_{2n^{\p}}(\mb R)))$-module.
\end{defn}
\begin{conj}\label{qconj2}  Suppose $n+n^{\p}+1 \geq p+q$ and $n^{\prime} \geq n$.
If $\mc Q(2n;p,q;2n^{\p})(\pi) \neq 0$, then $\mc Q(2n;p,q;2n^{\p})(\pi)$ is an irreducible admissible representation of
$Mp_{2n^{\p}}(\mb R)$. If $\pi$ is unitary, then $\mc Q(2n;p,q;2n^{\p})(\pi)$ is also unitary.
\end{conj}
Clearly, under the assumptions of Conjecture \ref{qconj1}, $\mc Q(p,q;2n;p^{\p},q^{\p})(\pi)$ defines an induction process from a subset of $\Pi(O(p,q))$ to $(\f g, K)$-modules of $O(p^{\p}, q^{\p})$.  Similarly, under the assumptions of Conjecture \ref{qconj2}, $\mc Q(2n;p,q;2n^{\p})$ defines an induction process from a subset of $\Pi(Mp_{2n}(\mb R))$ to $(\f g, K)$-modules of $Mp_{2n^{\p}}(\mb R)$.  It will be shown in the next section that these induction processes  produce irreducible unitary representations under some restrictions.   Hence we call $\mc Q(*)$ quantum induction. 
\section{Unitary Quantum Induction $Q$} 
$\mc Q(*)$ is closely related to a composition of two $\theta_s$.  The way it is defined allows us to bypass several technicalities which we shall address in this section. One difficult problem arising from composing two $\theta_s$ is that the range of one $\theta_s$ may not be in the domain of the other $\theta_s$.
In \!~\cite{basic}, we study the leading exponents of $\theta_s(*)(\pi)$. Under certain assumptions, one can 
compose $\theta_s$ with another $\theta_s$. Surprisingly, the strongly semistable range plays a crucial role in the composability of $\theta_s$. By identifying $\mc Q$ with a certain composition of $\theta_s$,
we established the unitarity and irreducibility of $\mc Q(\pi)$ for $\pi$ in the strongly semistable range.   Unitary quantum induction then enables us to construct certain irreducible unitary representations whose existence has been
conjectured by Arthur and Barbasch-Vogan. 

\begin{thm}[Theorem 7.1.1 \!~\cite{basic}]~\label{quantuminduction1}
Let $\pi$ be an irreducible unitary representation in $$R_{ss}(MO(p,q),\omega(p,q;2n)).$$ 
Suppose
\begin{itemize}
\item $q^{\p} \geq p^{\p} > n$;
\item $p^{\p}+q^{\p}-2n \geq 2n-(p+q)+2 \geq 1$;
\item $p+q \equiv p^{\p}+q^{\p} \qquad \pmod 2$. 
\end{itemize}
Then
\begin{enumerate}
\item $\theta_s(p,q;2n)(\pi) \in R_{ss}(MSp_{2n}(\mb R), \omega(p^{\p},q^{\p}; 2n))$;
\item  
 $\mc Q(p,q;2n;p^{\p},q^{\p})(\pi) \cong \theta_s(2n;p^{\p},q^{\p})\theta_s(p,q;2n) (\pi);$ 
\item If  $\mc Q(p,q;2n;p^{\p},q^{\p})(\pi) \neq 0$, then $\mc Q(p,q;2n;p^{\p},q^{\p})(\pi)$ is unitary.
\end{enumerate}
\end{thm}
Let me say a few words about the proof of this theorem. (1) is proved in \!~\cite{basic} through an estimate
on the matrix coefficients of $\theta_s(p,q;2n)(\pi)$. By (1),
$\theta_s(2n;p^{\p},q^{\p})\theta_s(p,q;2n) (\pi)$ is well-defined. In fact, by Definition \!~\ref{thetas},
the Harish-Chandra module of $\theta_s(p,q;2n)(\pi)$ consists of \lq\lq distributions\rq\rq of the following form
$$\int_{MO(p,q)} \omega(p,q;2n)(g_1) \phi \otimes \pi^c(g_1) u \, d g_1 \qquad (\phi \in V(\omega(p,q;2n)), u \in V(\pi)).$$
For the same reason, the Harish-Chandra module of $\theta_s(2n;p^{\p},q^{\p})\theta_s(p,q;2n) (\pi)$ consists of \lq\lq distributions\rq\rq of the following form
$$\int_{MSp_{2n}(\mathbb R)} \omega(p^{\p},q^{\p};2n)(g_2) \psi \otimes \theta_s(p,q;2n)(\pi)^c(g_2) v \, d g_2 \qquad
(\psi \in V(\omega(p^{\p},q^{\p};2n)), v \in V(\theta_s(p,q;2n)(\pi))).$$
Combining these two statements, $V(\theta_s(2n;p^{\p},q^{\p})\theta_s(p,q;2n) (\pi))$ consists of \lq\lq distributions\rq\rq of the form
\begin{equation}~\label{doubleintegral}
\int_{g_2 \in MSp_{2n}(\mathbb R)} \int_{g_1 \in MO(p,q)} \omega(p^{\p},q^{\p};2n)(g_2) \psi \otimes 
\omega(p,q;2n)^c(g_2 g_1) \phi \otimes \pi(g_1) u \, d g_1 \, d g_2 
\end{equation}
$$ (\psi \in V(\omega(p^{\p},q^{\p};2n)),\phi \in V(\omega(p,q;2n)), u \in V(\pi)).$$
Notice that $\omega(p,q;2n)^c \cong \omega(q,p;2n)$ and 
$$\omega(p^{\p},q^{\p};2n) \otimes \omega(q,p;2n) \cong
\omega(p^{\p}+q,q^{\p}+p;2n).$$
So Equation $(\!~\ref{doubleintegral})$ becomes
$$\int_{g_2 \in MSp_{2n}(\mathbb R)} \int_{g_1 \in MO(p,q)} \omega(p^{\p}+q,q^{\p}+p;2n)(i(g_1) g_2)(\psi \otimes \phi)  \otimes  \pi(g_1) u \, d g_1 \, d g_2 .$$
Here $i: g_1 \rightarrow i(g_1)$ is a canonical embedding of $MO(p,q)$ into $MO(p^{\prime}+q, q^{\prime}+p)$.
By the theorems of Howe-Li,
$$\{ \int_{g_2 \in MSp_{2n}(\mathbb R)} \omega(p^{\p}+q, q^{\p}+p;2n)(g_2) (\psi \otimes \phi) \, d g_2 \mid \psi \in V(\omega(p^{\p},q^{\p};2n)),\phi \in V(\omega(q,p;2n)) \}$$
can be identified with the Harish-Chandra module of $\mathcal E(p^{\p}+q,q^{\p}+p;2n)$.
Therefore, $$V(\theta_s(2n;p^{\p},q^{\p})\theta_s(p,q;2n) (\pi))$$
consists of $$
\int_{MO(p,q)} \mathcal E(p^{\p}+q,q^{\p}+p;2n)(i(g_1)) \eta \otimes \pi(g_1) u \, d g_1 \qquad ( \eta \in V(\mathcal E(p^{\p}+q,q^{\p}+p;2n)), u \in V(\pi)) .$$
We have assumed that $\pi(\epsilon)=-1$ from the beginning. The integral over $MO(p,q)$ is just twice the integral over $O(p,q)$. It follows that $V(\theta_s(2n;p^{\p},q^{\p})\theta_s(p,q;2n) (\pi))$ consists of vectors of the following form
 $$
 \int_{O(p,q)} \mathcal E(p^{\p}+q,q^{\p}+p;2n)(i(g_1)) \eta \otimes \pi(g_1) u \, d g_1 \qquad ( \eta \in V(\mathcal E(p^{\p}+q,q^{\p}+p;2n)), u \in V(\pi)). $$
 The absolute convergences proved in \!~\cite{basic} allow us to  change the order of the integrals. Thus
$$  \theta_s(2n;p^{\p},q^{\p})\theta_s(p,q;2n) (\pi) \cong \mc Q(p,q;2n;p^{\p},q^{\p})(\pi) . $$
(3) follows easily by Theorems \!~\ref{unitary1} and \!~\ref{unitary2}. $\Box$
\begin{defn}~\label{quantum0} Under the hypotheses in Theorem \!~\ref{quantuminduction1}, define
$$Q(p,q;2n;p^{\p},q^{\p})(\pi)= \theta_s(2n;p^{\p},q^{\p})\theta_s(p,q;2n)(\pi).$$
\end{defn}
We call $Q$ {\it unitary quantum induction}.
By Theorem \!~\ref{quantuminduction1}, $Q(p,q;2n;p^{\p},q^{\p})(\pi)$ is well-defined and equivalent to
$\mc Q(p,q;2n;p^{\p},q^{\p})(\pi)$. A similar statement holds for $\mc Q(2n;p,q;2n^{\p})(\pi)$.
\begin{thm}[Theorem 7.2.1 \!~\cite{basic}]~\label{quantuminduction2} 
Let $\pi$ be a unitary representation in $R_{ss}(MSp_{2n}(\mb R), \omega(p,q;2n))$. Suppose
\begin{itemize}
\item $2n^{\p}-p-q \geq p+q-2n-2$;
\item $n < p \leq q$.
\end{itemize}
Then
\begin{enumerate}
\item $\theta_s(2n;p,q)(\pi) \in R_{ss}(MO(p,q), \omega(p,q;2n^{\p}));$
\item $\mc Q(2n;p,q;2n^{\p})(\pi) \cong \theta_s(p,q;2n^{\p})\theta_s(2n;p,q) (\pi);$
\item
if $\mc Q(2n;p,q;2n^{\p})(\pi) \neq 0$, then
$\mc Q(2n;p,q;2n^{\p})(\pi)$ is unitary. 
\end{enumerate}
\end{thm}
\begin{defn}~\label{quantum} Under the hypotheses of Theorem \!~\ref{quantuminduction2}, define
$$Q(2n;p,q;2n^{\p})(\pi)=\theta_s(p,q;2n^{\p})\theta_s(2n;p,q)(\pi).$$
\end{defn}
By Theorem \!~\ref{quantuminduction2}, $ Q(2n;p,q;2n^{\p})(\pi)$ is well-defined and equivalent to $\mc Q(2n;p,q;2n^{\p})(\pi)$. \\
\\
In summary, $\mc Q$ is a generalization of $Q= \theta_s \circ \theta_s$. Without the requirement of unitarity, we can loose some restrictions posed in Theorems \!~\ref{quantuminduction1} and \!~\ref{quantuminduction2}. We only state the Theorem for $\theta_s(p,q;2n^{\p})\theta_s(2n;p,q)$ since we will need it in the next Chapter.
\begin{thm}~\label{quantuminduction3}
Let $\pi$ be a unitary representation in $R_{ss}(MSp_{2n}(\mb R), \omega(p,q;2n))$. Suppose
$2n^{\p}-p-q \geq p+q-2n-2$.
Then
$\theta_s(2n;p,q)(\pi) \in R_{ss}(MO(p,q), \omega(p,q;2n^{\p}));$
 and $$\mc Q(2n;p,q;2n^{\p})(\pi) \cong \theta_s(p,q;2n^{\p})\theta_s(2n;p,q) (\pi).$$
\end{thm}
This theorem follows directly from the estimate established in Theorem 6.3.1 and  Definition \!~\ref{qi} of \!~\cite{basic}. The proof is omitted. \\
\\
The advantage of introducing $\mc Q$ will be manifested when we begin to relate  quantum induction to parabolic induction. Based on the studies on
parabolic induction, we will establish the nonvanishing of $\theta_s(2n;p,q)(\pi)$ under some restrictions.
The nonvanishing of $\theta_s(2n;p,q)(\pi)$ is otherwise very hard to establish if $p$ or $q$ is less than $2n$.

\section{Moment Map and $\Theta(G_1, G_2)$}
From now on, let $G_1=O(p,q)$ and $G_2=Sp_{2n}(\mb R)$. 
Write $${G_1}_{\mb C}=O(p+q,\mb C), \qquad {G_2}_{\mb C}=Sp_{2n}(\mb C).$$
Let $\Mat(p+q;2n)$ be the set of $p+q$ by $2n$ real matrices. Let 
$$W_n=\arr{0 & I_n \\ -I_n & 0}, \qquad I_{p,q}= \arr{ I_p & 0 \\ 0 & -I_q}.$$
\begin{defn}~\label{moment}
Define the moment maps
$$m_1: x \in \Mat(p+q, 2n) \rightarrow I_{p,q} x W_n x^t \in \f o(p,q),$$
$$m_2: x \in \Mat(p+q,2n) \rightarrow W_n x^t I_{p,q} x \in {{\f{sp}}}_{2n}(\mb R).$$
We say that $x$ is nilpotent if
$m_1(x)$ is nilpotent.
\end{defn}
Observe that $m_1(x)^l=I_{p,q} x m_2(x)^{l-1} W_n x^t$. We have
\begin{lem}
$m_1(x)$ is nilpotent if and only if $m_2(x)$
is nilpotent. 
\end{lem}
The group $Sp_{2n}(\mb R) \times O(p,q)$ acts on the set of nilpotent element in $\Mat(p+q,2n)$.
Thus we obtain a set of nilpotent $Sp_{2n}(\mb R) \times O(p,q)$-orbits
in $\Mat(p+q, 2n)$. Clearly, $m_i$ induces a
map from the nilpotent orbits in $\Mat(p+q, 2n)$ to nilpotent orbits in
$\f g_i$. \\
\\
Similarly, we define the complex moment maps, still denoted by $m_i$. This should not cause any notational problem. For example, the usage of $\bold D$ usually points to the real moment map and the usage of $\bold d$ usually points to the complex moment map.\\
\\
For a complex nilpotent orbit $\mc O_1$ in $\f {g_1}_{\mb C}$, consider
the closure of
$$m_2(m_1^{-1} (\mbox{Closure of} \ \mc O_1)).$$
By a Theorem of Daszkiewicz-Kraskiewicz-Przebinda, the closure of
$$m_2(m_1^{-1} (\mbox{Closure of} \ \mc O_1))$$
is the closure of a unique nilpotent orbit in $\f {g_2}_{\mb C}$ (see \!~\cite{dkp}). This yields a correspondence between certain nilpotent orbits in ${\f g_1}_{\mb C}$ and certain nilpotent orbits in ${\f g_2}_{\mb C}$. We denote this correspondence by
$\Theta({G_1}_{\mb C}, {G_2}_{\mb C})$. Similarly, we define $\Theta({G_2}_{\mb C}, {G_1}_{\mb C})$. The reader
should be warned that $\Theta({G_1}_{\mb C}, {G_2}_{\mb C})$ is NOT the inverse of $\Theta({G_2}_{\mb C}, {G_1}_{\mb C})$.
For the nilpotent orbits we are concerned with, $\Theta({G_1}_{\mb C}, {G_2}_{\mb C})$ and $\Theta({G_2}_{\mb C}, {G_1}_{\mb C})$ are quite easy to describe (\!~\cite{dkp}). 
\begin{lem}[see Thm. 4.2 \!~\cite{dkp}]~\label{ass0}
Let $\mc O_{\bold D} \in \mc U$ be a nilpotent orbit of either $O(p,q)$ or $Sp_{2n}(\mb R)$ (see Definition ~\ref{u}).
Construct an alternating sequence of complex orthogonal orbits and complex symplectic orbits
$$\mc O_{\bold d}, \mc O_{\bold d-\bold 1}, \mc O_{\bold d-\bold 2}, \ldots \mc O_{\bold d-\bold{d_1}+\bold 1}$$
and a corresponding alternating sequence of complex orthogonal groups and symplectic groups:
$$G(\mc O_{\bold d}), G(\mc O_{\bold d- \bold 1}), \ldots, G(\mc O_{\bold d-\bold{d_1}+\bold 1}).$$
Then $\forall \ j$,
$$\Theta(G(\mc O_{\bold d-\bold j}), G(\mc O_{\bold d-\bold j+\bold 1}))(cl(\mc O_{\bold d-\bold j}))=
cl(\mc O_{\bold d-\bold j+\bold 1}).$$
\end{lem}
We define $\Theta({G_2}, {G_1})$ and $\Theta(G_1, G_2)$ for real nilpotent orbits in the same fashion. Pan showed that
$\Theta(G_2, G_1)(\mc O_{\bold D})$ is the closure of at most two real
nilpotent orbits (see Theorem 8.11 \!~\cite{pan}). We give two lemmas that can be easily deduced from the descriptions of $\Theta(G_2, G_1)$ and $\Theta(G_1, G_2)$ in \!~\cite{pan}.
\begin{lem}[see 6.4, 8.6 \!~\cite{pan}]~\label{panlemma}
Suppose $p+q \leq 2n$. Then the real nilpotent orbit $\mc O_{\bold D}$ occurs in the image of $m_2$ if and only if
$G(\mc O_{\bold D-1})=O(r,s)$ with $r \leq p$ and $s \leq q$.
\end{lem}
\begin{lem}[see 6.4, 8.1, 8.11 \!~\cite{pan}] ~\label{panlemma0}
Let $\mc O_{\bold D} \in \mc U$ be a real nilpotent orbit of either $O(p,q)$ or $Sp_{2n}(\mb R)$.
Construct an alternating sequence of real nilpotent orbits of symplectic groups and orthogonal groups
$$\mc O_{\bold D}, \mc O_{\bold {D-1}}, \mc O_{\bold {D-2}}, \ldots \mc O_{\bold {D-d_1+1}}$$
and a corresponding sequence of real groups
$$G(\mc O_{\bold D}), G(\mc O_{\bold {D-1}}), G(\mc O_{\bold {D-2}}), \ldots G(\mc O_{\bold {D-d_1+1}}).$$
Then $\forall \ j$,
$$\Theta(G(\mc O_{\bold D-\bold j}), G(\mc O_{\bold{D-j+1}}))(cl(\mc O_{\bold {D-j}}))=
cl(\mc O_{\bold {D-j+1}}).$$
\end{lem}

\section{Przebinda's Results on $\mc V(Ann\;\theta(\pi))$}
In \!~\cite{pr1}, Przebinda proved
that under some conditions, the following diagram commutes:
\begin{equation}
\begin{CD}
\pi @>{\theta} >> \theta(MG, MG^{\p})(\pi) \\
@VV{}V @VV{}V \\
\mc V(Ann\;\pi)  @>{\Theta(G_{\mb C}, G^{\p}_{\mb C})}>> \mc V(Ann\;\theta(MG, MG^{\p})(\pi))
\end{CD}
\end{equation}
(See Notation 5). 
We state Przebinda's theorem for $O(p,q)$ and $Mp_{2n}(\mathbb R)$. Recall from Theorem III.1 \!~\cite{mi} and Theorem 8.48 of \!~\cite{knapp} that
$\pi$ has the rate of growth $\gamma$ if and only if for every leading exponent $v$ of $\pi$,
$$\Re(v) \preceq (\gamma-1) \rho(G)$$
(see also Theorem 4.5 of \!~\cite{pr1}). 
\begin{thm}[Theorem 0.9 \!~\cite{pr1}]~\label{ass1}
Let
$\pi$ be an irreducible unitary representation of $O(p,q)$. Suppose 
\begin{enumerate}
\item $2n+1 \geq p+q$;
\item $\Theta_{\pi}$ has rate of growth $\gamma$ with
$\gamma+1 < \frac{2n}{p+q-2}$,
equivalently, every leading exponent $v$ of $\pi$ satisfies
$$\Re(v) \prec (\frac{2n}{p+q-2}-2) \rho(O(p,q));$$
\item $\gamma \geq 0$, equivalently, $\gamma$ is not a discrete series representation;
\item $(\ , \ )_{\pi}$ converges and does not vanish;
\item $(\ , \ )_{\pi}$ is positive semidefinite;
\item there exists a full rank element $x \in \Mat(p+q,2n)$ such that
the $O(p,q)$-orbit generated by $m_1(x)$ is of maximal dimension in $AS(\Theta_{\pi})$.
\end{enumerate}
Then $\theta(p,q;2n)(\pi)$ is unitary and its associated variety is the complex orbit 
$$\Theta(O(p+q, \mb C), Sp_{2n}(\mb C))(\mc V(Ann\;\pi)).$$
\end{thm}
Let me make two remarks here. Firstly, the assumption (1) does not appear in Przebinda's original theorem.   Conditions (2) and (3) imply $2n > p+q-2$. We add (1) as an assumption for the sake of clarity. The reader should also notice that condition (1) guarantees the nonvanishing of $(\ , \ )_{\pi}$ or
$(\ , \ )_{\det \otimes \pi}$ by Theorem \!~\ref{nonvanishing}. Secondly, condition (6) is satisfied in the setting of Lemma \!~\ref{panlemma0}. More precisely, 
\begin{lem}~\label{50}
Suppose that $\pi$ is an irreducible unitary representation of $O(p,q)$, that $2n+1 \ge p+q$, and that $\bold D$ is a signed Young tableau for $Sp_{2n}(\mb R)$. Assume that
$\mc O_{\bold D-\bold 1}$ is of maximal dimension in $AS(\Theta_{\pi})$. Then there is an element $x \in \Mat(p+q, 2n)$ of rank $p+q$ with the properties
\begin{enumerate}
\item $O(p,q) m_1(x)= \mc O_{\bold D-\bold 1}$;
\item $Sp_{2n}(\mb R) m_2(x)=\mc O_{\bold D}$.
\end{enumerate}
\end{lem}
Proof: According to Lemma \!~\ref{panlemma}, $\mc O_{\bold D}$ occurs in the image of $m_2$. By Lemma \!~\ref{panlemma0}, there exists $x$  in $\Mat(p+q,2n)$ such that $m_2(x) \in \mc O_{\bold D}$ and $m_1(x) \in cl(\mc O_{\bold D-1})$. Hence $$rank( W_n x^t I_{p,q} x)=rank(m_2(x))=p+q.$$
But $x \in \Mat(p+q,2n)$. So $x$ must be of full rank. For $x$ of full rank, observe that
$$[m_2(x)]^r= W_n x^t [m_1(x)]^{r-1} I_{p,q} x.$$
By writing out the defining equations for $m_2(x) \in \mc O_{\bold D}$, we obtain
$m_1(x) \in \mc O_{\bold {d-1}}$. 
But $m_1(x) \in cl(\mc O_{\bold D-1})$. It follows that
$m_1(x) \in \mc O_{\bold {D-1}}$.  $\Box$\\
\\
Similarly, Przebinda's theorem for $Mp_{2n}(\mb R)$ can be formulated as follows.
\begin{thm}[Theorem 0.9 \!~\cite{pr1}]~\label{ass2}
Let
$\pi$ be an irreducible unitary representation of $Mp_{2n}(\mb R)$. Suppose
\begin{enumerate}
\item $p+q >2n$;
\item $\Theta_{\pi}$ has the rate of growth $\gamma$ with
$\gamma+1 < \frac{p+q}{2n}$,
equivalently, every leading exponent $v$ of $\pi$ satisfies
$$\Re(v) \prec (\frac{p+q}{2n}-2) \rho(Mp_{2n}(\mb R));$$
\item $\gamma \geq 0$, equivalently, $\pi$ is not a discrete series representation;
\item $(\ , \ )_{\pi}$ converges and does not vanish;
\item $(\ , \ )_{\pi}$ is positive semidefinite;
\item there exists a full rank element $x \in \Mat(p+q,2n)$ such that
the $Sp_{2n}(\mb R)$-orbit generated from $m_2(x)$ is of maximal dimension in $AS(\Theta_{\pi})$.
\end{enumerate}
Then $\theta(2n;p,q)(\pi)$ is unitary and its associated variety is the complex nilpotent orbit 
$$\Theta(Sp_{2n}(\mb C), O(p+q,\mb C))(\mc V(Ann\;\pi)).$$
\end{thm}
It is easy to see that (2) and (3) imply (1). We leave (1) there for the sake of clarity.
Similar to Lemma \!~\ref{50}, we have
\begin{lem}~\label{51}
Suppose that $\pi$ is an irreducible unitary representation of $Mp_{2n}(\mb R)$, that $ p+q> 2n$, and that $\bold D$ is a signed Young tableau for $O(p,q)$. Assume that
$\mc O_{\bold D-\bold 1}$ is of maximal dimension in $AS(\Theta_{\pi})$. Then there is an element $x \in \Mat(p+q, 2n)$ of rank $2n$ with the properties
\begin{enumerate}
\item $O(p,q) m_1(x)= \mc O_{\bold D}$;
\item $Sp_{2n}(\mb R) m_2(x)=\mc O_{\bold D-\bold 1}$.
\end{enumerate}
\end{lem} 
So $(6)$ in Theorem \!~\ref{ass2} holds if there is an orbit in $AS(\Theta_{\pi})$ of the form $\mc O_{\bold D-\bold 1}$ for $\bold D \in \mc YD(p,q)$.\\
\\
Last let us recall the following result of Przebinda.
\begin{thm}[Corollary 2.8, \!~\cite{pr1}]~\label{pr2.8}
Consider the dual pair $(O(p,q), Sp_{2n}(\mb R))$.
\begin{itemize}
\item
If $\pi \in \mc R(MO(p,q), \omega)$, then 
$WF(\pi)$ is a subset of $m_1(\Mat(p+q;2n)).$
\item If $\pi \in \mc R(MSp_{2n}(\mb R), \omega)$, then
$WF(\pi)$ is a subset of $m_2(\Mat(p+q;2n))$.
\end{itemize}
\end{thm}

\chapter{A Nonvanishing Theorem}\label{nonvan}
In \!~\cite{non}, we established a nonvanishing theorem for $\theta_s(p,q;2n)$ with
$p+q \leq 2n+1$. This theorem is cited as Theorem \!~\ref{nonvanishing} in the last Chapter.
In this chapter, we will prove a nonvanishing theorem for $\theta_s(2n;p,q)$.
These two nonvanishing theorems will then be used to construct a packet of irreducible unitary representations $\mc N(\mc O)$ attached to a real special rigid orbit $\mc O$ and more generally any $\mc O \in \mc U$ (see Definition \!~\ref{u}).
We start with a trivial lemma.
\begin{lem}~\label{triv}
Suppose that 
$$\theta_s(p,q;2n_2)\theta_s(2n_1;p,q)(\pi)=\mc Q(2n_1;p,q;2n_2)(\pi).$$
If $\mc Q(2n_1;p,q;2n_2)(\pi) \neq 0$, then
$\theta_s(2n_1;p,q)(\pi) \neq 0$.
\end{lem}
Combining with Theorem \!~\ref{quantuminduction3}, we obtain
\begin{thm}~\label{notzero}
 Suppose $2n_2-p-q \geq p+q-2n_1-2$.
Let $\pi$ be an irreducible admissible representation in $R_{ss}(MSp_{2n_1}(\mb R), \omega(p,q;2n_1))$ ( see Def. \!~\ref{qi}).
If $\mc Q(2n_1;p,q;2n_2)(\pi) \neq 0,$ then
$\theta_s(2n_1;p,q)(\pi) \neq 0$.
\end{thm}
At first glance, $\mc Q(2n_1;p,q;2n_2)(\pi)$ seems to be more difficult to treat than $\theta_s(2n_1;p,q)$. This is true with one exception. As predicted in \!~\cite{basic}, when $p+q=n_1+n_2+1$, quantum induction can be obtained from parabolic induction. In this chapter, we will examine $\mc Q(2n_1;p,q;2n_2)(\pi)$ for $n_1+n_2+1=p+q$.
We will prove that precisely for  $n_1+n_2+1=p+q$ and $n_1 \leq n_2$, $\mc Q(2n_1;p,q;2n_2)(\pi)$ can be obtained from parabolically induced representations. As a consequence, we can determine the nonvanishing of $\mc Q(2n_1;p,q;2n_2)(\pi)$.
\begin{no}~\label{maxparabolic} Let $m, n$ be two nonnegative integers such that $m+n \neq 0$. Let $P_{m,n}$ be the standard maximal parabolic subgroup of $Sp_{2m+2n}(\mb R)$ with Levi subgroup $GL(m,\mathbb R) \times Sp_{2n}(\mathbb R)$.
\end{no}
\section{Nonvanishing of $\mc Q(2n_1;p,q;2(p+q-n_1-1))(\pi)$: Overview and Main Result}
Suppose that $n_2 \geq n_1$. 
Let $Sp_{2n_1}(\mb R) \times Sp_{2n_2}(\mb R)$ be a subgroup of $Sp_{2n_1+2n_2}(\mb R)$. Let $U(n_1+n_2)$ be a maximal compact subgroup of $Sp_{2n_1+2n_2}(\mb R)$ such at $U(n_1+n_2) \cap Sp_{2n_1}(\mb R)$ and $U(n_1+n_2) \cap
Sp_{2n_2}(\mb R)$ are maximal compact subgroups of $Sp_{2n_1}(\mb R)$ and $Sp_{2n_2}(\mb R)$ respectively.\\
\\
First of all, a theorem of Kudla and Rallis says that, for some $\alpha$, the parabolic induced representation 
$$I^{\alpha}=\Ind_{MP_{n_1+n_2,0}}^{Mp_{2n_1+2n_2}(\mb R)} \chi^{\alpha} $$
(see Definition ~\ref{inalpha})
decomposes into a direct sum of irreducible representations, namely $\mc E(2n_1+2n_2; p,q)$ for some $p$ and $q$ (\!~\cite{kr}, \!~\cite{lz}).  Regard $I^{\alpha}$ as a representation of $Mp_{2n_1}(\mb R) \times Mp_{2n_2}(\mb R)$ by restriction. For  $\pi \in \Pi(Mp_{2n_1}(\mb R))$ satisfying a certain growth condition, $I^{\alpha}$ induces a representation of $Mp_{2n_2}(\mb R)$, $I^{\alpha}(\pi)$ (see Definition \!~\ref{ipi}). Using the theorem of Kudla and Rallis, we show that  $I^{\alpha}(\pi)$ decomposes into a direct sum of $\mc Q(2n_1;p,q;2n_2)(\pi)$ for some $p+q=n_1+n_2+1$. \\
\\
Next, we study the restriction of $I^{\alpha}$ to $Mp_{2n_1}(\mb R) \times Mp_{2n_2}(\mb R)$. Notice that the Harish-Chandra module $V(I^{\alpha})$ consists of sections of a homogeneous line bundle over the Lagrangian Grassmannian of
$\mathbb R^{2n_1+2n_2}$. By analyzing the action of $Mp_{2n_1}(\mathbb R) \times Mp_{2n_2}(\mathbb R)$ on the Grassmanian, we prove that
$$I^{\alpha}(\pi) \cong \Ind_{ MP_{n_2-n_1, n_1}}^{Mp_{2n_2}(\mathbb R)} \pi^{\tau} \otimes \chi^{\alpha},$$
whenever the left hand side is well-defined (see Theorem \!~\ref{induction1}). In order that $I^{\alpha}(\pi)$ is well-defined, $\pi$ must satisfy a certain growth condition. The details are given in Section \!~\ref{ialphainduced}.\\
\\
Now, let us assume that $I^{\alpha}(\pi)$ is well-defined. Then we know that
$$ \Ind_{ MP_{n_2-n_1, n_1}}^{Mp_{2n_2}(\mathbb R)} \pi^{\tau} \otimes \chi^{\alpha} \simeq \bigoplus \mc Q(2n_1;p,q;2n_2)(\pi)$$
for a set of pairs $(p,q)$.  The remaining question is 
which $\mc Q(2n_1;p,q;2n_2)(\pi)$ is nonvanishing. This turns out to be a hard question. We don't have any answer in general. But for some $\pi$, we can detect the nonvanishing of $\mc Q(2n_1;p,q;2n_2)(\pi)$ precisely, thanks to the notion of wave front sets. 
Notice that the decomposition of $I^{\alpha}(\pi)$ into $\mc Q(*)(\pi)$ induces a decomposition on the wave front level, namely,
\begin{equation}~\label{twowf}
WF(\Ind_{ MP_{n_2-n_1, n_1}}^{Mp_{2n_2}(\mathbb R)} \pi^{\tau} \otimes \chi^{\alpha})= \bigcup WF(\mc Q(2n_1;p,q; 2n_2)(\pi)).
\end{equation}
 On the one hand, we know that 
$$WF(\Ind_{ MP_{n_2-n_1, n_1}}^{Mp_{2n_2}(\mathbb R)} \pi^{\tau} \otimes \chi^{\alpha})=\Ind^{{\f{sp}}_{2n_2}(\mathbb R)}_{{\f{gl}}(n_2-n_1)(\mb R) \times {\f{sp}}_{2n_1}(\mb R)} WF(\pi^{\tau}).$$
As we have pointed out earlier, the induced orbit consists of a few irreducible components which can be parametrized easily. For the $\pi$'s we are concerned with, Theorems \!~\ref{induced1} and \!~\ref{tauinduction} provide us with the precise information about the nilpotent orbits contained in the left hand side of Equation \!~\ref{twowf}.\\
\\
On the other hand, for some $\pi$, for example, those specified in Theorem \!~\ref{quantuminduction3}, the right hand side of Equation \!~\ref{twowf} equals 
$$\bigcup WF(\theta_s(p,q; 2n_2) \theta_s(2n_1;p,q)(\pi)).$$
Theorem \!~\ref{pr2.8} of Przebinda gives some control on which nilpotent orbit occurs in the wave front set
of $\theta(p,q;2n_2)(*)$ (\!~\cite{pr1}). Lemma \!~\ref{panlemma} of Pan gives the parametrization of these nilpotent orbits (\!~\cite{pan}). \\
\\
We have now two sets of nilpotent orbits, one from the induced orbit, the other from the moment map
$m_2$ with respect to various $(O(p,q), Sp_{2n_2}(\mb R))$. By matching these two sets of wave front sets,  we obtain a nonvanishing theorem for 
$WF(\mc Q(2n_1;p,q; 2n_2))$, which in turn gives us  the nonvanishing of $\mc Q(2n_1;p,q;2n_2)$ and  $\theta_s(2n_1;p,q)$. Some of the results here will
be used in Chapter 6. In this Chapter, we will prove the following generic nonvanishing theorem.

\begin{thm}~\label{non2}
Consider the group $Mp_{2n_1+2n_2}(\mb R)$ with $n_1 < n_2$.
Let $\pi$ be a unitary representation in $R_{ss}(Mp_{2n_1}(\mb R), \omega(n_1+n_2+1,0; 2n_1))$
(see Def. \!~\ref{qi}).
Let $\mc O_{\bold D}$ be a nilpotent orbit of maximal dimension in 
$$\Ind_{\f{sp}_{2n_1}(\mb R) \oplus \f{gl}(n_2-n_1, \mb R)}^{\f{sp}_{2n_2}(\mb R)} WF(\pi^{\tau}).$$ 
Let $S$ be the subset of $ \{(p, q) \mid  p+q=n_1+n_2+1, p \ fixed \ parity \},$ consisting of those $(p,q)$ such that
$\mc O_{\bold D}$ is contained
the image of the moment map $m_2$ associated with $(O(p,q), Sp_{2n_2}(\mb R))$.
Then there exists a pair $(p,q) \in S$ such that
$\mc Q(2n_1; p,q; 2n_2)(\pi) \neq 0$. For such a pair,
$\theta_s(2n_1;p,q)(\pi) \neq 0$.
\end{thm}
Notice that
$$R_{ss}(Mp_{2n_1}(\mb R),\omega(n_1+n_2+1,0;2n_1))=R_{ss}(Mp_{2n_1}(\mb R), \omega(p,q;2n_1))$$
(see Def. \!~\ref{qi}) for every $p+q=n_1+n_2+1$. The most favorable situation is when the set $S$, which depends on $\bold D$, consists of only one pair $(p,q)$. By our theorem, we will have $\theta_s(2n_1;p,q)(\pi) \neq 0$. The proof of Theorem \!~\ref{non2} will be given in Ch. 4.5.\\
\\
In this chapter, we take $\theta$ and $\theta_s$ as a correspondence
between representations of $Mp_{2n}(\mb R)$ and representations of $O(p,q)$.
Throughout this Chapter, we assume $n_2 \geq n_1$.

\section{Results of Kudla-Rallis and Lee-Zhu}
Consider the symplectic group $Sp_{2n}(\mb R)$. Let
$MAN$ be the Langlands decomposition of $P_{n,0}$ ( Notation \!~\ref{maxparabolic}). The Levi factor 
$MA \cong GL(n,\mathbb R)$. Let $Mp_{2n}(\mb R)$ be the metaplectic covering
of $Sp_{2n}(\mb R)$. The group $MGL(n, \mb R)$ can be expressed as
$$\{(\xi, x) \mid x \in GL(n, \mb R), \xi^2= \det x, \xi \in \mathbb C\}.$$
Put $\chi(\xi, x)= \frac{\xi}{|\xi|}$. Then $\chi$ is a character  of $MGL(n, \mb R)$ of order 4. We extend $\chi$ trivially
on $N$. $\chi$ becomes a character of $MP_{n,0}$.
\begin{defn}~\label{inalpha}
$I_n^{\alpha}=\Ind_{MP_{n,0}}^{Mp_{2n}(\mb R)} \chi^{\alpha}$. 
\end{defn}
$I_n^{\alpha}$ is a unitarily
induced representation.

\begin{thm}[Kudla-Rallis, \!~\cite{kr}]
Suppose $p+q=n+1$ and $\alpha \equiv p-q \pmod 4$. 
Then $\theta(p,q;2n)(\trivial)$ is a subrepresentation of $I_n^{\alpha}$.
\end{thm}
By studying the $K$-types in $\theta(p,q;2n)(\trivial)$ for various $p+q=n+1$, one has
\begin{thm}[\!~\cite{kr}, \!~\cite{lz}]~\label{krlz}
$$I_n^{\alpha}= \bigoplus_{p+q=n+1, \alpha \equiv p-q \pmod 4} \theta(p,q;2n)(\trivial).$$
More precisely, for $n$ odd,
$$I_n^0=\Ind_{P_{n,0}}^{Sp_{2n}(\mb R)} \trivial =\bigoplus_{p+q=n+1, p \equiv q \pmod 4} \theta(p,q;2n)(\trivial),$$
$$I_n^2=\Ind_{P_{n,0}}^{Sp_{2n}(\mb R)} \sgn(\det)=\bigoplus_{p+q=n+1, p-q \equiv 2 \pmod 4} \theta(p,q;2n)(\trivial).$$
For $n$ even,
$$I_n^1=\bigoplus_{p+q=n+1, 1 \equiv p-q \pmod 4} \theta(p,q;2n)(\trivial),$$
$$I_n^3=\bigoplus_{p+q=n+1, 3 \equiv p-q \pmod 4} \theta(p,q;2n)(\trivial).$$
\end{thm}
This decomposition theorem can be derived directly from the results in \!~\cite{kr}. It was explicitly stated in \!~\cite{lz}. 
\begin{no} We use
$I^{\alpha}$ to denote $I^{\alpha}_{n_1+n_2}$ which is the main object of study in this section.
\end{no}
\begin{defn}~\label{ipi}
Let $Sp_{2n_1}(\mb R) \times Sp_{2n_2}(\mb R)$ be embedded diagonally into $Sp_{2n_1+2n_2}(\mb R)$. Let $U(n_1+n_2)$ be a maximal compact subgroup of $Sp_{2n_1+2n_2}(\mb R)$ such at $U(n_1+n_2) \cap Sp_{2n_1}(\mb R)$ and $U(n_1+n_2) \cap
Sp_{2n_2}(\mb R)$ are maximal compact subgroups of $Sp_{2n_1}(\mb R)$ and $Sp_{2n_2}(\mb R)$ respectively. The embedding $Sp_{2n_1}(\mb R) \rightarrow Sp_{2n_1+2n_2}(\mb R)$ lifts to
$Mp_{2n_1}(\mb R) \rightarrow Mp_{2n_1+2n_2}(\mb R)$ (see for example \!~\cite{unit}). Consider $I^{\alpha}$. Let $\pi \in \Pi(Mp_{2n_1}(\mb R))$ such that the form $(\ , \ )_1$  
$$( \phi_1 \otimes u_1,  \phi_2 \otimes u_2)_1=\int_{Mp_{2n_1}(\mb R)}(\pi(g) u_1, u_2)(I^{\alpha}(g) \phi_1, \phi_2) \, d g $$
converges for every $u_1, u_2 \in V({\pi}), \phi_1, \phi_2 \in V({I^{\alpha}})$. Define $I^{\alpha}(\pi)$  to be $ V({I^{\alpha}}) \otimes V(\pi)$ modulo the radical of $(\ , \ )_1$. $I^{\alpha}(\pi)$ is a $({\f sp}_{2n_2}(\mb R), MU(n_1+n_2) \cap Mp_{2n_2}(\mb R))$-module.
\end{defn}
From the definition of quantum induction $\mc Q$ and Theorem \!~\ref{krlz}, we obtain
\begin{thm}~\label{induction}
Suppose $\alpha \equiv n_1+n_2+1 \pmod 2$. Then as $(\f g, K)$-modules
$$I^{\alpha}(\pi)=\bigoplus_{p+q=n_1+n_2+1, p-q \equiv \alpha \pmod 4} \mc Q(2n_1;p,q;2n_2)(\pi)$$
whenever one side is well-defined.
\end{thm}
Let me remind the reader that
$\theta(p,q;2n)(\trivial)$ is denoted by $\mc E(2n;p,q)$ in the definition of $\mc Q$.\\
\\
Suppose that $\pi \in \Pi_u(Mp_{2n_1}(\mb R))$. It is known, at least under the assumptions of Theorem \!~\ref{quantuminduction2}, that $\mathcal Q(2n_1;p,q;2n_2)(\pi)$ is a unitarizable Harish-Chandra module of $Mp_{2n_2}(\mathbb R)$. Theorem \!~\ref{induction} then implies that $I^{\alpha}(\pi)$ is also a unitarizable Harish-Chandra module in this case. Later, in Ch. 4.4, we will show directly that $I^{\alpha}(\pi)$ is a unitarizable Harish-Chandra module of $Mp_{2n_2}(\mathbb R)$ in a much more general setting.

\section{$Sp_{2n_1}(\mb R) \times Sp_{2n_2}(\mb R)$ action on Lagrangian Grassmannian $X_{n_1+n_2}$}
We are interested in the $Mp_{2n_2}(\mb R) \times Mp_{2n_1}(\mb R)$-action on 
$$I^{\alpha}=\Ind_{MP_{n_1+n_2,0}}^{Mp_{2n_1+2n_2,0}(\mb R)} \chi^{\alpha}$$
(see Notations ~\ref{m},  ~\ref{maxparabolic}).
Recall that representation $I^{\alpha}$ consists of sections of a certain homogeneous line bundle over the Lagrangian Grassmannian
$Sp_{2n_1+2n_2}(\mb R)/ P_{n_1+n_2,0}$ (see Notation ~\ref{maxparabolic}).
\begin{no}
Denote $Sp_{2n_1+2n_2}(\mb R)/ P_{n_1+n_2,0}$ by $X$. 
\end{no}
$X$ parametrizes
Lagrangian subspaces of a symplectic space of dimension $2n_1+2n_2$. We are thus led to the problem of analyzing the $Sp_{2n_1}(\mb R) \times Sp_{2n_2}(\mb R)$-action on $X$.\\
\\
The existence of an open dense $Sp_{2n_1}(\mb R) \times Sp_{2n_2}(\mb R)$-orbit 
in $X$ can be proved using algebraic group theory. But this is not sufficient for our purpose. Not only do we need to know the existence of an open dense orbit, but also the detailed structure of this orbit. For this reason, we give two elementary lemmas and one theorem in this section. We shall use them in the next section to identify the restriction of $I^{\alpha}$ to $Mp_{2n_1}(\mb R) \times Mp_{2n_2}(\mb R)$. The main result in this section is Theorem \!~\ref{iso}.\\
\\
Suppose from now on $n_2 > n_1 >0$. The case $n_1=n_2$ has been studied in \!~\cite{he00}, \!~\cite{he1}.
Let me begin by recalling the following result.
\begin{thm}[Theorem 0.3 \!~\cite{he00}]~\label{l}
Let $L(4n)$ be the Lagrangian Grassmannian of $(\mb R^{4n}, \Omega)$. Then
$Sp_{2n}(\mb R) \times Sp_{2n}(\mb R)$ acts on $L(4n)$ with $2n+1$ orbits.
Furthermore, there is an open dense orbit in $L(4n)$ which can be identified with $Sp_{2n}(\mb R)$.
\end{thm}

\subsection{Notations}
\begin{no}~\label{v}
Let $(\mb R^{2n_2+2n_1}, \Omega)$ be a symplectic space. Let
$X$ be the Lagrangian Grassmannian of $(\mb R^{2n_2+2n_1}, \Omega)$.
Fix a standard real basis $$\{e_1,e_2, \ldots, e_{n_2+n_1}, f_1, f_2, \ldots f_{n_2+n_1} \}$$ on $\mb R^{2n_1+2n_2}$ such that
$$\Omega(e_i, e_j)=0, \qquad \Omega(f_i, f_j)=0 $$
$$\Omega(f_j, e_i)= \delta_i^j, $$
where $\delta_i^j$ is the Kronecker symbol.
Write
$$\mathbb R^{2n_1}= span\{e_1, \ldots, e_{n_1}, f_1, \ldots, f_{n_1} \}$$
$$\mathbb R^{2n_2}= span\{e_{n_1+1}, \ldots, e_{n_1+n_2}, f_{n_1+1}, \ldots, f_{n_1+n_2} \}.$$ 
Let $Sp_{2n_1}(\mb R)$ be the subgroup of $Sp_{2n_1+2n_2}(\mb R)$ fixing
$e_j, f_j$ for every $j > n_1$. Let $Sp_{2n_2}(\mb R)$ be the subgroup of $Sp_{2n_1+2n_2}(\mb R)$ fixing
$e_j, f_j$ for every $j \leq n_1$. The group $Sp_{2n_1}(\mb R) \times Sp_{2n_2}(\mb R)$ is then
diagonally embedded into $Sp_{2n_1+2n_2}(\mb R)$. 
Let 
$$V_1=span\{e_1-e_{n_2+1}, e_2-e_{n_2+2}, \ldots e_{n_1}-e_{n_2+n_1}\}$$
$$V_2=span\{f_{n_2+1}+f_1, f_{n_2+2}+ f_2, \ldots, f_{n_1+n_2}+ f_{n_1} \}$$
$$V_0=span\{ e_{n_1+1}, \ldots, e_{n_2}\}$$
$$V_0^{\p}=span\{ f_{n_1+1}, \ldots f_{n_2} \}.$$
 Put
$V=V_1 \oplus V_2 \oplus V_0.$
\end{no}
\subsection{$Sp_{2n_2}(\mb R)$-action and the Generic Orbit $X_0$}
The space $V$ is a Lagrangian subspace of  $(\mb R^{2n_1+2n_2}, \Omega)$, or $V \in X$.
Consider the action of $Sp_{2n_2}(\mb R)$ on $V$.
Let $X_0$ be the $Sp_{2n_2}(\mb R)$-orbit of $V$. Then $X_0$ is a subset of $X$.
\begin{lem}~\label{simpleiso} Let $P_{n_2-n_1, n_1}$ be the maximal parabolic subgroup in $Sp_{2n_2}(\mb R)$ stabilizing $V_0 \subset \mathbb R^{2n_2}$. Let $Sp_{2n_1}(\mb R) 
GL(n_2-n_1) $ be the Levi subgroup stabilizing both $V_0$ and $V_0^{\p}$. Let
$N $ be the nilradical of $P_{n_2-n_1, n_1}$.
Then $(Sp_{2n_2}(\mb R))_{V} = GL(n_2-n_1) N$ and 
$$X_0 \cong Sp_{2n_2}(\mb R)/ GL(n_2-n_1) N.$$
Furthermore, $X_0$ is open in $X$.
\end{lem}
Proof: 
Let $g  \in Sp_{2n_2}(\mb R)$. Then $g$ fixes $e_i$, $f_i$ for every $i \leq n_1$. Suppose $g V =V$. 
\begin{enumerate}
\item Since $g$ stabilizes $V$ and $\mb R^{2n_2}$, $g$ stabilizes
$V_0=V \cap \mb R^{2n_2}$. So $g \in P_{n_2-n_1,n_1}$.
\item Let $j \leq n_1$. Since $g e_j=e_j$ and $g \mathbb R^{2n_2}= \mb R^{2n_2}$, 
$$ g (\mathbb R e_j \oplus \mb R^{2n_2})= \mb R e_j \oplus \mb R^{2n_2}.$$
It follows that
$$ g [(\mathbb R e_j \oplus \mb R^{2n_2}) \cap V] = [(\mb R e_j \oplus \mb R^{2n_2}) \cap V ].$$
Notice that
$$ (\mathbb R e_j \oplus \mb R^{2n_2}) \cap V= (\mathbb R e_j \oplus \mb R^{2n_2}) \cap (V_1 \oplus V_2 \oplus V_0)= \mb R (e_j- e_{n_2+j}) \oplus V_0.$$
So $$g (\mb R (e_j- e_{n_2+j}) \oplus V_0)= \mb R (e_j- e_{n_2+j}) \oplus V_0.$$
Since $g e_j=e_j$,
$$ g (e_j-e_{n_2+j}) \in (e_j-e_{n_2+j})+ V_0.$$
Hence for every $j \leq n_1$,
$g e_{n_2+j} \in  e_{n_2+j}+ V_0$.
\item Similarly, since $g f_j=f_j$ for $j \leq n_1$,
$$g f_{n_2+j} \in f_{n_2+j}+V_0, \qquad (j \leq n_1).$$
\item Thus, for every 
$$v \in {\rm span} \{e_{n_2+1},\ldots, e_{n_2+n_1}, f_{n_2+1}, \ldots f_{n_2+n_1} \} \simeq \mathbb R^{2n_1} ,$$
$g(v) \in  v + V_0.$  This shows that $g \in GL(n_2-n_1) N$.
\item Conversely, if $g \in GL(n_2-n_1) N$, then $g \in P_{n_2-n_1,n_1} \subset Sp_{2n_2}(\mb R)$. So $g$ fixes $e_j, f_j$ for every $j \leq n_1$ and $g$ stabilizes $V_0$. Since $g \in GL(n_2-n_1)N$, g preserves $e_{n_2+j}+V_0$ and $f_{n_2+j}+V_0$ for every $j \leq n_1$. So $g$ preserves
$e_j-e_{n_2+j}+V_0, f_{n_2+j}+f_j+V_0$ for every $j \leq n_1$.  It follows that $ g V=V$.
\end{enumerate}
Therefore, $(Sp_{2n_2}(\mb R))_{V} = GL(n_2-n_1) N$ and
$X_0 \cong Sp_{2n_2}(\mb R)/ GL(n_2-n_1) N.$
To show that $X_0$ is open in
$X$, we compute the dimensions of $X_0$ and $X$:
\begin{equation}
\begin{split}
\dim X_0= & \dim(Sp_{2n_1}(\mb R)) +\dim(Sp_{2n_2}(\mb R)/Sp_{2n_1}(\mb R) 
GL(n_2-n_1)N ) \\
=&  \dim(Sp_{2n_1}(\mb R)) +\frac{1}{2} \dim(Sp_{2n_2}(\mb R)/Sp_{2n_1}(\mb R) 
GL(n_2-n_1) ) \\
= &
2n_1^2+n_1+ \frac{2n_2^2+n_2-2n_1^2-n_1-(n_2-n_1)^2}{2} \\
= & \frac{4n_1^2+2n_1+2n_2^2+ n_2-2n_1^2- n_1-n_1^2+2n_1n_2-n_2^2}{2} \\
= & \frac{n_1^2+n_2^2+2n_1 n_2+(n_1+n_2)}{2} \\
= & \frac{2(n_1+n_2)^2+(n_1+n_2)-(n_1+n_2)^2}{2} \\
= & \frac{1}{2} \dim(Sp_{2n_1+2n_2}(\mb R)/ GL(n_1+n_2, \mb R))\\
= & \dim(Sp_{2n_1+2n_2}(\mb R)/ P_{n_1+n_2,0}) \\
= & \dim X.
\end{split}
\end{equation}
$\Box$. \\
\\
For $n_2=n_1$, $P_{n_2-n_1,n_1}=Sp_{2n_2}(\mb R)$. Our Lemma says that $X_0$ can be identified with
$Sp_{2n_2}(\mathbb R)$. This case is already treated in \!~\cite{he00} (see Theorem \!~\ref{l}).
\subsection{$Sp_{2n_1}(\mb R) \times Sp_{2n_2}(\mb R)$-Action}
Next, consider the $Sp_{2n_1}(\mb R) \times Sp_{2n_2}(\mb R)$ action on $X$.
Again, our symplectic space is $\mathbb R^{2n_1} \oplus \mathbb R^{2n_2}$.
The basis for $\mathbb R^{2n_1}$ is
$$\{e_1, e_2, \ldots e_{n_1}, f_1, f_2, \ldots, f_{n_1} \}.$$
Let 
$$ U=span\{e_{n_2+1}, e_{n_2+2}, \ldots, e_{n_2+n_1}, f_{n_2+1}, f_{n_2+2}, \ldots, f_{n_2+n_1} \}.$$
Identify $\mathbb R^{2n_1}$ with $U$ by mapping
$$ e_j \rightarrow e_{n_2+j}, f_j \rightarrow f_{n_2+j}, \qquad (j \in [1,n_1]) .$$
For each $g \in Sp_{2n_1}(\mb R)$, define $i(g) \in Sp_{2n_2}(\mb R)$ such that
$$i(g)|_{V_0 \oplus V_0^{\p}}=identity; \qquad i(g)|_{U}= g.$$
In other words, $Sp_{2n_1}(\mb R)$ is embedded into $Sp_{2n_2}(\mb R)$ as the first factor of the Levi subgroup
$Sp_{2n_1}(\mb R) GL(n_2-n_1)$. 
Recall for every $g \in Sp_{2n_1}(\mb R)$, $\tau(g)$ is defined to be
$$\arr{I_{n_1} & 0 \\ 0 & -I_{n_1} } g \arr{I_{n_1} & 0 \\ 0 & -I_{n_1}}.$$
\begin{lem} Let $g \in Sp_{2n_1}(\mb R)$. Let 
 $i(g)$ be an element in $ Sp_{2n_2}(\mb R)$ such that $i(g)|_{V_0 \oplus V_0^{\p}}=1$. Regarding $(g, \tau(i(g)))$ as an element in $Sp_{2n_1}(\mb R) \times Sp_{2n_2}(\mb R)$, we have
$$(g, \tau( i(g))) V= V.$$
\end{lem}
Proof: Let $g \in Sp_{2n_1}(\mb R)$. We borrow some ideas from \!~\cite{he00}. 
Define an isomorphism
$\phi:  \mathbb R^{2n_1} \rightarrow U $
by
$$\phi(e_j)= -e_{n_2+j} \qquad (j=1,2,\ldots, n_1)$$
$$\phi(f_j)= f_{n_2+j} \qquad (j=1,2,\ldots, n_1).$$
Then $V_1 \oplus V_2$ is the graph of $\phi$ and $\tau(i(g)) = \phi g \phi^{-1}$.
Notice that 
\begin{equation}
\begin{split}
 (g, \tau(i(g)))(V_1 \oplus V_2)= & \{(g u, \tau(i(g)) \phi(u) ) \mid u \in \mathbb R^{2n_1} \} \\
=& \{(u, \tau(i(g)) \phi (g^{-1} u)) \mid u \in \mb R^{2n_1}  \}\\
=& \{(u, \phi u) \mid u \in \mathbb R^{2n_1} \}\\
=& V_1 \oplus V_2.
\end{split}
\end{equation}
Thus $V_1 \oplus V_2$ is preserved by the actions of $(g, \tau(i(g)))$. Clearly, $(g, \tau(i(g)))$ also preserves vectors in $V_0$. It follows that
$(g, \tau(i(g))) V=V.$
$\Box$ \\
\\
Combining with Lemma \!~\ref{simpleiso}, we obtain
\begin{thm}~\label{iso} Let $P_{n_2-n_1,n_1}$ be the maximal parabolic subgroup of $Sp_{2n_2}(\mb R)$ which preserves $V_0$. Let $Sp_{2n_1}(\mb R) GL(n_2-n_1) N$
be the Langlands decomposition of $P_{n_2-n_1,n_1}$ as in Lemma \!~\ref{simpleiso}.
Then
the isotropy group
$$(Sp_{2n_1}(\mb R) \times Sp_{2n_2}(\mb R))_{V}=\{( \tau(g), gh) \mid
g \in Sp_{2n_1}(\mb R), h \in GL(n_2-n_1) N \}.$$
Moreover, $$X_0 \cong (Sp_{2n_1}(\mb R) \times Sp_{2n_2}(\mb R))/ (Sp_{2n_1}(\mb R) \times Sp_{2n_2}(\mb R))_{V}$$
and $X_0$ is open and dense in $X$.
\end{thm}
Here $Sp_{2n_1}(\mb R)$ and $Sp_{2n_2}(\mb R)$ are of the standard matrix form. The
$Sp_{2n_1}(\mb R)$ in the Levi factor $Sp_{2n_1}(\mb R) GL(n_2-n_1)$ of $Sp_{2n_2}(\mb R)$
is identified with the standard $Sp_{2n_1}(\mb R)$. We avoid using $i$ again.\\
\\
Proof: We only need to prove that $X_0$ is dense in $X$. This follows from an argument similar to the argument in \!~\cite{he00}. Briefly, for every Lagrangian subspace $W$, define
$$Ind(W)=(\dim(W \cap \mb R^{2n_1}), \dim (W \cap \mb R^{2n_2})).$$
Notice that 
$\dim (V \cap \mb R^{2n_1})=0$ and $\dim(V \cap \mb R^{2n_2})=\dim (V_0)=n_2-n_1$. In general,  for every $W \in X$, we have
$$\dim (W \cap \mb R^{2n_1}) \geq 0, \qquad \dim(W \cap \mb R^{2n_2}) \geq n_2-n_1.$$
We say a Lagrangian subspace $W$ is generic if $Ind(W)=(0, n_2-n_1)$.
Observe that the group action $Sp_{2n_1}(\mathbb R) \times Sp_{2n_2}(\mathbb R)$ preserves $Ind$. Therefore, every element in $X_0$ is generic. \\
\\
Conversely, 
every generic Lagrangian subspace $W$  must be an element in $X_0$. Let $W$ be generic. 
Since $W \cap \mb R^{2n_1} = \{0\}$, the projection
$$p_2: \mathbb R^{2n_1} \oplus \mathbb R^{2n_2} \rightarrow \mathbb R^{2n_2}$$
restricted to $W$ is injective.  So $\dim(p_2(W))=n_1+n_2$. Similarly, since
$\dim(W \cap \mb R^{2n_2})=n_2-n_1= \dim(W)-\dim(\mathbb R^{2n_1})$, the projection
$$p_1: \mathbb R^{2n_1} \oplus \mathbb R^{2n_2} \rightarrow \mathbb R^{2n_1}$$
restricted to $W$ must be surjective. It can then be shown that $W$ is uniquely determined by a symplectic map from $p_2(W)$ to $\mathbb R^{2n_1}$, namely $p_1 p_2^{-1}$. The kernel of this symplectic map is exactly
$W \cap \mathbb R^{2n_2}$, which is an isotropic subspace.  So $p_1 p_2^{-1}$ is reduced to a symplectic isomorphism from $p_2(W)/W \cap \mb R^{2n_2}$ to $\mathbb R^{2n_1}$. By basic linear group theory, the group $Sp_{2n_1}(\mb R) \times Sp_{2n_2}(\mb R)$ acts on the set of symplectic maps from $p_2(W)$ to $\mathbb R^{2n_1}$ transitively.  So the orbit $X_0$ consists of all generic Lagrangian subspaces. Therefore $X_0$ is dense in $X$. $\Box$

\section{ $I^{\alpha}(\pi)$ and $\Ind_{MP_{n_2-n_1, n_1}}^{Mp_{2n_2}(\mb R)} \pi^{\tau} \otimes \chi^{\alpha}$ }~\label{ialphainduced}
Following Notation \!~\ref{v}, Lemma \!~\ref{simpleiso}, and  Theorem \!~\ref{iso}, we set
\begin{no}~\label{embeding}
$$n_2 > n_1, \qquad P_{n_1+n_2, 0}= Sp_{2n_1+2n_2}(\mb R)_V.$$
Fix  $Sp_{2n_1}(\mb R) \times Sp_{2n_2}(\mb R) \subset Sp_{2n_1+2n_2}(\mb R)$ as in Notation \!~\ref{v}. Write
$$Sp_{2n_2}(\mb R)_{V}= GL(n_2-n_1) N.$$
Fix a maximal compact subgroup $MU(n_1+n_2)$ in $Mp_{2n_1+2n_2}(\mb R)$ such that
$$MU(n_1+n_2) \cap Mp_{2n_1}(\mb R), MU(n_1+n_2) \cap Mp_{2n_2}(\mb R)$$
are maximal compact subgroups for $Mp_{2n_1}(\mb R)$ and $Mp_{2n_2}(\mb R)$ respectively.
 Let $$MU(n_1)=MU(n_1+n_2) \cap Mp_{2n_1}(\mb R), \qquad MU(n_2)=MU(n_1+n_2) \cap Mp_{2n_2}(\mb R).$$
\end{no}
 Let $\pi \in \Pi(Mp_{2n_1}(\mb R))$. Recall from Definition \!~\ref{ipi} that
if  
$$( \phi_1 \otimes u_1,  \phi_2 \otimes u_2)_{1}=\int_{Mp_{2n_1}(\mb R)}(\pi(g_1) u_1, u_2)(I^{\alpha}(g_1) \phi_1, \phi_2) \, dg_1 $$
converges for every $u_1, u_2 \in V({\pi}), \phi_1, \phi_2 \in V({I^{\alpha}})$, then $I^{\alpha}(\pi)$  is defined to  be $V({I^{\alpha}}) \otimes V(\pi)$ modulo the radical of $(\ , \ )_{1}$. $I^{\alpha}(\pi)$ is only a $(\f{sp}_{2n_2}(\mb R), MU(n_2))$-module. In addition, $(\ , \ )_1$ descends into a nondegenerate $\f{sp}_{2n_2}(\mb R)$-invariant Hermitian form on $I^{\alpha}(\pi)$. The main purpose of this section is to show that
$(I^{\alpha}(\pi), (\ , \ )_1)$ is equivalent to the unitarily induced representation $$(V(\Ind_{MP_{n_2-n_1,n_1}}^{Mp_{2n_2}(\mb R)} \pi^{\tau} \otimes \chi^{\alpha}), (\ ,\ )).$$
Here $(\ , \ )$ denotes the inner product in the compact picture.
We do this in three steps.\\
\\
 The first step is to show that
$$I^{\alpha}|_{Mp_{2n_1}(\mb R) \times Mp_{2n_2}(\mb R)} \cong L^2(Mp_{2n_2}(\mb R) \times_{MGL(n_2-n_1) N} |\det|^{\frac{n_1+n_2}{2}} \otimes \chi^{\alpha}, X_0)$$
where $X_0$ is equipped with a $MU(n_2) \times Sp_{2n_1}(\mb R)$-invariant measure. On the right hand side of the equation, $Mp_{2n_2}(\mb R)$ acts on the left and $Mp_{2n_1}(\mb R)$ acts on the right with an involution $\tau$. See Theorem \!~\ref{mix} for the precise statement. The main ingredient of the proof is Theorem \!~\ref{iso}.\\
\\
From Lemma \!~\ref{simpleiso}, we have a principal fibration
$$Mp_{2n_1}(\mb R)/\{1, \epsilon \} \rightarrow X_0 \rightarrow Mp_{2n_2}(\mb R)/ Mp_{2n_1}(\mb R) MGL(n_2-n_1) N.$$
Let $C_c^{\infty}(X_0, I^{\alpha})$ be the set of compactly supported smooth sections of
$$Mp_{2n_2}(\mb R) \times_{MGL(n_2-n_1) N} |\det|^{\frac{n_1+n_2}{2}} \otimes \chi^{\alpha} \rightarrow X_0.$$
The next step is to define a map
$$\mc L^0: C_{c}^{\infty}(X_0, I^{\alpha}) \otimes V(\pi) \rightarrow \Ind_{MP_{n_2-n_1,n_1}}^{Mp_{2n_2}(\mb R)} \pi^{\tau} \otimes \chi^{\alpha},$$
This map is essentially integration on the fiber $Mp_{2n_1}(\mb R)$ (see Definition \!~\ref{l0}). The image of $\mc L^0$ are sections of
$$Mp_{2n_2}(\mb R) \times_{Mp_{2n}(\mb R) \times MGL(n_2-n_1) N} \pi^{\tau} \otimes |\det|^{\frac{n_1+n_2}{2}} \otimes \chi^{\alpha} \rightarrow X_0$$
in a proper sense.
\\
\\
Last step is to extend $\mc L^0$ to $V(I^{\alpha}) \otimes V(\pi)$ under the assumption that $(\ , \ )_1$ is well-defined. We will further show that $(\, ,\,  )_1$ on $I^{\alpha}(\pi)$ coincides with the Hilbert structure on
$$\Ind_{MP_{n_2-n_1,n_1}}^{Mp_{2n_2}(\mb R)} \pi^{\tau} \otimes \chi^{\alpha}.$$
The main result is Theorem \!~\ref{induction1}. Let us first analyze the restriction of $I^{\alpha}$ to $Mp_{2n_2}(\mb R)$.
\subsection{$I^{\alpha}$ restricted to $Mp_{2n_1} (\mb R) \times Mp_{2n_2}(\mb R)$}
The representation $I^{\alpha}$ is a parabolically induced representation.
Parabolic induction is a powerful tool in the study of representations of reductive Lie groups (see for example \!~\cite{knapp} and \!~\cite{wallach}). There are the compact picture
and the noncompact picture. To study $I^{\alpha}$, we recall
 the compact picture of $I^{\alpha}$.
\begin{no} Let $P_{n_2+n_1,0}$ be the Siegel parabolic subgroup of $Sp_{2n_1+2n_2}(\mb R)$ (See Notation \!~\ref{embeding}). Define a character on $P_{n_2+n_1,0}=GL(n_2+n_1) N$ by
$$\exp(-\rho_0(g n))=|\det g|^{-\frac{n_1+n_2+1}{2}} \qquad (g \in GL(n_1+n_2), n \in N).$$
Lift $\exp(-\rho_0)$ trivially to a character of $MP_{n_2+n_1,0}$. By abusing notation, denote this character by $\exp(-\rho_0)$.
\end{no}
The Hilbert space of $I^{\alpha}$ consists of $L^2$ sections of the linear bundle
$$ Mp_{2n_1+2n_2}(\mb R) \times_{MGL({n_1+n_2}) N} \mathbb C [\chi^{\alpha} \otimes \exp(\rho_0)] \rightarrow X.$$
Equivalently, 
\begin{defn}~\label{compactpic} 
The Hilbert space of $I^{\alpha}$ consists of  functions on $Mp_{2n_1+2n_2}(\mb R)$ such that
$$ f (g h)=  \chi(h)^{-\alpha} \exp(-\rho_0(h)) f(g) \qquad (\, \, \forall \, \, h \in MP_{n_2+n_1,0}),$$
$$ \int_{MU(n_1+n_2)} \| f(k) \|^2 \, d k < \infty. $$
The inner product between $f_1, f_2 \in I^{\alpha}$ is
$$(f_1, f_2)= \int_{MU(n_1+n_2)} f_1(k) \overline{f_2(k)} \, d k.$$
 Notice that the restriction map 
 $$ f \in I^{\alpha} \rightarrow f|_{MU(n_1+n_2)} \in \Ind_{MO(n_1+n_2)}^{MU(n_1+n_2)} \chi^{\alpha} $$
 defines an isomorphism between $I^{\alpha}$ and $\Ind_{MO(n_1+n_2)}^{MU(n_1+n_2)} \chi^{\alpha}$. 
We call the restrictions of $f$ to $MU(n_1+n_2)$ the compact picture. 
\end{defn}
Restrict the $L^2$-sections to the open dense subset $X_0$. Fix the base point $V$ (see Notation \!~\ref{v}). From Theorem \!~\ref{iso}, $X_0$ is a single $Sp_{2n_2}(\mb R) $-orbit of the form
$$Sp_{2n_2}(\mb R)/ GL(n_2-n_1) N.$$
We may model $I^{\alpha}|_{Mp_{2n_1}(\mb R) \times Mp_{2n_2}(\mb R)}$ on the sections of
$$Mp_{2n_2}(\mb R) \times_{M GL(n_2-n_1) N} \exp(\rho_0) \otimes \chi^{\alpha} \rightarrow X_0.$$
Since $X_0$ is open and dense in $X$, a section $f$ is uniquely determined by the section $ f|_{X_0}$ up to a set of measure zero. 
\begin{no} Fix the base point $V$ (see Notation \!~\ref{v}) on $X_0$. Regard $f \in I^{\alpha}$ as a section
on
$$Mp_{2n_2}(\mb R) \times_{M GL(n_2-n_1) N} \exp(\rho_0) \otimes \chi^{\alpha} \rightarrow X_0.$$
Let
$f^{\sharp}$ be the function on $Mp_{2n_2}(\mb R)$ corresponding to this section. 
\end{no}
 The function $f^{\sharp}$ satisfies
 \begin{equation}~\label{mp2n2model}
 f^{\sharp} (g h)=  \chi(h)^{-\alpha} \exp(-\rho_0(h))f^{\sharp} (g) \qquad (\forall \, \, h \in MGL(n_2-n_1) N) .
 \end{equation}
Let us consider the action of $Mp_{2n_1}(\mb R) \times Mp_{2n_2}(\mb R)$ on $I^{\alpha}$. Keep in mind that
$Mp_{2n_1}(\mb R) \cap Mp_{2n_2}(\mb R) = \{1, \epsilon\}$ and
$$I^{\alpha}(\epsilon, \epsilon) = Identity.$$
\begin{lem}~\label{g1g2} For $f \in I^{\alpha}$, $g_1 \in Mp_{2n_1}(\mb R)$, $g_2 \in Mp_{2n_2}(\mb R)$, 
$(g_1, g_2)$ embedded in $Mp_{2n_1+2n_2}(\mb R)$ diagonally, we have
$$f((g_1, g_2))=f^{\sharp}( g_2 \tau(g_1)^{-1}).$$
For every $h_1 \in Mp_{2n_1}(\mb R)$, we have
$$[I^{\alpha}(h_1)  f]^{\sharp}( g_2)=
f^{\sharp}(g_2 \tau(h_1)).$$
\end{lem}
Proof: By Theorem \!~\ref{iso},
\begin{eqnarray}
f((g_1, g_2))= & f((1, g_2 \tau(g_1)^{-1})(g_1, \tau(g_1))\\
=& f^{\sharp}(g_2 \tau(g_1)^{-1})
\chi((g_1,\tau(g_1))^{-\alpha} \exp(-\rho_0(g_1, \tau(g_1)))\\
= & f^{\sharp}(g_2 \tau(g_1)^{-1}).
\end{eqnarray}
As to the action of $h_1 \in Mp_{2n_1}(\mb R)$,
$$[I^{\alpha}(h_1) f]^{\sharp}( g_2)= (I^{\alpha}(h_1) f)((1,g_2)) =f((h_1^{-1}, g_2))=
f^{\sharp}(g_2 \tau(h_1)).$$
$\Box$.\\
\\
So for every $h_1 \in Mp_{2n_1}(\mb R)$, $[I^{\alpha}(h_1) f]^{\sharp}=R(\tau(h_1)) f^{\sharp}$. Here $R$ stands for the right regular action. Clearly, the action of $Mp_{2n_2}(\mb R)$ on $f^{\sharp}$ is the left regular action.
 Now we shall determine the Hilbert structure  on 
$f^{\sharp}$. More precisely, for $f_1, f_2 \in I^{\alpha}$, we want to express
$(f_1, f_2)$ in terms of $f_1^{\sharp}, f_2^{\sharp}$.

\begin{no} 
Let $d k_2$ be the normalized invariant measure for the maximal compact subgroup $MU(n_2)$. Let $d [k_2]$ be the induced measure on 
$$ MU(n_2)/M(O(n_2-n_1)U(n_1)) \cong Mp_{2n_2}(\mb R)/ MP_{n_2-n_1, n_1} .$$
Fix an invariant measure 
$ d g_1$ on $Mp_{2n_1}(\mb R)$. Every element $x \in  Mp_{2n_2}(\mb R)/ MGL(n_2-n_1) N$ will be represented by  
$k_2 g_1 $ for some $k_2 \in MU(n_2)$ and some $g_1 \in Mp_{2n_1}(\mb R)$.
This is the case since $ Sp_{2n_2}(\mb R)$ has a $U(n_2) Sp_{2n_1}(\mb R) GL(n_2-n_1) N$ decomposition. Let $ d g_1 d [k_2]$ be the induced measure on $Mp_{2n_2}(\mb R)/ MGL(n_2-n_1) N$.
 \end{no}
\begin{defn}~\label{sharpproduct} Let 
$$L^2(Mp_{2n_2}(\mb R) \times_{MGL(n_2-n_1) N} |\det|^{\frac{n_1+n_2}{2}} \otimes \chi^{\alpha},  d g_1 d [k_2])$$
be the set of sections of 
$$Mp_{2n_2}(\mb R) \times_{MGL(n_2-n_1) N} |\det|^{\frac{n_1+n_2}{2}} \otimes \chi^{\alpha}
\rightarrow Mp_{2n_2}(\mb R)/ MGL(n_2-n_1) N$$
such that
$$\|f\|^2=\int_{Mp_{2n_2}(\mb R)/ MGL(n_2-n_1) N} |f(k_2 g_1)|^2  \, d [k_2] \, d g_1 < \infty.$$
Define
$$(f_1, f_2)^{\sharp}= \int_{Mp_{2n_2}(\mb R)/ MGL(n_2-n_1) N} f_1(k_2 g_1) \overline{f_2(k_2 g_1)} \, d [k_2] \, dg_1.$$
This integral does not depend on the choices of $k_2$. Let $g_2 \in Mp_{2n_2}(\mb R)$ act on 
$f$
by the left regular action $L(g_2)$ and  $g_1 \in Mp_{2n_1}$ act on $f$ by the right regular action $R(\tau(g_1))$.
\end{defn}

\begin{thm}~\label{mix} Up to a proper choice of $d g_1$, the map $f \rightarrow f^{\sharp}$ defines an isometry between
$I^{\alpha}$ and 
$$L^2(Mp_{2n_2}(\mb R) \times_{MGL(n_2-n_1) N} |\det|^{\frac{n_1+n_2}{2}} \otimes \chi^{\alpha},  d g_1 d [k_2]).$$
In addition, this isometry induces an equivalence of $Mp_{2n_1}(\mb R) \times Mp_{2n_2}(\mb R)$-representations:
$$I^{\alpha}|_{Mp_{2n_1}(\mb R) Mp_{2n_2}(\mb R)} \cong L^2(Mp_{2n_2}(\mb R) \times_{MGL(n_2-n_1) N} |\det|^{\frac{n_1+n_2}{2}} \otimes \chi^{\alpha},  d g_1 d [k_2]).$$
\end{thm}
We call this the $Mp_{2n_2}(\mb R)$-model of $I^{\alpha}$.\\
\\
Proof: Let $f_1, f_2 \in I^{\alpha}$. Clearly, $(f_1, f_2)$ can be expressed as an integral on
$X_0$ by ignoring the measure zero set $X-X_0$. Then
$(f_1, f_2)$ can be expressed as
$$\int_{Mp_{2n_2}(\mb R)/ MGL(n_2-n_1) N} f_1^{\sharp}(k_2 g_1) \overline{f_2^{\sharp}(k_2 g_1)} \Delta([k_2], g_1) \, d [k_2] \, dg_1.$$
Here $\Delta([k_2], g_1)$ comes from a product of the Jacobian of a coordinate transform 
$$ Mp_{2n_2}(\mb R)/ MGL(n_2-n_1) N \ni ([k_2], g_1) \longrightarrow [k] \in MU(n_1+n_2)/MO(n_1+n_2) \cong X $$
and an $\exp(-\rho_0)$ term. Since $(f_1,f_2)=(I^{\alpha}(k_2 g_1) f_1, I^{\alpha}(k_2 g_1) f_2)$ for 
$k_2 \in MU(n_2)$ and $g_1 \in Mp_{2n_1}(\mb R)$, $\Delta([k_2], g_1) )$ must be left $MU(n_2)$ invariant and right $Mp_{2n_1}(\mb R)$-invariant. Hence $\Delta([k_2], g_1) )$ is a constant. The first assertion is proved. \\
\\
The second assertion follows from Lemma \!~\ref{g1g2}. $\Box$
\subsection{$\mc L^0: C_c^{\infty}(X_0, I^{\alpha}) \otimes V(\pi) \rightarrow \Ind_{MP_{n_2-n_1,n_1}}^{Mp_{2n_2}(\mb R)} \pi^{\tau} \otimes \chi^{\alpha}$}
 \begin{defn}
Let $f$ be an element in the Hilbert space of $I^{\alpha}$. Let $f^{\sharp}$ be the restriction of $f$ to $Mp_{2n_2}(\mb R)$.
We say that a vector $f$ is compactly supported on $X_0$, if $f^{\sharp}$ is compactly supported on $X_0 \cong Mp_{2n_2}(\mb R)/ MGL(n_2-n_1)N$. Let $C_c^{\infty}(X_0, I^{\alpha})$ be the set of smooth vectors with compact support in $X_0$. Let 
$C_{c}^{\infty}(X_0, I^{\alpha})_{MU(n_2)}$ be the $MU(n_2)$-finite vectors in $C_c^{\infty}(X_0, I^{\alpha})$.
\end{defn}
\begin{defn}~\label{l0} Let 
$\pi \in \Pi_u(Mp_{2n_1}(\mb R))$ and $u \in V(\pi)$.
For each $f \in C_c^{\infty}(X_0, I^{\alpha})$, define
$\mc L^0(f \otimes u)$ to be the function
$$g_2 \in Mp_{2n_2}(\mb R) \rightarrow \int_{Mp_{2n_1}(\mb R)} ( I^{\alpha}(g_1) f)^{\sharp}(g_2) \ \pi(g_1) u \ d g_1 \in \mc H_{\pi}.$$
\end{defn}
Since $f^{\sharp}$ is compactly supported, for each $g_2 \in Mp_{2n_2}(\mb R)$,
$$\int_{Mp_{2n_1}(\mb R)} ( I^{\alpha}(g_1) f)^{\sharp}(g_2) \ \pi(g_1) u \ d g_1=
\int_{Mp_{2n_1}(\mb R)}  f^{\sharp}(g_2 \tau(g_1) ) \ \pi(g_1) u \ d g_1.$$
So $\mc L^0(f \otimes u)$ is well-defined as a vector valued function on  $Mp_{2n_2}(\mb R)$.
\begin{thm}~\label{inte} Let $f \in C_c^{\infty}(X_0, I^{\alpha})$. Let 
$\pi \in \Pi_u(Mp_{2n_1}(\mb R))$ and $u \in V(\pi)$. Then
$\mc L^0(f \otimes u)$
is an element in the Hilbert space of 
$\Ind_{MP_{n_2-n_1,n_1}}^{Mp_{2n_2}(\mb R)} \pi^{\tau} \otimes  \chi^{\alpha}.$
Furthermore,  
$$\mc L^0(C_{c}^{\infty}(X_0, I^{\alpha})_{MU(2)} \otimes V(\pi)) =
 V(\Ind_{MP_{n_2-n_1,n_1}}^{Mp_{2n_2}(\mb R)} \pi^{\tau} \otimes \chi^{\alpha})$$
and $\mc L^0$ preserves the $({\f{sp}}_{2n_2}(\mb R), MU(n_2))$ action.
\end{thm}
Proof:  Let $f \in C_c^{\infty}(X_0, I^{\alpha})$ and $u \in V(\pi)$. 
\begin{enumerate}
\item Notice that for $g_2 \in Mp_{2n_2}(\mb R)$ and $h_1$ in the $Mp_{2n_1}(\mb R)$ factor of $MP_{n_2-n_1,n_1}$, 
\begin{equation}
\begin{split}
\mc L^0(f \otimes u)(g_2 h_1)= & \int_{Mp_{2n_1}(\mb R)} ( I^{\alpha}(g_1) f)^{\sharp}(g_2 h_1) \ \pi(g_1) u \  d g_1 \\
= & \int_{Mp_{2n_1}(\mb R)} f^{\sharp}( g_2 h_1 \tau(g_1)) \  \pi(g_1) u \ d g_1 \\
= & \int_{Mp_{2n_1}(\mb R)} f^{\sharp}( g_2 h_1 g_1) \ \pi(\tau(g_1))u \ d g_1 \\
=& \int_{Mp_{2n_1}(\mb R)} f^{\sharp}(g_2 g_1) \ \pi(\tau(h_1^{-1} g_1)) u \ d g_1 \\
=& \int_{Mp_{2n_1}(\mb R)} f^{\sharp}(g_2 g_1) \ \pi(\tau(h_1)^{-1}) \pi(\tau( g_1)) u \ d g_1\\
=& \pi(\tau(h_1)^{-1}) \{ \int_{Mp_{2n_1}(\mb R)} ( I^{\alpha}(\tau(g_1)) f)^{\sharp}(g_2) 
\ \pi(\tau( g_1)) u \ d g_1 \} \\
=& \pi(\tau(h_1)^{-1})[ \int_{Mp_{2n_1}(\mb R)} ( I^{\alpha}(g_1) f)^{\sharp}(g_2) 
\ \pi( g_1) u \ d g_1 ].   
\end{split}
\end{equation}
So $\mc L^0(f \otimes u)(g_2 h_1)=\pi(\tau(h_1)^{-1}) \mc L^0(f \otimes u)(g_2)$.
\item  Recall that for every $h \in MGL(n_2-n_1)N \subset MP_{n_2-n_1, n_1}$, we have
$$ f^{\sharp} (g_2 h)= f^{\sharp} (g_2) \chi(h)^{-\alpha} \exp(-\rho_0(h)). $$
It follows that
$$ \int_{Mp_{2n_1}(\mb R)} ( I^{\alpha}(g_1) f)^{\sharp}(g_2 h) \ \pi(g_1) u \ d g_1= \chi(h)^{-\alpha} \exp(-\rho_0(h)) \int_{Mp_{2n_1}(\mb R)} ( I^{\alpha}(g_1) f)^{\sharp}(g_2 ) \ \pi(g_1) u \ d g_1.$$
So $\mc L^0(f \otimes u)(g_2 h) =\chi(h)^{-\alpha} \exp(-\rho_0(h)) \mc L^0(f \otimes u)(g_2) $.
\item 
Since $f^{\sharp}$ is smooth and compactly supported, the restriction of
$$g_2 \rightarrow \int_{Mp_{2n_1}(\mb R)} ( I^{\alpha}(g_1) f)^{\sharp}(g_2) \ \pi(g_1) u \ d g_1$$
onto the maximal compact subgroup $MU(n_2)$ must be continuous and bounded. We see now that $\mc L^0(f \otimes u)$ is a continuous section of
$$Mp_{2n_2}(\mb R) \times_{Mp_{2n_1}(\mb R) MGL(n_2-n_1) N } \pi^{\tau} \otimes \exp(\rho_0) \otimes \chi^{\alpha} \rightarrow Mp_{2n_2}(\mb R)/Mp_{2n_1}(\mb R) GL(n_2-n_1) N.$$
Notice that $\exp(\rho_0)= |\det|^{\frac{n_1+n_2}{2}}$ on the $GL$ factor is exactly the normalizing factor.
Therefore
$\mc L^0(f \otimes u)$
is an element in the Hilbert space of 
$\Ind_{MP_{n_2-n_1,n_1}}^{Mp_{2n_2}(\mb R)} \pi^{\tau} \otimes \chi^{\alpha}.$
\item Since the right action of $Mp_{2n_1}(\mb R)$ commutes with the left action of $Mp_{2n_2}(\mb R)$, $\mc L^0$ preserves the action of ${ Mp}_{2n_2}(\mb R)$. In particular, $\mc L^0$ preserves the action of $MU(n_2)$. Therefore
$$\mc L^0: C_{c}^{\infty}(X_0, I^{\alpha})_{MU(n_2)} \otimes V(\pi) \rightarrow V(\Ind_{MP_{n_2-n_1,n_1}}^{Mp_{2n_2}(\mb R)} \pi^{\tau} \otimes \chi^{\alpha}).$$
\item Suppose that the map above is not surjective. There must exist a nonzero $MU(n_2)$-finite vector
$$\phi: Mp_{2n_2}(\mb R) \rightarrow \mc H_{\pi^{\tau}}$$
such that 
$$\phi(g_2 p_2)=(\pi^{\tau} \otimes \chi^{\alpha})(p_2^{-1}) \phi(g_2) \qquad (\forall \ \ p_2 \in MP_{n_2-n_1,n_1}),$$
$$(\phi, \mc L^{0}(f \otimes u))=0 \qquad (\forall \ \ f \in C_{c}^{\infty}(X_0, I^{\alpha})_{MU(n_2)}, u \in V(\pi)).$$
By an easy computation,
$$0=(\phi, \mc L^0(f \otimes u))=\int_{Mp_{2n_2}(\mb R)/MGL(n_2-n_1)N} (\phi( k_2 g_1), u) \overline{f(k_2 g_1)} \, d [k_2] \, dg_1.$$
We see that the $MU(n_2)$-finite function $k_2 g_1 \rightarrow (\phi(k_2 g_1), u)$ is perpendicular to all compactly supported $MU(n_2)$-finite functions on $Mp_{2n_2}(\mb R)/MGL(n_2-n_1)N$ with respect to the measure $ d[k_2] d g_1$. So $(\phi(k_2 g_1), u)=0$ for all $k_2 \in MU(n_2), g_1 \in Mp_{2n_1}(\mb R), u \in V(\pi^{\tau})$. Therefore
$\phi=0$. We reach a contradiction. Hence
$$\mc L^0(C_{c}^{\infty}(X_0, I^{\alpha})_{MU(2)} \otimes V(\pi)) =
 V(\Ind_{MP_{n_2-n_1,n_1}}^{Mp_{2n_2}(\mb R)} \pi^{\tau} \otimes \chi^{\alpha}).$$
\end{enumerate}
$\Box$

\subsection{The Inner Product $(\ , \ )_1$}
We have shown that $\mc L^0$ defines a map from $C_{c}^{\infty}(X_0, I^{\alpha}) \otimes V(\pi)$ to $\Ind_{MP_{n_2-n_1,n_1}}^{Mp_{2n_2}(\mb R)} \pi^{\tau} \otimes \chi^{\alpha}$.
In this section, we shall extend the definition of $\mc L^0$ to $V(I^{\alpha}) \otimes V(\pi)$ under the assumption that $(\ , \ )_1$ is well-defined for $V(I^0) \otimes V(\pi)$. Let us begin with the following lemma.
\begin{lem}~\label{1k} Let $MU(n_1+n_2)$ be a maximal compact subgroup of $Mp_{2n_1+2n_2}(\mb R)$. Let $1_{MU(n_1+n_2)}$ be the spherical vector in $I^{0}$. Then 
\begin{enumerate}
\item $g \rightarrow (I^0(g) 1_{MU(n_1+n_2)}, 1_{MU(n_1+n_2)})$ is positive;
\item Each $\phi \in V(I^{\alpha})$, as a function on $Mp_{2n_1+2n_2}(\mb R)$, is bounded by a multiple of $1_{MU(n_1+n_2)}$;
\item For every $\phi, \psi \in V(I^{\alpha})$, there exists a constant $C_{\phi,\psi}$ such that $$|(I^{\alpha}(g) \phi, \psi)| \leq C_{\phi, \psi}(I^0(g) 1_{MU(n_1+n_2)}, 1_{MU(n_1+n_2)})$$ for every $g \in Mp_{2n_1+2n_2}(\mb R)$;
\item If $(1_{MU(n_1+n_2)} \otimes u , 1_{MU(n_1+n_2)} \otimes v)_1$ converges absolutely for all $u, v \in V(\pi)$, then $(\ , \ )_1$ is well-defined for $V(I^{\alpha}) \otimes V(\pi)$.
\end{enumerate}
\end{lem}
Proof: In the compact picture, $1_{MU(n_1+n_2)}$ is the constant function $1$ on $MU(n_1+n_2)$. $I^{0}(g) 1_{MU(n_1+n_2)}$ is also a positive function. Therefore $(I^0(g) 1_{MU(n_1+n_2)}, 1_{MU(n_1+n_2)}) > 0$. By continuity,  every $MU(n_1+n_2)$-finite function in the compact picture of $I^{\alpha}$ is bounded by a positive multiple of $1_{MU(n_1+n_2)}$. Hence every $MU(n_1+n_2)$-finite matrix coefficient of $I^{\alpha}$
is bounded by a multiple of $(I^0(g) 1_{MU(n_1+n_2)}, 1_{MU(n_1+n_2)})$, i.e.,
$$| (I^{\alpha}(g) u, v) | \leq  C_u C_v (I^0(g) 1_{MU(n_1+n_2)}, 1_{MU(n_1+n_2)}) \qquad (u, v \in V(I^{\alpha}))$$
for some constants $C_u$ and $C_v$. As a result, if $(u \otimes 1_{MU(n_1+n_2)} , v \otimes 1_{MU(n_1+n_2)})_1$ converges absolutely for all $u, v \in V(\pi)$ and $1_{MU(n_1+n_2)}$ in $I^0$, then $(\ , \ )_1$ is well-defined for $V(I^{\alpha}) \otimes V(\pi)$.
$\Box$
\begin{lem}~\label{uniformsharp}
Let $\phi \in V(I^{\alpha})$. Then there exists a positive constant $C_{\phi}$ such that $|\phi^{\sharp}(k_2 g_1)| \leq C_{\phi} 1^{\sharp}_{MU(n_1+n_1)}(g_1)$ for all $k_2 \in MU(n_2)$ and $g_1 \in Mp_{2n_1}(\mb R)$.
\end{lem}
Proof: First of all $1_{MU(n_1+n_2)}^{\sharp}(k_2 g_1)= \exp (-\rho_0(p(g_1))$ where
$g_1=k(g_1) p(g_1)$ under the decomposition $MU(n_1+n_2) MP_{n_1+n_2, 0}$ of $Mp_{2n_1+2n_2}(\mb R)$. Here $g_1 \in Mp_{2n_1}(\mb R) \subseteq Mp_{2n_2}(\mb R) \subseteq Mp_{2n_1+2n_2}(\mb R )$ (see Notation \!~\ref{embeding}). Similarly,
$\phi^{\sharp}(k_2 g_1)=\phi(k_2 k(g_1)) \exp(-\rho_0(p(g_1))$. By Lemma \!~\ref{1k} (2),
$\phi(k_2 k(g_1))$ is bounded by a constant $C_{\phi}$ for all $k_2 \in MU(n_2), k(g_1) \in MU(n_1+n_2)$. So 
$$|\phi^{\sharp}(k_2 g_1)| \leq C_{\phi} 1^{\sharp}_{MU(n_1+n_1)}(g_1)$$
for all $k_2 \in MU(n_2)$ and $g_1 \in Mp_{2n_1}(\mb R)$. $\Box$.
\begin{lem}~\label{1k2}
Let $\pi \in \Pi_u(Mp_{2n_1}(\mb R)$.
The following are equivalent.
\begin{enumerate}
\item $I^0(\pi)$ is well-defined.
\item $(1_{MU(n_1+n_2)} \otimes u , 1_{MU(n_1+n_2)} \otimes v)_1$ converges absolutely for all $u, v \in V(\pi)$.
\item For every $u, v \in V(\pi)$, the following integral
$$\int_{Mp_{2n_1}(\mb R)} (R(\tau(h_1)) 1_{MU(n_1+n_2)}^{\sharp}(g_1), 1_{MU(n_1+n_2)}^{\sharp}(g_1))_{Mp_{2n_1}(\mb R)} (\pi(h_1) u, v) \, d h_1$$
converges absolutely. Here $(\ , \ )_{Mp_{2n_1}(\mb R)}$ is the inner product for  $L^2(Mp_{2n_1}(\mb R))$.
\item For every $u, v \in V(\pi)$, the function
$$(g_1, h_1) \mapsto 1_{MU(n_1+n_2)}^{\sharp}(g_1 \tau(h_1)) 1_{MU(n_1+n_2)}^{\sharp}(g_1) |(\pi(h_1) u, v)|$$
is in $L^1(Mp_{2n_1}(\mb R) \times Mp_{2n_1}(\mb R), d g_1 d h_1)$.
\item $I^{\alpha}(\pi)$ is well-defined for every $\alpha$.
\end{enumerate}
\end{lem}
Proof: By Definition \!~\ref{sharpproduct} and Theorem \!~\ref{mix}, we have $$(L(h_1) 1_{MU(n_1+n_2)}, 1_{MU(n_1+n_2)})= (R(\tau(h_1)) 1^{\sharp}_{MU(n_1+n_2)}, 1^{\sharp}_{MU(n_1+n_2))})^{\sharp}.$$
The right hand side is equal to
$$ C (R(\tau(h_1)) 1_{MU(n_1+n_2)}^{\sharp}(g_1), 1_{MU(n_1+n_2)}^{\sharp}(g_1))_{Mp_{2n_1}(\mb R)}$$
where $C$ is a constant related to the measure of $Mp_{2n_2}(\mb R)/MP_{n_2-n_1, n_1}$.
We obtain that 
\begin{equation}
\begin{split}
 & (1_{MU(n_1+n_2)} \otimes u , 1_{MU(n_1+n_2)} \otimes v)_1\\
= & C \int_{Mp_{2n_1}(\mb R)} (R(\tau(h_1)) 1_{MU(n_1+n_2)}^{\sharp}(g_1), 1_{MU(n_1+n_2)}^{\sharp}(g_1))_{Mp_{2n_1}(\mb R)} (\pi(h_1) u, v) \, d h_1.
\end{split}
\end{equation}
So $(2) \leftrightarrow (3)$. $(3) \leftrightarrow (4)$ follows from Fubini's theorem for nonnegative functions. \\
\\
Suppose $I^0(\pi)$ is well-defined. Then $(1_{MU(n_1+n_2)} \otimes u , 1_{MU(n_1+n_2)} \otimes v)_1$ converges absolutely for all $u, v \in V(\pi)$. So $(1) \rightarrow (2)$. By Lemma \!~\ref{1k} (4), $I^{\alpha}(\pi)$ is well-defined for every $\alpha$. So $(2) \rightarrow (5)$. $(5) \rightarrow (1)$ is obvious. We have proved $ (1) \rightarrow (2) \rightarrow (5) \rightarrow (1)$. $\Box$\\
\\
Suppose that $I^0(\pi)$ is well-defined. Let $\phi, \psi \in V(I^{\alpha})$. By Lemma \!~\ref{uniformsharp} and Lemma \!~\ref{1k2} (4),  
$$ \int \int_{ Mp_{2n_1}(\mb R) \times Mp_{2n_1}(\mb R)} | \phi^{\sharp}(k_2 g_1) \overline{\phi^{\sharp}(k_2 h_1) }  (\pi(\tau(g_1)) u ,  \pi(\tau(h_1)) u)|  \, d h_1 \, dg_1  < \infty $$
for all $k_2 \in MU(n_2)$. 
\begin{lem}~\label{fa} Let $\mc H$ be a separable Hilbert space and $\Phi$ be a $\mc H$-valued continuous function on $Mp_{2n_1}(\mb R)$. If $$\int \int_{x,y \in Mp_{2n_1}(\mb R)}  | (\Phi(x), \Phi(y))| \,d x \, d y < \infty,$$
 then
$\int \Phi(x) \, d x \in \mc H$ is well-defined, as the unique limit of $\int_K \Phi(x) \, d x$ as $K \rightarrow Mp_{2n_1}(\mb R)$. Furthermore,
$$(\int_{x \in Mp_{2n_1}(\mb R)} \Phi(x) d x, \int_{y \in Mp_{2n_1}(\mb R)} \Phi(y) d y)=\int_{x,y \in Mp_{2n_1}(\mb R)}  (\Phi(x), \Phi(y)) d x dy$$
\end{lem}
The proof of this will be given in the Appendix. See Theorem \!~\ref{ivh1}.
\begin{defn} Suppose that $I^0(\pi)$ is well-defined for a $\pi \in \Pi_u(Mp_{2n_1}(\mb R))$.  Let $\phi \in V(I^{\alpha})$ and $u \in V(\pi)$. Consider the $\mc H_{\pi}$-valued function
$$ g_1 \in Mp_{2n_1}(\mb R) \rightarrow \phi^{\sharp}(k_2 g_1) \pi(\tau(g_1)) u \qquad (k_2 \in MU(n_2))$$
Define
$$\mc L^0(\phi \otimes u):  MU(n_2) \ni k_2 \mapsto \int_{g_1} \phi^{\sharp}(k_2 g_1) \pi(\tau(g_1)) u \, d g_1 \in \mc H_{\pi}.$$
By  Lemma \!~\ref{fa}, $\mc L^0(\phi \otimes u)$ is a $\mc H_{\pi}$-valued function on $MU(n_2)$. Now
$\mc L^0$ is defined on $C_{c}^{\infty}(X_0, I^{\alpha})  \otimes V(\pi) + V(I^{\alpha}) \otimes V(\pi)$.
\end{defn}
\begin{thm}~\label{innerproduct} Let $\pi \in \Pi_u(Mp_{2n_1}(\mb R))$. 
Suppose that $(\ , \ )_1$ is well-defined for $V(I^{0}) \otimes V(\pi)$.
Then for every $u,v \in V(\pi)$ and $\phi, \psi \in V(I^{\alpha})$, we have
\begin{enumerate}
\item  $\mc L^0( \phi \otimes u), \mc L^0( \psi \otimes v) \in 
\Ind_{MP_{n_2-n_1,n_1}}^{Mp_{2n_2}(\mb R)} \pi^{\tau} \otimes \chi^{\alpha};$ 
\item
$ 2 ( \phi \otimes u,  \psi \otimes v)_{1}=(\mc L^0( \phi \otimes u) , \mc L^0( \psi \otimes v) )_{\Ind_{\,} \pi^{\tau} \otimes \chi^{\alpha}}.$
\end{enumerate}
In particular, the Hermitian form $(\ , \ )_1$ is positive semidefinite and $\mc L^0(V(I^{\alpha}) \otimes V(\pi)) \cong I^{\alpha}(\pi).$
\end{thm}
Proof: To prove $(1)$, we will construct a sequence $\phi_m \in C_{c}^{\infty}(X_0, I^{\alpha}) \rightarrow \phi$ such that $\mc L^0(\phi_m \otimes u)(k)$ is a sequence approaching $\mc L^0(\phi \otimes u)(k)$ in $\mc H_{\pi}$ uniformly for all $k \in MU(n_2)$. We start with a sequence of functions $\mu_m \in C_c^{\infty}(Mp_{2n_1}(\mb R))$ such that
\begin{enumerate}
\item $\mu_m(k_1 g k_1^{\prime})=\mu_m(g), \forall \ k_1, k_1^{\prime} \in MU(n_1)$;
\item $\mu_m(\exp H)=1$ for all $\|H\| \leq m$ where $H \in \f a$ and $\| H \|^2=|(H, H)|$ (see 1.1).
\item $0 \leq \mu_m(g) \leq 1$ for all $g$.
\end{enumerate}
Define $\phi_m \in C_c^{\infty}(X_0, I^{\alpha})$ such that $\phi_m^{\sharp}(k_2 g_1)= \phi^{\sharp}(k_2 g_1) \mu_m(g_1)$. By Lemma \!~\ref{uniformsharp} and the proof of Lemma \!~\ref{fa} for the function $g_1 \rightarrow 1_{MU(n_1+n_2)}^{\sharp}(g_1) \pi(\tau(g_1)) u$, we have
\begin{equation}
\begin{split}
 & \|\mc L^0(\phi_{m} \otimes u)(k_2)- \mc L^0(\phi \otimes u)(k_2)\|^2\\
 =& \int \int_{Mp_{2n_1}(\mb R) \times Mp_{2n_1}(\mb R)} [\phi_m^{\sharp}(k_2 g_1)-\phi^{\sharp}(k_2 g_1)] \overline{[\phi_m^{\sharp}(k_2 h_1)-\phi^{\sharp}(k_2 h_1)]} (\pi(\tau(g_1)) u, \pi(\tau(h_1)) u) \, dg_1 \, d h_1\\
\leq & 4 C_{\phi}^2 \int \int_{\|g_1 \| \geq m \ or \  \|h_1\| \geq m} |(1^{\sharp}_{MU(n_1+n_2)}(g_1) \pi(\tau(g_1)) u, 1^{\sharp}_{MU(n_1+n_2)}(h_1) \pi(\tau(h_1)) u)|  \, dg_1 \, d h_1\\
 \rightarrow & 0
\end{split}
\end{equation}
where $\|g_1 \|= \|H \|$ for $g_1=k_1 \exp H k_1^{\prime}$. So
$$\mc L^0(\phi_m \otimes u)(k_2) \rightarrow \mc L^0(\phi \otimes u)(k_2)$$
uniformly for all $k_2$.
By the proof of Theorem \!~\ref{inte}, each $\mc L^0(\phi_{m} \otimes u)$ is continuous in $$\Ind_{MP_{n_2-n_1,n_1}}^{Mp_{2n_2}(\mb R)} \pi^{\tau} \otimes \chi^{\alpha}.$$
Therefore $\mc L^0(\phi \otimes u)$ is also continuous in $\Ind_{MP_{n_2-n_1,n_1}}^{Mp_{2n_2}(\mb R)} \pi^{\tau} \otimes \chi^{\alpha}$.\\
\\
To prove (2), we compute
\begin{equation}~\label{long}
\begin{split}
 & 2 ( \phi \otimes u,  \psi \otimes v)_{1} \\
= &  2 \int_{g_1 \in Mp_{2n_1}(\mb R)} (\pi(g_1) u, v) ((I^{\alpha}(g_1) \phi)(x) , \psi(x)) \, dg_1 \\
= &  \int_{g_1 \in Mp_{2n_1}(\mb R)} (\pi(g_1) u, v) \int_{Mp_{2n_2}(\mb R)/MGL(n_2-n_1)N} \phi^{\sharp}( k_2 h_1 \tau(g_1)) \overline{\psi^{\sharp}}( k_2 h_1) \, d [k_2] \, d h_1 \ d g_1 \\
= &   \int_{g_1 \in Mp_{2n_1}(\mb R)} (\pi(\tau(g_1))u, v)  \int_{[k_2] \in Mp_{2n_2}(\mb R)/MP_{n_2-n_1,n_1}} \int_{ h_1 \in Mp_{2n_1}(\mb R)}   \phi^{\sharp}(k_2 h_1 g_1) \overline{\psi^{\sharp}}(k_2 h_1)  \, d h_1 \, d [k_2] \, dg_1
\end{split}
\end{equation}
By Lemma \!~\ref{1k} (1), the triple integral (\!~\ref{long}) 
converges absolutely for $\phi, \psi=1_{MU(n_1+n_2)} \in I^0$. By Lemma \!~\ref{1k} (2)(3), the triple integral (\!~\ref{long}) converges absolutely for all  $\phi, \psi \in V(I^{\alpha})$. So by Fubini's Theorem, changing the order of integrations in Integral (\!~\ref{long}), we obtain
\begin{equation}
\begin{split}
 & 2 (\phi \otimes u, \psi \otimes v)_1\\
= & \int_{[k_2] \in Mp_{2n_2}(\mb R)/ MP_{n_2-n_1}}  \int \int_{ Mp_{2n_1}(\mb R) \times Mp_{2n_1}(\mb R)} \phi^{\sharp}(k_2 h_1 g_1) \overline{\psi^{\sharp}(k_2 h_1) } (\pi(\tau(g_1)) u ,  v)  \, d h_1 \, dg_1 \  \, d [k_2].\\
= & \int_{[k_2] \in Mp_{2n_2}(\mb R)/ MP_{n_2-n_1}}  \int \int_{ Mp_{2n_1}(\mb R) \times Mp_{2n_1}(\mb R)} \phi^{\sharp}(k_2 g_1) \overline{\psi^{\sharp}(k_2 h_1) } (\pi(\tau(g_1)) u ,  \pi(\tau(h_1)) v)  \, d h_1 \, dg_1 \  d [k_2]\\
= & \int_{[k_2] \in Mp_{2n_2}(\mb R)/ MP_{n_2-n_1}} (\mc L^0(\phi \otimes u)(k_2), \mc L^0(\psi \otimes v)(k_2)) \, d [k_2]\\
= & (\mc L^0(\phi \otimes u), \mc L^0(\psi \otimes v))_{\Ind_{\,} \pi^{\tau} \otimes \chi^{\alpha}}.
\end{split}
\end{equation}
It follows that $\mc L^0(V(I^{\alpha}) \otimes V(\pi)) \cong I^{\alpha}(\pi).$
$\Box$
\begin{cor}~\label{innerproduct2}
 Let $\pi \in \Pi_u(Mp_{2n_1}(\mb R))$. 
Then for every $u,v \in V(\pi)$ and $\phi, \psi \in C_c^{\infty}(X_0,I^{\alpha})$, 
$$ 2 ( \phi \otimes u,  \psi \otimes v)_{1}=(\mc L^0( \phi \otimes u) , \mc L^0( \psi \otimes v) )_{\Ind_{\,} \pi^{\tau} \otimes \chi^{\alpha}}.$$
Suppose that $I^0(\pi)$ is well-defined. Then for every $\phi, \psi \in C_c^{\infty}(X_0,I^{\alpha}) \oplus V(I^{\alpha})$, 
$$ 2 ( \phi \otimes u,  \psi \otimes v)_{1}=(\mc L^0( \phi \otimes u) , \mc L^0( \psi \otimes v) )_{\Ind_{\,} \pi^{\tau} \otimes \chi^{\alpha}}.$$
\end{cor}
By Theorem \!~\ref{innerproduct}, $I^{\alpha}(\pi)$ can be identified with $\mc L^0( V(I^{\alpha}) \otimes V(\pi))$, which is a subspace of
$$\Ind_{MP_{n_2-n_1,n_1}}^{Mp_{2n_2}(\mb R)} \pi^{\tau} \otimes \chi^{\alpha}.$$
It is also $K$-finite with respect to the maximal compact subgroup $MU(n_2)$.
\subsection{A Density Theorem }~\label{ialphainduced1}
In this section, we will show that 
$$\mc L^0( V(I^{\alpha}) \otimes V(\pi)) = V(\Ind_{MP_{n_2-n_1,n_1}}^{Mp_{2n_2}(\mb R)} \pi^{\tau} \otimes \chi^{\alpha}).$$
\begin{thm}~\label{induction1}
Let $\tau$ and $\pi^{\tau}$ be as in Definition \!~\ref{tauauto}.
Let $\pi \in \Pi_u (Mp_{2n_1}(\mb R))$. Suppose that $I^{0}(\pi)$ is well-defined.
Suppose $n_1 < n_2$ and $\pi(\epsilon) \chi^{\alpha}(\epsilon)=1$.
Let $P_{n_2-n_1,n_1}=Sp_{2n_1}(\mb R) GL(n_2-n_1, \mb R) N$ be as in Lemma \!~\ref{simpleiso}.
Realize $MP_{n_2-n_1,n_1}$ as
the quotient group $$Mp_{2n_1}(\mb R) \times MGL(n_2-n_1)N/\{(1,1), (\epsilon, \epsilon) \}.$$
Then $$I^{\alpha}(\pi) \cong \mc L^0(V(I^{\alpha}) \otimes V(\pi))=V(\Ind_{MP_{n_2-n_1, n_1}}^{Mp_{2n_2}(\mb R)} \pi^{\tau} \otimes \chi^{\alpha}).$$ 
Furthermore, $\Ind_{MP_{n_2-n_1,n_1}}^{Mp_{2n_2}(\mb R)} \pi^{\tau} \otimes \chi^{\alpha}$ is equivalent to the completion of
$I^{\alpha}(\pi)$ with respect to the inner product $(\ , \ )_1$. 
\end{thm}
By Theorems \!~\ref{inte} and \!~\ref{innerproduct}, $ I^{\alpha}(\pi) \cong \mc L^0(V(I^{\alpha}) \otimes V(\pi))$ is a subspace of $V(\Ind_{MP_{n_2-n_1,n_1}}^{Mp_{2n_2}(\mb R)} \pi^{\tau} \otimes \chi^{\alpha})$ and $ \mc L^0(C_c^{\infty}(X_0,I^{\alpha}) \otimes V(\pi))$ is a dense subspace of $\Ind_{MP_{n_2-n_1,n_1}}^{Mp_{2n_2}(\mb R)} \pi^{\tau} \otimes \chi^{\alpha}$.
Since  every element in $ C_c^{\infty}(X_0, I^{\alpha}) \otimes V(\pi)$ can be approximated by vectors in $V(I^{\alpha}) \otimes V(\pi)$, it is reasonable to believe that $\mc L^0(V(I^{\alpha}) \otimes V(\pi))$ is dense in $\Ind_{MP_{n_2-n_1,n_1}}^{Mp_{2n_2}(\mb R)} \pi^{\tau} \otimes \chi^{\alpha}.$  This may be a direct consequence of some theorem in functional analysis about the averaging operator 
$\mc L^0$. Unfortunately, we are not aware of any such theorem. We are forced to go back and dig into the compact picture in which the Harish-Chandra module is visible.\\
\\
{\bf Proof of Theorem \!~\ref{induction1}}: It suffices to prove that $\mc L^0(V(I^{\alpha}) \otimes V(\pi))$ is dense in
$$\Ind_{MP_{n_2-n_1,n_1}}^{Mp_{2n_2}(\mb R)}  \pi^{\tau} \otimes \chi^{\alpha}.$$
We shall prove that every element in $\mc L^0( C_c^{\infty}(X_0, I^{\alpha}) \otimes V(\pi))$ can be approximated
by a sequence of elements in $\mc L^0(V(I^{\alpha}) \otimes V(\pi))$. \\
\\
Consider $ \psi \otimes u \in  C_c^{\infty}(X_0, I^{\alpha}) \otimes V(\pi)$. Regard $\psi$ as a smooth function on
$MU(n_1+n_2)$. By the Stone-Weierstrass Theorem, there exists a sequence of functions $\phi_i \in V(I^{\alpha})$, such that
$\phi_i \rightarrow \psi$ under the sup norm. In other words, for $i$ sufficiently large,
$$ |\phi_i(k)- \psi(k)| \leq \delta 1_{MU(n_1+n_2)}(k) \qquad (\, \forall \, \, k \in K).$$
By  Theorem \!~\ref{innerproduct} and Corollary \!~\ref{innerproduct2},
$$(\mc L^0((\phi_i-\psi) \otimes u), \mc L^0( (\phi_i-\psi) \otimes u))_{\Ind \pi^{\tau} \otimes \chi^{\alpha}}=2 ((\phi_i-\psi) \otimes u,  (\phi_i-\psi) \otimes u)_1$$
The latter is no greater that $ 2 \delta^2(  1_{MU(n_1+n_2)} \otimes u,  1_{MU(n_1+n_2)} \otimes u)_1$ with respect to $I^0 \otimes \pi$.
Therefore, $\mc L^0( \phi_i \otimes u) \rightarrow \mc L^0( \psi \otimes u)$ under the Hilbert norm. By Theorem \!~\ref{inte}, $\mc L^0(C_c^{\infty}(X_0, I^{\alpha}) \otimes V(I^{\alpha}))$ is already dense in $\Ind_{MP_{n_2-n_1,n_1}}^{Mp_{2n_2}(\mb R)} \pi^{\tau} \otimes \chi^{\alpha}.$ So $\mc L^0(V(I^{\alpha}) \otimes V(\pi))$ is dense in
$\Ind_{MP_{n_2-n_1,n_1}}^{Mp_{2n_2}(\mb R)} \pi^{\tau} \otimes \chi^{\alpha}.$ \\
\\
By an easy $K$-finiteness argument with respect to $MU(n_2)$,
$$\mc L^0( V(I^{\alpha}) \otimes V(\pi))= V(\Ind_{MP_{n_2-n_1,n_1}}^{Mp_{2n_2}(\mb R)} \pi^{\tau} \otimes \chi^{\alpha}).$$
We have proved that
$$I^{\alpha}(\pi) \cong \Ind_{MP_{n_2-n_1,n_1}}^{Mp_{2n_2}(\mb R)} \pi^{\tau} \otimes \chi^{\alpha}.$$
$\Box$\\
\\
Combining this with Theorem \!~\ref{induction}, we obtain 
\begin{thm}~\label{11} Assume that $\alpha \equiv n_1+n_2+1 \pmod 2$. If $\pi \in \Pi_u(Mp_{2n_1}(\mb R))$ and $I^0(\pi)$ is well-defined, then
$$\Ind_{MP_{n_2-n_1,n_1}}^{Mp_{2n_2}(\mb R)} \pi^{\tau} \otimes \chi^{\alpha}
\cong \bigoplus_{ p+q=n_1+n_2+1, \ p-q \equiv \alpha \pmod 4} \mc Q(2n_1;p,q;2n_2)(\pi)$$
as unitary representations.
\end{thm}

\section{Wave Front Set and the Nonvanishing Theorem}
\begin{thm}[\!~\cite{bv} \!~\cite{vogan01}]~\label{voganlemma}
Let $P=MAN$ be a parabolic subgroup of a reductive Lie group $G$.
Let $\sigma$ be an irreducible representation of $L=MA$. Then
$$WF(\Ind_{P}^{G} \sigma)= cl(\Ind_{\f l}^{\f g} WF(\sigma)).$$
\end{thm}
Notice that $ \Ind_{\f l}^{\f g} (WF(\sigma))$ may contain several irreducible components of the same dimension. 
Applying Theorem \!~\ref{voganlemma} to $I^{\alpha}$ for $Mp_{2n}(\mb R)$, we obtain
\begin{cor}
The wave front set $WF(I^{\alpha})=\Ind_{\f{gl}(n,\mb R)}^{{\f{sp}}_{2n}(\mb R)} \mc O_{[1^n]}$.
\end{cor}
Again assume $p+q=n+1$. By the theorem of Kudla-Rallis-Lee-Zhu, we have
$$I^{\alpha}=\bigoplus_{p+q=n+1, p-q \equiv \alpha \pmod 4} \theta(p,q;2n)(\trivial).$$
Consider the wave front sets. 
From Corollary \!~\ref{induced0}, $\Ind_{\f{gl}(n,\mb R)}^{{\f{sp}}_{2n}(\mb R)} \mc O_{[1^n]}$ consists of $n+1$ components $\bold D^{(j)}$ and is equal to the closure of $\mc O_{[2^n]} \cap {\f{sp}}_{2n}(\mb R)$. We have
$$\text{\Large cl}(\bigcup D^{(j)})= \bigcup_{p+q=n+1, p-q \equiv \alpha \pmod 4} WF(\theta(p,q;2n)(\trivial)).$$
This implies that these
$WF(\theta(p,q;2n)(\trivial))$ must be in the closure of $\mc O_{[2^n]} \cap {\f{sp}}_{2n}(\mb R)$ (\!~\cite{bv0}). In fact, they are distinct.
\begin{thm}~\label{twocom}
Assume $p+q=n+1$.
Let $[2^n]^{(i)}$ be the signed Young diagram with $i$ rows starting with $+$. Then
$$WF(\theta(p,q;2n)(\trivial)) = cl(\mc O_{[2^n]^{(p)}} \cup \mc O_{[2^n]^{(p-1)}}).$$
Here we take $\mc O_{[2^n]^{(-1)}}$ and $\mc O_{[2^n]^{(n+1)}}$ to be the empty set.
In particular,
$$WF(I^{\alpha})=\bigcup_{p+q=n+1, p-q \equiv \alpha \pmod 4} WF(\theta(p,q;2n)(\trivial))$$
\end{thm}
Proof: From Lemma \!~\ref{panlemma}, only 
$cl(\mc O_{[2^n]^{(p)}} \cup \mc O_{[2^n]^{(p-1)}})$
can occur in the image of the moment map $m_2$ for $(O(p,q), Sp_{2n}(\mb R))$. Thus
by Theorem \!~\ref{pr2.8} of Przebinda,
$$WF(\theta(p,q;2n)(\trivial)) \subseteq cl(\mc O_{[2^n]^{(p)}} \cup \mc O_{[2^n]^{(p-1)}}).$$
But by the decomposition Theorem \!~\ref{krlz},
$$\bigcup_{j=0}^n cl(\mc O_{[2^n]}^{(j)}) = WF(I^{\alpha})=\bigcup_{p+q=n+1, p-q \equiv \alpha \pmod 4} WF(\theta(p,q;2n)(\trivial)) \subseteq \bigcup_{j=0}^{n} cl(\mc O_{[2^n]^{(j)}})$$
We must have
$WF(\theta(p,q;2n)(\trivial)) = cl(\mc O_{[2^n]^{(p)}} \cup \mc O_{[2^n]^{(p-1)}}).$
$\Box$ \\
\\
This theorem partly explains the decomposition theorem of Kudla-Rallis-Lee-Zhu in terms of the orbit philosophy. Generalizing this idea to $I^{\alpha}(\pi)$, we give a proof of the Generic Nonvanishing Theorem.\\
\\
Proof of Theorem \!~\ref{non2}: 
By Lemma \!~\ref{1k} (3), every $MU(n_1+n_2)$-finite matrix coefficient of $I^{\alpha}$ is bounded by a multiple of $$(I^{0}(g) 1_{MU(n_1+n_2)}, 1_{MU(n_1+n_2)}),$$ which is a spherical function. 
Notice that the infinitesimal character of $I^0$ is equal to
$$\mu=\overbrace{(\frac{n_1+n_2-1}{2}, \frac{n_1+n_2-1}{2}, \frac{n_1+n_2-3}{2}, \frac{n_1+n_2-3}{2}, \ldots, \frac{n_1+n_2-1}{2}-[\frac{n_1+n_2-1}{2}])}^{n_1+n_2},$$
and $$\mu-\rho(Mp_{2n_1+2n_1}(\mb R)) \prec \overbrace{(\frac{-n_1-n_2}{2}, \frac{-n_1-n_2}{2}+1, \ldots, \frac{-n_1-n_2}{2}+[\frac{n_1+n_2}{2}], 0, \ldots, 0)}^{n_1+n_2}.$$
By Theorem \!~\ref{mat} or by Harish-Chandra's spherical expansion, $(I^{0}(g) 1_{MU(n_1+n_2)}, 1_{MU(n_1+n_2)})$ is
bounded by a multiple of $$a(g)^{(\frac{-n_1-n_2}{2}, \frac{-n_1-n_2}{2}+1, \ldots, \frac{-n_1-n_2}{2}+[\frac{n_1+n_2}{2}], 0, \ldots, 0)}.$$
In particular, $(I^{0}(g_1) 1_{MU(n_1+n_2)}, 1_{MU(n_1+n_2)})$ is
bounded by a multiple of $$a(g_1)^{(\frac{-n_1-n_2}{2}, \frac{-n_1-n_2}{2}+1, \ldots, \frac{n_1-n_2}{2}+1)} \qquad (g_1 \in Mp_{2n_1}(\mb R).$$
Since $\pi \in R_{ss}(Mp_{2n_1}(\mb R), \omega(n_1+n_2+1,0; 2n_1))$ (Def. \!~\ref{qi}), for every leading exponent $v$ of $\pi$,
$$\Re(v) \preceq \bold{\frac{n_1+n_2+1}{2}}-\bold{n_1}-\bold{1} - \rho(Sp_{2n_1}(\mb R)).$$
We obtain
 \begin{equation}~\label{welldefine}
 \begin{split}
 & \Re(v)-(\frac{n_1+n_2}{2}, \frac{n_1+n_2}{2}-1, \ldots, \frac{n_2-n_1}{2}+1) + 2 \rho(Mp_{2n_1}(\mb R))\\
 \preceq & \bold{\frac{n_2-n_1-1}{2}}- \rho(Mp_{2n_1}(\mb R))-(\frac{n_1+n_2}{2}, \frac{n_1+n_2}{2}-1, \ldots, \frac{n_2-n_1}{2}+1) + 2 \rho(Mp_{2n_1}(\mb R)) \\
 \preceq & \bold{\frac{n_2-n_1-1}{2}}-(\frac{n_1+n_2}{2}, \frac{n_1+n_2}{2}-1, \ldots, \frac{n_2-n_1}{2}+1)+(n_1, n_1-1, \ldots 1)\\
 \prec  & 0.
 \end{split}
\end{equation}
Inequality \!~\ref{welldefine} and Theorem \!~\ref{mat} imply that $(\ ,\  )_1$ is well-defined for $V(I^{0} \otimes V(\pi))$.
So $I^{\alpha}(\pi)$ is well-defined for any $\alpha$. Suppose that $\alpha \equiv n_1+n_2+1 \pmod 2$. Theorem \!~\ref{11} applies. We have 
$$\Ind_{MP_{n_2-n_1,n_1}}^{Mp_{2n_2}(\mb R)} \pi^{\tau} \otimes \chi^{\alpha} =\bigoplus_{p+q=n_1+n_2+1, p-q \equiv \alpha \pmod 4} \mc Q(2n_1;p,q;2n_2)(\pi).$$
By Theorem \!~\ref{quantuminduction3},
$$\Ind_{MP_{n_2-n_1,n_1}}^{Mp_{2n_2}(\mb R)} \pi^{\tau} \otimes \chi^{\alpha}=\bigoplus_{p+q=n_1+n_2+1, p-q \equiv \alpha \pmod 4} \theta_s(p,q;2n_2) \theta_s(2n_1;p,q)(\pi).$$
It follows that
$$ \Ind_{\f{sp}_{n_1}(\mb R) \f{gl}(n_2-n_1)}^{\f{sp}_{2n_2}(\mb R)} WF(\pi^{\tau}) 
= \bigcup_{p+q=n_1+n_2+1, p-q \equiv \alpha \pmod 4} WF(\theta_s(p,q;2n_2) \theta_s(2n_1;p,q) (\pi)) .$$
Let $\mc O_{\bold D}$ be a nilpotent orbit of maximal dimension in $$ \Ind_{\f{sp}_{2n_1}(\mb R) \f{gl}(n_2-n_1)}^{\f{sp}_{2n_2}(\mb R)} WF(\pi^{\tau}).$$
Define $S \subseteq \{(p, q) \mid p+q=n_1+n_2+1,  p \ \text{fixed parity} \}$ to be the set of those $(p,q)$ such that
$\mc O_{\bold D}$ occurs in
the image of the moment map $m_2$ associated with $(O(p,q), Sp_{2n_2}(\mb R))$ (see Definition \!~\ref{moment}).
By Theorem \!~\ref{pr2.8}, $\mc O_{\bold D}$ can only appear in $WF(\theta_s(p,q;2n_1) \theta_s(2n_1;p,q) (\pi))$ with $(p,q) \in S$. There must exist at least a $(p_0,q_0) \in S$ such that
$$\mc O_{\bold D} \subseteq WF(\theta_s(p_0,q_0;2n_2) \theta_s(2n_1;p_0,q_0)  (\pi)).$$
Hence $\mc Q(2n_1;p_0,q_0;2n_2)(\pi) \neq 0$ and $\theta_s(2n_1;p_0,q_0)(\pi) \neq 0$.
$\Box$ \\
\\
Under certain favorable circumstances, $\mc O_{\bold D}$ only occurs in the image of $m_2$ for a unique $(p,q)$ with $p+q=n_1+n_2+1$. In this case, we are able to determine for which $(p,q)$,
$\theta_s(2n_1;p,q)(\pi) \neq 0$. 
\begin{cor}~\label{non20} Consider the group $Mp_{2n_1+2n_2}(\mb R)$ with $n_1 < n_2$.
Let $\pi$ be a unitary representation in $$R_{ss}(Mp_{2n_1}(\mb R), \omega(n_1+n_2+1,0; 2n_1)),$$
(see Def. \!~\ref{qi}).
Let $\mc O_{\bold D}$ be a nilpotent orbit of maximal dimension in 
$$\Ind_{\f{sp}_{2n_1}(\mb R) \oplus \f{gl}(n_2-n_1, \mb R)}^{\f{sp}_{2n_2}(\mb R)} WF(\pi^{\tau}).$$ 
Fix a parity of $p$. Suppose that $\mc O_{\bold D}$ is not contained
the image of the moment map $m_2$ associated with $(O(p,n_1+n_2+1-p), Sp_{2n_2}(\mb R))$ for $p \neq p_0$.
Then 
$\mc Q(2n_1; p_0, n_1+n_2+1-p_0; 2n_2)(\pi) \neq 0$ and
$\theta_s(2n_1;p_0, n_1+n_2+1-p_0)(\pi) \neq 0$.
\end{cor}

\chapter{Construction of Unipotent Representations, Unitarity and Infinitesimal Character}\label{construct}
Start with a real nilpotent orbit 
$\mc O=\mc O_{\bold D}$ in $\mc U(O(p,q))$ or $\mc U(Mp_{2n}(\mb R))$ (see Definition \!~\ref{u}). First, construct an alternating sequence of nilpotent orbits of metaplectic groups and nilpotent orbits of orthogonal groups:
$$\mc O(d_1)=\mc O_{\bold D}, \mc O(d_1-1)=\mc O_{\bold D-\bold 1}, \ldots, \mc O(1)=
\mc O_{\bold D-\bold d_1+\bold 1}.$$
Write
$$G(1)=G(\mc O(1)), G(2)=G(\mc O(2)), \ldots, G(d_1)=G(\mc O(d_1)).$$
Clearly, $\mc O(k) \in \mc U$ for every $k$ (see Cor. \!~\ref{inu}). Let $\mc O^{\prime}$ be the negative of $\mc O$
$$\mc O^{\prime}=\{ -x \mid x \in \mc O \}.$$
$\mc O^{\prime}$ may or may not equal to $\mc O$.
\begin{defn}~\label{unip}
We define $\mc N(\mc O_{\bold D})$ inductively.
\begin{enumerate} 
\item For $\mc O(1)$, let $\mc N(\mc O(1))$ be the set of all one-dimensional unitary characters of $G(1)$. 
\item
Suppose $\mc N(\mc O(k))$ is defined.
\begin{enumerate}
\item For $G(k+1)$ orthogonal, define
$\mc N_0(\mc O(k+1))$ to be
$$\{\pi(k+1) \otimes \eta
 \mid \pi(k+1)=\theta_s(G(k), G(k+1))(\pi), \pi \in \mc N(\mc O(k)), \eta {\text \ a \ character \ of \ } O(p,q) \}.$$
Let $\mc N(\mc O(k+1))= \mc N_0(\mc O(k+1)) \cup \{ \pi^* \mid \pi \in \mc N_0(\mc O(k+1)^{\prime}) \}.$
\item For $G(k+1)$ metaplectic, define
$\mc N_0(\mc O(k+1))$ to be 
$$\{\pi(k+1) \mid \pi(k+1)=\theta_s(G(k), G(k+1))(\pi), \pi \in \mc N(\mc O(k))\}.$$
Let $\mc N_1(\mc O(k+1))=\mc N_0(\mc O(k+1)) \cup \{\pi^{\tau} \mid \pi \in \mc N_0(\tau(\mc O(k+1))) \}$ and
$$\mc N(\mc O(k+1))=\mc N_1(\mc O(k+1)) \cup \{ \pi^* \mid \pi \in \mc N_1(\mc O(k+1)^{\prime}) \}.$$
\end{enumerate}
\end{enumerate}
\end{defn} 
Any of the operations in the definition of $\mc N$, tensoring with a character, involution by $\tau$ or $*$, do not alter the infinitesimal character, the real part of the leading exponent, the unitarity, or the associated variety of the annihilator.
In most cases, we will ignore these operations. The reader should be aware of the fact that $*$ and $\tau$ operations alter the wave front set by the involution $\prime$ and $\tau$. The author would like to thank the referee for suggesting this definition.
\begin{thm}~\label{exi} For $\mc O_{\bold D} \in \mc U(G)$,
the set $\mc N(\mc O_{\bold D})$ is not empty.
\end{thm}
For technical reasons, we will postpone the proof till the next chapter. In this chapter, we work under the assumption that $\mc N(\mc O_{\bold D}) \neq \emptyset$. We will show that
\begin{thm}~\label{un}
Let $\mc O_{\bold D}$ be in $\mc U(O(p,q))$ or $\mc U(Sp_{2n}(\mb R))$. Then any $\pi(k) \in \mc N(\mc O(k))$ is in $R_{ss}(G(k), \omega)$ (Def. \!~\ref{qi}) for the pair $(G(k), G(k+1))$. The set $\mc N(\mc O_{\bold D})$ is well-defined. The representations in $\mc N(\mc O_{\bold D})$ are all unitary.
\end{thm}
Thus 
$\mc N(\mc O_{\bold D}) \subset \Pi_{u}(G(\mc O_{\bold D})).$ We will further
determine the infinitesimal character of $\pi \in \mc N(\mc O_{\bold D})$. It turns out that $\mc I(\pi)$ only depends on $\bold d$ and consists of $\rho$-like segments. 
\section{Unitarity of $\mc N(\mc O)$}
Let $\mc O_{\bold D}$ be in $\mc U(O(p,q))$ or $\mc U(Sp_{2n}(\mb R))$. 
The group $G(\mc O_{\bold D})$ can be read off from $D^+$, the number of positive boxes in $\bold D$, and from $D^-$, the number of negative boxes in $\bold D$. Construct an alternating sequence of symplectic signed Young diagrams and orthogonal signed Young diagrams:
$$\bold D(d_1)=\bold D, \bold D(d_1-1)=\bold D(d_1)-\bold 1, \ldots, \bold D(2)=\bold D(3)-\bold 1, \bold D(1)=\bold D(2)-\bold 1.$$
Let $\bold D(0)=0$. Recall that $\| \bold  d \|$ stands for the number of boxes in $\bold d$. We have the following.
\begin{lem}~\label{pm}
Let $\mc O_{\bold D} \in \mc U$. Then
\begin{enumerate}
\item  For every $k \in [1,d_1-1]$,
$$D(k+1)^+ - D(k)^+ \geq D(k)^- - D(k-1)^-, \qquad D(k+1)^- - D(k)^- \geq D(k)^+ - D(k-1)^+.$$
Thus,
$$ D(k+1)^+ + D(k-1)^- \geq || \bold d(k) ||, \qquad D(k+1)^-+ D(k-1)^+ \geq || \bold d(k) ||.$$
\item For every $k \in [1, d_1-1]$,
$$ D(k+1)^+ \geq D(k)^- ,\qquad D(k+1)^- \geq D(k)^+.$$
\item  For every $k \in [1, d_1-1]$,
$$ \| \bold d(k+1) \|+\| \bold d(k-1) \| \geq  2\| \bold d(k) \|+2$$
if $\bold D(k)$ is symplectic;
\item 
for every $k \in [1, d_1-1]$,
$$ \| \bold d(k+1) \|+\| \bold d(k-1) \| \geq  2\| \bold d(k) \|$$
if $\bold D(k)$ is orthogonal.
\item $\| \bold d(k+2) \| \equiv \|\bold d(k) \| \, \pmod 2 $ for all $k>0$.
\end{enumerate}
\end{lem}
Proof: All the statements can be read off from $\bold D$. The number of $+$ boxes on the $k+1$-th column must be greater than or equal to the number of $-$ boxes on the $k$-th column. We obtain (1).
Adding all columns together, we obtain (2). (3) and (4) are consequences of (2.a)and (3.a) from Definition \!~\ref{u}. (5) follows from (2.c) and (3.c) from Definition \!~\ref{u}. $\Box$\\
\\
Proof of Theorem \!~\ref{un}: Let $\mc O_{\bold D} \in \mc U$. Consider the sequences $\bold D(k)$ and $\mc O(k)$. We shall prove that $\pi(k)$ is unitary and lies in $R_{ss}(G(k), \omega))$ for the pair $(G(k), G(k+1))$. We will use induction.
\begin{enumerate}
\item For $k=1$, by definition, every $\pi(1)$ in $\mc N(\mc O(1))$ is unitary. 
If $G(1)$ is metaplectic, then we have 
\begin{itemize}
\item $D(1)^+=D(1)^-$; 
\item $ D(2)^+ \geq 2 D(1)^-,  D(2)^- \geq 2D(1)^+,$ by Lemma \!~\ref{pm} (1);
\item $D(2)^+ + D(2)^- \geq 4 D(1)^{-}+2,$ by Lemma \!~\ref{pm} (3).
\end{itemize}
 if $G(1)$ is orthogonal, then we have 
$D(2)^+ = D(2)^- \geq  D(1)^{-}+ D(1)^+.$
It follows from the Definition \!~\ref{qi}, that every unitary representation is in $R_{ss}(G(1), \omega)$ with respect to the dual pair $(G(1), G(2))$.
\item Suppose that $\mc N(\mc O(k-1))$ is well-defined and $\mc N(\mc O(k-1)) \subseteq
\Pi_u(G(k-1))$.  Suppose that $\mc N(\mc O(k-1)) \subseteq R_{ss}(G(k-1), \omega)$ with respect to the pair $(G(k-1), G(k))$. This is our induction hypothesis. \\
\\
If $G(k)$ is metaplectic, let
$$p^{\p}=D(k+1)^+, q^{\p}=D(k+1)^-, n= D(k)^+=D(k)^-, p= D(k-1)^+, q=D(k-1)^-.$$
Then by Lemma \!~\ref{pm}, $p+q \equiv p^{\p}+q^{\p} \pmod 2$, $p^{\p}+q^{\p} -2n \geq 2n-p-q+2$, $p^{\p} \geq n$ and $q^{\p} \geq n$.  Clearly $\max(p^{\p},q^{\p}) > n$. If $n \neq \min(p^{\p},q^{\p})$, then Theorem \!~\ref{quantuminduction1} holds and $\pi(k)=\theta_s(\pi(k-1))$ is well-defined and unitary. Furthermore,
$\pi(k) \in R_{ss}(G(k), \omega)$ with respect to $(G(k), G(k+1))$. If $n=\min(p^{\p},q^{\p})$, without loss of generality, assume, $n=p^{\p}$. We shall prove that this case does not occur if $\mc O_{\bold D} \in \mc U$.\\
\\
Since $D(k)^{-} = n=p^{\p}=D(k+1)^{+}$, the diagram $\bold D(k+1)$ must have its first column marked with $-$. By the structure of signed Young diagram, $\bold D(k)$ must have its first column marked with $+$. 
$\bold D(k)$ must have its second column marked with $-$. Furthermore, the first column of $\bold D(k)$ and the second column of $\bold D(k)$ must have the same length. Otherwise, there are at least two rows of length one in $\bold D(k)$. At least one $\fm$ would occur in the first column of $\bold D(k)$ according to our definition of $\mc YD_{-}(n,n)$. This
contradicts the fact that the first column of $\bold D(k)$ is all marked  with $+$.  Thus,
$$\| \bold d(k)\|- \| \bold d(k-1) \|=\| \bold d(k-1) \|-\|\bold d(k-2) \|.$$
By the same argument, $\bold D(k)$ must be of the following shape and sign pattern: 
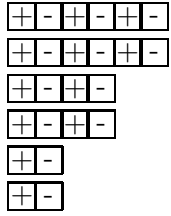
\begin{center}
\fp \fm \fp \fm \fp  \fm \\
\fp   \fm \fp \fm \fp \fm  \\
\fp  \fm \fp \fm  \mbo \mbo \\
\fp \fm  \fp \fm  \mbo  \mbo \\
\fp \fm \mbo  \mbo  \mbo  \mbo \\
\fp \fm \mbo \mbo \mbo \mbo 
\end{center}
In other words, the $j$-th column of $\bold D(k)$ must be marked by $(-1)^{j}$
and the $2j+1$-th column of $\bold D(k)$ must be of the same length as  the
$2j+2$-th column. This kind of $\bold D(k)$ is excluded from $\mc U$. So the case
$n=\min(p^{\p},q^{\p})$ does not occur. \\
\\
If $G(k)$ is orthogonal, let
$$n^{\p}=D(k+1)^+=D(k+1)^-, \qquad p= D(k)^+, q=D(k)^-, \qquad n=D(k-1)^+=D(k-1)^-.$$
Then by Lemma \!~\ref{pm}, we have
$$2n^{\p}-p-q \geq p+q-2n > p+q-2n-2, \qquad p \geq n, \qquad q \geq n.$$
By the same argument as for  metaplectic $G(k)$, $n \neq \min(p,q)$. By Theorem \!~\ref{quantuminduction2}, $\pi(k)=\theta_s(\pi(k-1))$ is well-defined, unitary and $\pi(k) \in R_{ss}(G(k), \omega)$ with respect to the pair $(G(k), G(k+1))$.
\item It follows that  $\mc N(\mc O(k)) \subseteq R_{ss}(G(k), \omega)$ for the pair $(G(k), G(k+1))$ and
every $\pi(k) \in \mc N(\mc O(k))$ is unitary.
\end{enumerate}
$\Box$
\section{Infinitesimal Character}
\begin{thm}[Thm 1.19, \!~\cite{pr}] Let $\mc I(\pi)$ be the infinitesimal character of $\pi$. Suppose that $\pi$ occurs in the theta correspondence with respect to $(Sp_{2n}(\mb R),O(p,q))$.
\begin{enumerate}
\item Suppose $p+q < 2n$. Then $\mc I(\theta(p,q;2n)(\pi))$ can be obtained by augmenting $\mc I(\pi)$ by
$$(n-\frac{p+q}{2}, n-\frac{p+q}{2}-1, \ldots, 1+[\frac{p+q}{2}]-\frac{p+q}{2})$$
\item Suppose $2n+1 < p+q$. Then $\mc I(\theta(2n;p,q)(\pi))$ can be obtained by augmenting $\mc I(\pi)$ by
$$(\frac{p+q}{2}-n-1, \frac{p+q}{2}-n-2,\ldots, \frac{p+q}{2}-[\frac{p+q}{2}])$$
\item Suppose $p+q=2n$ or $p+q=2n+1$. Then $\mc I(\theta(p,q;2n)(\pi)
)$ is just $\mc I(\pi)$.
\end{enumerate}
\end{thm}
\begin{no}~\label{ipisegment} We define the orthogonal segment
$$\mc I_+(m)=\overbrace{(\frac{m}{2}-1, \frac{m}{2}-2, \ldots, \frac{m}{2}-[\frac{m}{2}])}^{[\frac{m}{2}]}$$
and the symplectic segment
$$\mc I_-(m)=\overbrace{(\frac{m}{2}, \frac{m}{2}-2, \ldots, \frac{m}{2}+1-[\frac{m+1}{2}])}^{[\frac{m+1}{2}]}.$$
\end{no}
The orthogonal segment $\mc I_+(m)$ is just $\rho(\f o(m, \mb C))$. For $m$ even, $\mc I_{-}(m)$ is $\rho(\f{sp}_{m}(\mb C))$. For $m$ odd, $\mc I_{-}(m)$ is the infinitesimal character of the oscillator representation of $Mp_{m+1}(\mb R)$.
\begin{thm}~\label{infi} Let $\mc O_{\bold D}$ be in $\mc U(O(p,q))$ or $\mc U(Sp_{2n}(\mb R))$. Let $\pi \in \mc N(\mc O_{\bold D})$ and
$$\bold d^t=(m_1 \geq m_2 \geq \ldots \geq m_{d_1}).$$
If $G(\bold D)$ is an orthogonal group, then 
$$\mc I(\pi)=(\mc I_{+}(m_1), \mc I_{-}(m_2), \mc I_{+}(m_3), \mc I_{-}(m_4), \ldots).$$
If $G(\bold D)$ is a symplectic group, then 
$$\mc I(\pi)=(\mc I_{-}(m_1), \mc I_{+}(m_2), \mc I_{-}(m_3), \mc I_{+}(m_4), \ldots).$$
\end{thm}
Proof: Let $\pi \in \mc N(\mc O_{\bold D})$. We prove this theorem by induction on $d_1$.  If $d_1=1$, by definition of $\mc N(\mc O(1))$, $\pi$ restricted to the identity component of $G(\mc O_{\bold D})$ must contain a trivial constituent.
Thus $\mc I(\pi)=\mc I_{\pm}(m_1)$. 
Suppose our assertion holds for any $\bold D$ with $d_1 \leq k$.
Let $d_1=k+1$. \\
\\
If $G(\bold D)$ is $O(p,q)$, then according to our definition of $\pi$,
$\pi|_{SO_0(p,q)}$ must be equivalent to $$\theta(G(\mc O_{\bold D-1}), G(\mc O_{\bold D}))(\sigma)|_{SO_0(p,q)}$$ or its contragredient for some
$\sigma \in \mc N(\mc O_{\bold D-1})$. 
$\mc I(\pi)$ can be obtained by augmenting $\mc I(\sigma)$ 
with
$$\overbrace{(\frac{m_1}{2}-1, \frac{m_1}{2}-2, \ldots, \frac{m_1}{2}-[\frac{m_1}{2}])}^{[\frac{m_1}{2}]}=\mc I_+(m_1).$$
If $G(\bold D)$ is $Mp_{2n}(\mb R)$,  observe that
$\mc I(\pi^{\tau})=\mc I(\pi)$. Then 
$$\mc I(\pi)=\mc I(\theta(G(\mc O_{\bold D-1}), G(\mc O_{\bold D}))(\sigma))$$
for some $\sigma \in \mc N(\mc O_{\bold D-1})$. $\mc I(\pi)$ can be obtained by augmenting $\mc I(\sigma)$ 
with
$$\overbrace{(\frac{m_1}{2}, \frac{m_1}{2}-1, \ldots, \frac{m_1}{2}-[\frac{m_1-1}{2}])}^{[\frac{m_1+1}{2}]}=\mc I_{-}(m_1).$$
$\Box$ \\
\begin{no}~\label{infd}
Let $$\bold d^t=(m_1 \geq m_2 \geq \ldots \geq m_{d_1}).$$
For each orthogonal Young diagram $\bold d$, we define
$\mc I_{+}(\bold d)$ 
to be
$$(\mc I_{+}(m_1), \mc I_{-}(m_2), \mc I_{+}(m_3), \mc I_{-}(m_4), \ldots).$$
For each symplectic Young diagram $\bold d$, we define $\mc I_{-}(\bold d)$ to be
$$(\mc I_{-}(m_1), \mc I_{+}(m_2), \mc I_{-}(m_3), \mc I_{+}(m_4), \ldots).$$
\end{no}
Theorem \!~\ref{infi} states that $\mc I(\pi)=\mc I(\bold d)$ for $\pi \in \mc N(\mc O_{\bold D})$. Thus $\mc I(\pi)$ only depends on $\bold d$, not on the signs of $\bold D$. 
\section{$\mc I_{\pm}(\bold d)$ and $\bold d$: An algorithm}
In this section, we will seek a direct way of obtaining $\mc I_{\pm}(\bold d)$ from the  Young diagram $\bold d$. We assume 
$\bold d^t$ is very odd or very even. The Young diagram $\bold d$ does not have to be pre-rigid or in $\mc U$. The algorithm can be described as follows.
\begin{enumerate}
\item First, cover the Young diagram $\bold d$ by horizontal and vertical dominoes
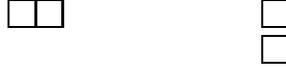
\begin{center}
\fb \fb \hspace{1.in} \fb \\
\mbo \mbo \hspace{1.in} \fb \\
\end{center}
as follows. We cover each column of the Young diagram $\bold d$ by consecutive vertical dominoes, starting from the bottom row. If $\bold d^t$ is very even, vertical dominoes covers $\bold d$ completely. If $\bold d^t$ is very odd, we have the first row left uncovered. We cover the first row of the Young diagram $\bold d$ by
consecutive horizontal dominoes, starting from the {\it right}. We may have the leftmost block uncovered. In that case, cover it with a horizontal domino anyway. We call this domino an open domino.
\item  For each vertical domino $DO$, we can enumerate the number of dominoes above it. A horizontal domino will be counted as $\frac{1}{2}$ domino and a vertical domino will be counted as a full domino. Thus we obtain a number $n(DO)$.
\item If $\bold d^t$ is odd, fill the open domino with no number and fill the other horizontal dominoes with $\frac{1}{2}$. For $\mc I_{+}(\bold d)$, if a vertical domino $DO$ is in the $k$-th column, we fill in $DO$ with the number $n(DO)+\frac{1+(-1)^k}{2}$. For $\mc I_{-}(\bold d)$, if a domino $DO$ is in the $k$-th column, we fill in $DO$ with the number $n(DO)+\frac{1-(-1)^k}{2}$.  Extracting  the numbers in all the dominoes, we obtain $\mc I(\bold d)$. 
\end{enumerate}
Let us consider the example. Let $\mc O_{\mathbf D}$ be a real nilpotent coadjoint orbit of $Sp_{36}(\mathbb R)$, with $\bold D=$
\clearpage
\begin{center}
\fm \fp \fm \fp \fm \fp  \\
\fm \fp \fm \fp \fm \mbo \\
\fp \fm \fp \fm \fp \mbo \\
\fm \fp \fm \fp \mbo \mbo  \\
\fm \fp \fm \fp \mbo \mbo  \\
\fm \fp \fm \mbo \mbo \mbo \\
\fp \fm \fp \mbo \mbo \mbo \\
\fm \fp \mbo \mbo \mbo \mbo \\
\fm \fp \mbo \mbo \mbo \mbo \\
\fm \mbo \mbo \mbo \mbo \mbo \\
\fp \mbo \mbo \mbo \mbo \mbo 
\end{center}
We will cover the first column with 5 vertical dominos, 2nd column with 4 vertical dominos, 3rd column with 3 vertical dominos, 4th column with 2 vertical dominos and 5th column with 1 vertical dominos. Then there is still the first row left out. We cover the first row with 3 horizontal dominos. Since $\mathbf d^t$ is odd, we fill in the horizontal dominos with $\{\frac{1}{2}, \frac{1}{2}, \frac{1}{2} \}$. Then for the vertical dominos, we fill the first column with
$$\{\frac{1}{2}, \frac{3}{2}, \frac{5}{2}, \frac{7}{2}, \frac{9}{2} \};$$
the 2nd column with
$$\{\frac{1}{2}, \frac{3}{2}, \frac{5}{2}, \frac{7}{2} \};$$
the 3rd column with
$$\{\frac{1}{2}, \frac{3}{2}, \frac{5}{2} \};$$
and the 4th column with 
$\{\frac{1}{2}, \frac{3}{2} \};$ and the fifth column with $\{\frac{1}{2} \}$.  Finally, we need to add $\frac{1-(-1)^1}{2}=1$ to the first column and 0 to the second column and so on. We obtain
the infinitesimal character for $\mc I_{-}({\mathbf d})$:
$$\{ \frac{1}{2}, \frac{1}{2}, \frac{1}{2}, \frac{3}{2}, \frac{5}{2}, \frac{7}{2}, \frac{9}{2}, \frac{11}{2}, \frac{1}{2}, \frac{3}{2}, \frac{5}{2}, \frac{7}{2}, \frac{3}{2}, \frac{5}{2}, \frac{7}{2}, \frac{1}{2}, \frac{3}{2} , \frac{1}{2} \}.$$
After reordering it, we obtain
$$\{\frac{11}{2}, \frac{9}{2}, \frac{7}{2}, \frac{5}{2},  \frac{3}{2}, \frac{1}{2},  \frac{7}{2}, \frac{5}{2}, \frac{3}{2},\frac{1}{2}, \frac{7}{2}, \frac{5}{2}, \frac{3}{2},\frac{1}{2}, \frac{3}{2},\frac{1}{2},\frac{3}{2}, \frac{1}{2}\}.$$
\begin{no}~\label{barorder}
Let 
$\lambda=(\lambda_1, \lambda_2,\ldots, \lambda_n) \in \mb R^n$.
We define $\overline{\lambda}$ to be the reordering of $\lambda$ such that
$$\overline{\lambda}_1 \geq \overline{\lambda}_2 \geq \ldots \geq \overline{\lambda}_n.$$
\end{no}
\section{Orderings and Reversal Phenomena}
Recall from \!~\cite{cm} that complex nilpotent orbits of classical groups have a partial ordering $\preceq$. 
\begin{no} Write
$\bold d_1 \preceq \bold d_2$ if $\mc O_{\bold d_1}$ is contained in the closure of $\mc O_{\bold d_2}$.
\end{no}
The ordering $\preceq$ for Young diagrams is different from the ordering $\preceq$ for numerical sequences (see Notation \!~\ref{preceq}). If one regards $\bold d_1^t$ and $\bold d_2^t$ as numerical sequences arranged in descending orders, then $\bold d_1 \preceq \bold d_2$ if and only if $\bold d_1^t \succeq \bold d_2^t$. 
\begin{thm}[Reversal Phenomena]~\label{rev}
Let $\mc O_{\bold D_1}$ and $\mc O_{\bold D_2}$ be nilpotent orbits of a fixed orthogonal group or symplectic group. Suppose that ${\bold d_1}^t$ and ${\bold d_2}^t$ are either both very even or both very odd. If 
$\bold d_1 \preceq \bold d_2$ then
$\overline{\mc I_{\pm}(\bold d_1)} \succeq \overline{\mc I_{\pm}(\bold d_2)}$ (see Notation \!~\ref{barorder}).
\end{thm}
 Later we will show that $\mc V(Ann\; \pi)=\mc O_{\bold d}$ for every $\pi \in \mc N(\bold D)$.
This theorem suggests that, among the representations $\mc N$ of a fixed group $G$,
the smaller the infinitesimal character, the bigger
is the associated variety and conversely. For example, the trivial representation has the smallest associated variety, namely $\{0\}$. Its 
infinitesimal character is given by $\rho(G)$, the greatest among all
$\mc I(\pi)$ with $\pi \in \mc N$. Theorem \!~\ref{rev} will be used in the Appendix to establish the estimates of $\mc I_{\pm}(\bold d)$ given in Theorem \!~\ref{infs} and Theorem \!~\ref{info}.\\
\\
We prove the reversal phenomena through a number of lemmas. {\bf The conditions in the theorem are always assumed in these lemmas}.
\begin{lem}
If the Young diagram $\bold d_1$ can be obtained from the Young diagram $\bold d_2$ by moving one domino to its lower left without changing other dominoes, then $\bold d_1 \preceq \bold d_2$. 
\end{lem}
We call this procedure a {\it move}.
\begin{lem}
If the Young diagram $\bold d_1$ can be obtained from $\bold d_2$ by a finite number of moves, then $\bold d_1 \preceq \bold d_2$ and vice versa.
\end{lem}
This lemma holds only under the assumption of the theorem. If $\bold d_1^t$ is very even and $\bold d_2^t$ is very odd, this lemma is no longer valid.  
\begin{lem}
If Young diagram $\bold d_1$ can be obtained from the Young diagram $\bold d_2$ 
by a move, then $\overline{\mc I_{\pm}(\bold d_1)} \succeq \overline{\mc I_{\pm}(\bold d_2)}$.
\end{lem}
This Lemma is an easy consequence of our algorithm  from $\bold d$ to $\mc I_{\pm}(\bold d)$. \\
\\
By the above three lemmas,
Theorem \!~\ref{rev} is proved. $\Box$ 
\begin{cor}
Let $\mc O_{\bold D_1}$ and $\mc O_{\bold D_2}$ be pre-rigid nilpotent orbits of a fixed orthogonal group or symplectic group (see Definition \!~\ref{prerigid}). Suppose that $\bold d_1^t$ and $\bold d_2^t$ are both very even or both very odd. Then
$\overline{\mc I_{\pm}(\bold d_1)} \succeq \overline{\mc I_{\pm}(\bold d_2)}$ if and only if
$\bold d_1 \preceq \bold d_2$.
\end{cor}
Proof: For pre-rigid orbits, one can reconstruct $\bold d^t$ from $\mc I(\bold d)$ in a unique way and the partial orderings are reversed in this reconstruction. $\Box$ 

\chapter{Associated Varieties and Existence of $\mc N(\mc O_{\bold D})$}
Let
$\mc O_{\bold D} \in \mc U(G)$ and $\pi \in \mc N(\mc O_{\bold D})$.
In this chapter, we study the matrix coefficients and the associated variety of $\pi$. We first bound $\overline{\mc I(\pi)}$ by a multiple of
$\rho(G)$.
By Theorem \!~\ref{kn}, we obtain a bound on the matrix coefficients of 
$\pi$. Then we apply Przebinda's theorem to show that
$$\mc V(Ann\;\pi)=cl(\mc O_{\bold d}).$$
With $\mc V(Ann\; \pi)$ in hand, we gain some control of $WF(\pi)$. Our next task is to prove Theorem \!~\ref{exi}.
Let me briefly describe the proof of Theorem \!~\ref{exi}.
Suppose $\pi(k) \in \mc N(\mc O(k))$. Then $\mc I(\pi(i))$ and $\mc V(Ann\;\pi(i))$ are all known for every $i \leq k$ based on our computation. If $G(\mc O(k))$ is orthogonal, then $\mc N(\mc O(k+1)) \neq \emptyset $ follows from Theorem \!~\ref{nonvanishing}. If $G(\mc O(k))$ is the metaplectic group, $\mc N(\mc O(k+1)) \neq \emptyset$ follows from
Theorem \!~\ref{non2} and Cor. \!~\ref{non20}. The details of the proof of \!~\ref{exi} are given at the end of this chapter. 

\section{Estimates on Infinitesimal Characters: I}
Let $G=Mp_{2n}(\mb R)$ and $\rho(G)=(n,n-1,\ldots, 1)$.
Let $\mc O_{\bold D}$ be in $\mc U(G)$. Let  
$$\bold d^t=(m_1 \geq m_2 \geq \ldots \geq m_{d_1})$$
Then we have $2n =\| \bold d\|$.
\begin{thm}~\label{infs}
Suppose $\mc O_{\bold D} \in \mc U(Mp_{2n}(\mb R))$. Let
$$\bold d^t=(m_1, m_2, \ldots, m_{d_1}).$$
Then
$$\overline{\mc I_{-}(\bold d)} \prec \frac{m_1+2}{2n} \rho(Mp_{2n}(\mb R)).$$
\end{thm}
A similar statement holds for $G=O(p,q)$.
\begin{thm}~\label{info}
Suppose $\mc O_{\bold d} \in \mc U(O(p,q))$. Let
$$\bold d^t=(m_1, m_2, \ldots, m_{d_1}).$$
Then
$$\overline{\mc I_{+}(\bold d)} \prec \frac{m_1+2}{p+q-2}(\frac{p+q}{2}-1, \frac{p+q}{2}-2, \ldots \frac{p+q}{2}-[\frac{p+q}{2}]).$$
\end{thm}
In fact, we will have that
\begin{equation}~\label{infstrong}
\begin{split}
\overline{\mc I_{-}(\bold d)} & \preceq  \frac{m_1}{2n} \rho(Mp_{2n}(\mb R)) \\
\overline{\mc I_{+}(\bold d)} & \preceq \frac{m_1}{p+q-2}
(\frac{p+q}{2}-1, \frac{p+q}{2}-2, \ldots \frac{p+q}{2}-[\frac{p+q}{2}]).
\end{split}
\end{equation}
for $\mc O_{\bold D} \in \mc U$. These two statements are slightly stronger than those of
Theorems \!~\ref{infs} and \!~\ref{info}. In the Appendix, we will give a proof for the first statement and Theorem \!~\ref{infs}. We skip the proof for Theorem \!~\ref{info}.
\section{Associated Varieties of $\mc N(\mc O_{\bold D})$}
\begin{thm}[Estimates on Leading Exponents]~\label{leadex}
Let $\mc O_{\bold D}$ be in $\mc U$. Let $\bold d^t=(m_1, m_2, \ldots, m_{d_1})$. Let $\pi \in \mc N(\mc O_{\bold D})$.
\begin{enumerate}
\item If $G(\mc O_{\bold D})$ is $Mp_{2n}(\mb R)$, then every leading exponent $v$ of $\pi$ satisfies
$$\Re(v) \prec (\frac{m_1+2}{2n}-1) \rho(Mp_{2n}(\mb R)).$$
\item If $G(\mc O_{\bold D})$ is $O(p,q)$, then every leading exponent $v$ of $\pi$ satisfies
$$\Re(v) \prec (\frac{m_1+2}{p+q-2}-1) \rho(O(p,q)).$$
\end{enumerate}
\end{thm}
Proof: Part (1) is a direct consequence of Theorem \!~\ref{infs} and Theorem \!~\ref{mat}. Part (2) is a direct consequence of Theorem \!~\ref{info} and Theorem \!~\ref{mat}. $\Box$ \\
\\
Theorem \!~\ref{leadex} can be established independently in the framework of \!~\cite{basic} without resorting to the infinitesimal character estimate in Theorems \!~\ref{infs} and \!~\ref{info}. However, the argument in \!~\cite{basic} relies heavily  on the estimates of twisted integral. It is more general but less intuitive. The estimate
of the the leading exponents using infinitesimal character only works for the unitarily small representation, namely unitary representations with infinitesimal character sitting inside the convex hull of $\{ w \rho \mid w \in W(\f g, \f h)) \}$. This method can be used to establish strongly semistable condition for unitarily small representations.
\begin{thm}~\label{as}
Let $\mc O_{\bold D} \in \mc U(G)$. Let $\pi \in \mc N(\mc O_{\bold D})$.
Then $\mc V(Ann\;\pi)=cl(\mc O_{\bold d})$.
\end{thm}
Proof: Notice that none of the operations in the definition of $\mc N(\mc O_{\bold D})$ changes $\mc I(\pi)$ or
$\mc V(Ann\;\pi)$. Without loss of generality, we assume
$\pi(k)=\theta_s(G(k-1), G(k))(\pi(k-1))$ and $\pi(d_1)=\pi$. 
Recall Theorem \!~\ref{ass1} and Theorem \!~\ref{ass2}.
\begin{enumerate}
\item Condition (1) is automatic by the definition of $(G(k-1), G(k))$. 
\item Conditions (2) from Theorem \!~\ref{ass1}
and \!~\ref{ass2} are satisfied due to Theorem \!~\ref{leadex}
and the definition of $(G(k-1), G(k))$.
\item $\pi(k-1)$ is not a discrete series representation. If $\pi(k-1)$ is in the discrete series and $G(k-1)$ is of rank $r$, then the infinitesimal character $\overline{\mc I(\pi(k-1))} \succeq (r-1, r-2, \ldots, 1,0)$ (for example, see Theorem 9.20 \!~\cite{knapp}). By Theorem \!~\ref{infs} and Inequalities (\!~\ref{infstrong}) this is only possible for several large $\mc I(\pi(k-1))$'s. For these large $\mc I(\pi(k-1))$'s, the representation $\pi(k-1)$ can not be in the discrete series.
\item Condition (4) is assumed.
\item Condition (5), which basically say that $\theta_s(\pi(k-1))$ is unitary, is readily verified by Theorem \!~\ref{un}.  
\item Conditions (6) are checked in Lemma \!~\ref{50} and Lemma \!~\ref{51}. 
\end{enumerate}
Thus Theorem \!~\ref{ass1} and Theorem \!~\ref{ass2} hold. Therefore, by Lemma \!~\ref{ass0} of DKP,
$$\mc V(Ann\;\pi(k))=\Theta(G(k-1), G(k))(\mc V(Ann\; \pi(k-1)))=cl(\mc O(k)).$$
By induction we have
 $\mc V(Ann\;\pi)=cl(\mc O_{\bold d}).$
$\Box$
\begin{cor}~\label{wavefront}
Fix an $\mc O_{\bold D} \in \mc U(Mp_{2n}(\mb R))$.
Let $O(p_0,q_0)=G(\mc O_{\bold D-1})$. Let $\pi \in \mc N_0(\mc O_{\bold D})$. Then there exists a nilpotent orbit $\mc O_{\bold D_0}$
in $WF(\pi)$ such that
\begin{enumerate}
\item $\bold d_0=\bold d$.
\item $G(\mc O_{\bold D_0-1}) = O(p_0,q_0)$.
\end{enumerate}
\end{cor}
Proof: Suppose that $\pi=\theta_s(p_0,q_0;2n)(\sigma)$ with $\sigma \in \mc N(\mc O_{\bold D-1})$. Then $\mc V(Ann\;\sigma)=\mc O_{\bold d-1}$ and $\mc V(Ann\;\pi)=\mc O_{\bold d}$. By Theorem \!~\ref{pr2.8}, $WF(\pi)$ must be in the image of $m_2$ with respect to $(O(p_0,q_0), Sp_{2n}(\mb R))$. Then Lemma \!~\ref{panlemma} says that every orbit $\mc O_{\bold S}$ in
$WF(\pi)$ satisfies $G(\mc O_{\bold S-1})=O(s,t)$ with $s \leq p_0$ and $t \leq q_0$. There must be at least one $\mc O_{\bold D_0}$ in $WF(\pi)$ such that
$G(\mc O_{\bold D_0-1})=O(p_0,q_0)$. Otherwise, $\mc V(Ann\;\pi)$, which is the complexification of $WF(\pi)$, will be strictly smaller than $cl(\mc O_{\bold d})$ and will not be equal to $cl(\mc O_{\bold D})$.
$\Box$
\section{Nonvanishing of $\mc N(\mc O_{\bold D})$}
Let $\mc O_{\bold D} \in \mc U$. Let $\pi \in \mc N_0(\mc O_{\bold D})$. We have proved that $\pi$ must be unitary and
$\mc V(Ann\;\pi)=\mc O_{\bold d}$. In other words, any representation in $\mc N(\mc O)$ must be unitary and $\mc V( Ann\; \pi)$ must be equal to $cl(\mc O_{\bold d})$. In this section, we will prove that $\mc N(\mc O_{\bold D})$ is nonempty. In fact, we will prove that $\mc N_0(\mc O_{\bold D})$ is nonempty.
\begin{thm}~\label{non3}
Let $\mc O_{\bold D} \in \mc U$.  Let $O(p_0,q_0)=G(k-1)$, $Mp_{2n_1}(\mb R)=G(k)$ and $O(p,q)=G(k+1)$.
Suppose 
$$\pi(k-1) \in \mc N(\mc O(k-1)), \mbox{ and } \pi(k)=\theta_s(p_0,q_0;2n_1)(\pi(k-1)) \neq 0.$$ Then
$\pi(k+1)=\theta_s(2n_1;p,q)(\pi(k)) \neq 0$.
\end{thm}
Proof: Since $\mc O_{\bold D} \in \mc U$, write $$\bold d(k+1)^t=(m_1 > m_2 \geq m_3 > m_4 \geq \ldots,).$$
Then 
$$\bold d(k)^t=(m_2 \geq m_3 > m_4 \geq \ldots, ), \qquad \bold d(k-1)^t=(m_3 > m_4 \geq \ldots,),$$
$$p+q=2n_1+m_1, \qquad 2n_1=p_0+q_0+m_2.$$
By Theorem \!~\ref{as}, $\mc V(Ann \, \pi(k))=\mc O_{\bold d(k)}$ and $\mc V(Ann \, \pi(k-1))=\mc O_{\bold d(k-1)}$.
Let $n_2=p+q-n_1-1$. Clearly $n_2-n_1=p+q-2n_1-1=m_1-1 > 0$ and $n_1+n_2+1 \geq p+q$. By Theorem \!~\ref{notzero} and Theorem \!~\ref{un}, it suffices to show that
$$\mc  Q(2n_1;p,q;2n_2)(\pi(k)) \neq 0.$$
\begin{enumerate}
\item By Corollary \!~\ref{wavefront}, there exists a nilpotent orbit $\mc O_{\bold D_0}$ of $Sp_{2n_1}(\mb R)$ such that
\begin{enumerate}
\item $\bold d_0= \bold d(k)$;
\item $\bold D_0-1$ contains $p_0$ \fp \ 's and $q_0$ \fm \ 's;
\item $\mc O_{\bold D_0} \subseteq WF(\pi(k))$.
\end{enumerate}
\item 
By Theorem \!~\ref{voganlemma},
$$ \Ind_{{\f{sp}}_{2n_1}(\mb R) \oplus \f{gl}(n_2-n_1, \mb R)}^{{\f{sp}}_{2n_2}(\mb R)} \tau(\mc O_{\bold D_0})
\subseteq WF(\Ind_{MP_{n_2-n_1,n_1}}^{Mp_{2n_2}(\mb R)} \chi^{\alpha} \otimes \pi^{\tau}). $$
Observe that 
$n_2-n_1=m_1-1$, and $m_1 -1 \geq m_2$. $m_2$ is the first entry of $\bold d_0^t=\bold d(k)^t$. 
By Theorem \!~\ref{tauinduction},
$$\Ind_{{\f{sp}}_{2n_1}(\mb R) \oplus \f{gl}(n_2-n_1, \mb R)}^{{\f{sp}}_{2n_2}(\mb R)} \tau(\mc O_{\bold D_0})$$
contains $m_1-m_2$ irreducible components
$$ \cup_{j=0}^{m_1-m_2-1} \mc O_{\bold S^{(j)}}.$$
And $\bold S^{(j)}$ can be obtained by
\begin{itemize}
\item merging $2$ columns of length ${m_1-1}$ to $\bold D_0$ from left;
\item extending the signs of $\bold D_0$ for the first $m_2$ rows;
\item assigning $j$ \fm \fp \ 's and $m_1-m_2-1-j$ \fp \fm \ 's to rows of length $2$. 
\end{itemize}
\item Notice that for every $j$,
$$(\bold s^{(j)})^t=\bold s^t= (m_1-1,m_1-1,m_2,m_3, \ldots, );$$
$$ (\bold s^{(j)}-\bold 1)^t=(\bold s-\bold 1)^t=(m_1-1, m_2,m_3, \ldots, ).$$
Consider the signature of
$\bold S^{(j)}-\bold 1$.
For the first $m_2$ rows, there are $p_0+m_2$ \fp \ 's and $q_0+m_2$ \fm \ 's.
For the last $m_1-m_2-1$ rows of length 1, there are $j$ \fp \ 's and $m_1-m_2-1-j$ \fm \ 's.
Thus the signature of $\bold S^{(j)}-1$ is 
$$(p_0+m_2+j, q_0+m_2+m_1-m_2-1-j)=(p_0+m_2+j, q_0+m_1-1-j)$$
\item Fix the parity of $p$. Let $j=p-p_0-m_2$. Here $p-p_0$ is equal to the number of {\fp}
in the first two columns of $\bold D(k+1)$. It is greater or equal to the number of boxes in the second column of $\bold D(k+1)$. It follows that $j=p-p_0-m_2 \geq 0$.
Since $p+q=m_1+m_2+p_0+q_0$,
 the signature of $\bold S^{(j)}- \bold 1$ equals 
$$(p_0+m_2+j, q_0+m_1-1-j)=(p,q_0+m_1-1+p_0+m_2-p)=(p, q-1).$$
From Lemma \!~\ref{panlemma}, $\mc O_{\bold S^{(j)}}$ 
occurs in the image of $m_2$ associated with $(O(p,q), Sp_{2n_2}(\mb R))$.
Also from Lemma \!~\ref{panlemma}, $\mc O_{\bold S^{(j)}}$  does not
occurs in the image of $m_2$ associated with 
$$ (O(p-2i,q+2i), Sp_{2n_2}(\mb R)), (O(p+2i,q-2i), Sp_{2n_2}(\mb R)) $$
for any $i \neq 0$.
\item 
By Theorem \!~\ref{un}, $\pi(k) \in R_{ss}(Mp_{2n_1}(\mb R), \omega(n_1+n_2+1,0;2n_1))$. So the generic nonvanishing theorem and consequently Cor. \!~\ref{non20} hold. Therefore,
$$\mc Q(2n_1;p,q; 2n_2)(\pi(k)) \neq 0.$$
and
 $$\pi(k+1)=\theta_s(2n_1;p,q)(\pi(k))\neq 0.$$
\end{enumerate}
$\Box$ \\
\\
Proof of Theorem \!~\ref{exi}: It suffices to prove that $\mc N_0(\mc O)$ is nonempty. We prove this by induction.  The existence of $\pi(1)$ is trivial. The existence of $\mc N_0(\mc O(2))$ is established by Jian-Shu Li. Suppose that $k \geq 2$ and $\pi(k) \in \mc N(\mc O(k))$. If $G(k)$ orthogonal, from Theorem \!~\ref{nonvanishing}, 
$\mc N(\mc O(k+1))$ is nonempty. If $G(k)$ is metaplectic, from Theorem  \!~\ref{non3},
$\pi(k+1)$ is not zero. Thus $\mc N_0(\mc O(k+1))$ is not empty. $\Box$
\chapter{Perspectives}\label{ch7}
In this paper, we only treat the group $O(p,q)$ and $Mp_{2n}(\mb R)$. With some modifications,  quantum induction can be defined for 
$$U(p,q), Sp(p,q), O^*(2n).$$ 
 We do not intend to work out the details here. Rather, we shall list some problems concerning quantum induction, unipotent representations and unitary dual. The main question is whether quantum induction can supplement parabolic induction (including complementary series) and cohomological induction (\!~\cite{k-v}) to produce a complete classification of $\Pi_u$, $\Pi_{auto}$ and $\Pi_{rama}$.
In any case, this
article should be viewed as a starting point for new development in that direction.\\
\\
Let us recall that unitary duals of the general linear groups $GL_n(\mathbb R)$, $GL_n(\mathbb C)$ and  $GL_n(\mathbb H)$ were classified by D. Vogan in 1980's (\cite{vo86}). Later, the unitary duals of $O(n, \mathbb C)$ and $Sp_{2n}(\mathbb C)$ were classified by  D. Barbasch (\cite{b89}).  The following are  the classical groups whose unitary dual is generally not classified:
$$ U(p,q), Sp(p,q), O^*(2n), Sp_{2n}(\mathbb R), O(p,q)$$
(except  for ranks $1$ or $2$ cases). We shall now call all the classical groups that are not general linear groups, classical groups of type I. Classical groups of type I share one common feature---they are the isometry groups of certain sesquilinear forms over $\mathbb R$, $\mathbb C$ or the quaternions $\mathbb H$ (\cite{li891}. The main point here is that quantum induction can be defined  for the classical groups for which the unitary dual is not classified. In what follows, we shall formulate our problems and conjectures mainly for $O(p,q)$ and $Mp_{2n}(\mb R)$ with the understanding that most of them can be formulated for type I classical groups like $U(p,q), Sp(p,q)$ and $O^*(2n)$.
\section{Wave Front Sets of $\mc N(\mc O)$}
Wave front sets under Howe's local theta correspondence are discussed in detail in \!~\cite{pr00}. There are still open questions which need to be answered. For the unipotent representations in $\mc N(\mc O)$, we have showed that their associated varieties are the complexification of $cl(\mc O)$. Now one may speculate that $cl(\mc O)$ is the wave from set of $\pi \in \mc N(\mc O)$.
This is far from the case. In fact, there are irreducible representations in $\mc N(\mc O)$ such that $WF(\pi) \neq cl(\mc O)$ as demonstrated by Theorem \!~\ref{twocom}. What should be true is that $\mc O \subseteq WF(\pi)$. We fall short of proving it.
\begin{conj}~\label{conj0}
For every $\pi \in \mc N(\mc O)$, $\mc O \subseteq WF(\pi)$.
\end{conj}
Let me give some hint on how one may proceed. 
Przebinda started a program attempting to obtain the Harish-Chandra character of $\theta_s(\pi)$ from the Harish-Chandra character of $\pi$. If his program works in our situation, we will have a confirmation of this conjecture. Let me point out a critical step.
\begin{conj}
Let $\pi \in R_{s}(MG_1, MG_2)$. Then
$$V(\theta_s(MG_1, MG_2)(\pi)) \cong [\Hom_{\f g_1, MK_1}(\omega^c, \pi)]_{MK_2}.$$
\end{conj}
Roughly speaking, this is saying that every $MK_2$-finite $\pi^c$-valued $MG_1$-equivariant distribution on $\omega$ can be obtained by integration over $MG_1$ as in \!~\cite{theta}. 
\begin{no}
Let $\mc O_{\bold D}$ be a nilpotent orbit of $G$. Let $\Pi_{\mc O_{\bold D}}(G)$ be the set of $\pi \in \Pi(G)$ such that $cl(\mc O_{\bold D})$ is of maximal dimension in $WF(\pi)$. Let $\Pi_{u,\mc O_{\bold D}}(G)=\Pi_{\mc O_{\bold D}}(G) \cap \Pi_u(G)$.
\end{no}
 According to a conjecture of Barbasch-Vogan (\!~\cite{bv}),
$$\{\pi \in \Pi_{\mc O_{\bold D}}(G) \mid \|\mc I(\pi)\| = \min\{ \|\mc I(\sigma)\| \mid \sigma \in \Pi_{\mc O_{\bold D}}(G) \} \subseteq \Pi_u(G).$$
Assume that Conjecture \!~\ref{conj0} is true. Then $\mc N(\mc O_{\bold D}) \subseteq \Pi_{u, \mc O_{\bold D}}(G)$. In fact, we conjecture that the representations in $\mc N(\mc O_{\bold D})$ are the ones that occur in the 
Barbasch-Vogan conjecture.
\begin{conj}~\label{conj5}
Let $G$ be either an orthogonal group or a metaplectic group. 
Let $\pi \in \Pi_{\mc O_{\bold D}}(G)$ with $\mc O_{\bold d}$ {\it rigid}. Let $\lambda$ be the shortest infinitesimal character
among all $\mc I(\pi)$ with $\pi \in \Pi_{\mc O_{\bold D}}(G)$.
Then $\lambda=\mc I(\bold d)$.
\end{conj}
This conjecture is false if $\mc O_{\bold d}$ fails to be rigid. 
Now assuming Conjecture \!~\ref{conj0}, we formulate the following conjecture which implies Conjecture ~\ref{conj5}.
\begin{conj}~\label{conj6}
 Let  $\mc O$ be a rigid orbit in $\mc U$. Suppose that  $G(1)$ is not $O(l) (l \geq 2)$. Then $\mc N(\mc O_{\bold D})$ exhausts all $\Pi_{u, \mc O_{\bold D}}(G)$ for $Mp_{2n}(\mb R)$ and $O(p,q)$.
\end{conj}
Notice that the rigid orbits in $\mc U$ include all special rigid orbits. 
For nonspecial rigid orbits of $O(p,q)$, by the works of Brylinski-Kostant (\!~\cite{bk}) and Huang-Li (\!~\cite{hl}), there might be unitary representations of the nontrivial covering of $SO_0(p,q)$ attached to them. Our construction does not cover their cases. However, for rigid nilpotent orbits in Conjecture \!~\ref{conj6}, our construction should be exhaustive.

\section{ Wave Front Sets and Induction Functors}
A critical observation in this paper is that $\bold D-2$ remains to be a symplectic (orthogonal) Young diagram if $\bold D$ is a sympletic (orthogonal) Young diagram.
On the philosophical level, this paper raises the question whether there always exists a functor
$$F(\mc O_{\bold D-2}, \mc O_{\bold D}): \Pi_{u, \mc O_{\bold{D-2}}}(G) \rightarrow \Pi_{u, \mc O_{\bold D}}(G).$$
This question does not always have a positive answer. For example, there are $\mc O_{\bold D}$ with $\Pi_{u, \mc O_{\bold D}}(G) = \emptyset$.
\begin{pro} For all real forms of symplectic and orthogonal groups,
construct functors
$$F(\mc O_{\bold D-2}, \mc O_{\bold D}): \Pi_{u, \mc O_{\bold{D-2}}}(G^0) \rightarrow \Pi_{u, \mc O_{\bold D}}(G)$$
whenever $\Pi_{u, \mc O_{\bold D}}(G) \neq \emptyset$. Here $G^0$ is a group attached to $\mc O_{\bold{D-2}}$ that is of the same type as $G$.
\end{pro}
We call each one of these functors, an {\it induction functor}.\\
\\
The second question is whether one can construct enough induction functors $\mc F(\mc O_{\bold D-2}, \mc O_{\bold D})$ such that all 
$\pi \in \Pi_{u, \mc O_{\bold D}}(G)$ can be obtained from $\Pi_{u,\mc O_{\bold D-2}}(G^0)$ through one of the induction functors.\\
\\
Indeed, if $\mc O_{\bold d}$ is induced from $\mc O_{\bold d-2}$, then unitary parabolic induction functor is a right candidate. Under different assumptions, complementary parabolic induction and cohomological induction can also be applied (\cite{bv83} \!~\cite{barbasch2} \!~\cite{k-v}). This situation has been studied intensively in the past. More recently, we gave an inductive construction of the discrete series representations of $U(p,q)$ by induction functors (\cite{he15}). It is feasible that similar ideas can be applied to produce enough induction functors $F(\mc O_{\bold D-2}, \mc O_{\bold D})$ when $\mc O_{\bold d}$ is induced from $\mc O_{\bold d-2}$.\\
\\
The case that $\mc O_{\bold D}$ is not induced from $\mc O_{\bold D-2}$ is the main focus of this paper. It seems that the quantum induction $\mc Q$ is the right choice. For example, one can study the following problem.
\begin{pro} Consider $Mp_{2n}(\mb R)$.
Let $\mc O_{\bold D}$ be a symplectic nilpotent orbit. Suppose that $\mc O_{\bold D}$ is not induced from $\mc O_{\bold D-2}$. Let $2n_2 = \| \bold d \|$, $(p,q)=((\bold D-1)^+, (\bold D-1)^-)$ and $ 2n_1=\| \bold d-2 \|$. Then is  
$$\mc Q(2n_1; p,q; 2n_2): \Pi_{u,\mc O_{\bold D-2}}(G) \rightarrow \Pi_{u, \mc O_{\bold D}}(G)?$$
\end{pro}
Notice that $\Pi_{u, \mc O_{\bold D}}(G)$ may be empty for some $\bold D$. Even in those cases, this problem seems to make sense, that is, $\mc Q(2n_1; p,q; 2n_2)$ will conjecturally map $\Pi_{u,\mc O_{\bold D-2}}(G)$ to zero. We shall come back to this problem in a subsequent paper.

\section{Unitary Dual, Automorphic Dual and Ramanujan Dual}
Let $G$ be a classical group of type I.
Let $\Pi_{herm}(G)$ be the set of irreducible Harish-Chandra modules with an invariant Hermitian structure. Let $\Pi_{temp}(G)$ be the set of irreducible
Harish-Chandra modules with almost $L^2$ matrix coefficients. $\Pi_{temp}(G)$ is the same as $\Pi_2(G)$, the reduced dual of $G$.
Recall that
$$\Pi(G) \supset \Pi_{herm}(G) \supset \Pi_u(G) \supset \Pi_{temp}(G). $$
The classification of the admissible dual $\Pi(G)$ is due to Langlands. Langlands proved that every $\pi \in \Pi(G)$ occurs as the unique quotient of certain induced representation from a tempered representation. The quotients are often known as Langlands quotients. The classification of $\Pi_{temp}$ is due to Knapp-Zuckerman. The classification of the Hermitian dual is also known. See for example Theorem 16.6 and its remark \!~\cite{knapp}.\\
\\
With these classifications in hand, classification of unitary dual can be translated as the determination of unitarity of the Langlands quotients. But this approach is of great mathematical complexity. It often involves the positivity of certain analytically defined intertwining operator. 
In \!~\cite{vogan86} and \!~\cite{vogan86book}, Vogan envisioned a more geometric approach based on Kirillov-Kostant's orbit philosophy. In the Interlude in \cite{vo87}, a program for constructing unitary representations is described. Starting with a unipotent representation, one may apply the construction of unitary parabolic induction, cohomological induction and conplementary series  to obtain all irreducible unitary representations. This approach turned out to be successful for the general linear groups since all unipotent representations for the general linear groups can be obtained by unitary parabolic induction (\cite{vo86}). Howerver, for other classical groups, unitary parabolic induction is not sufficient to obtain all unipotent representations. It fails to produce representations associated with rigid orbits. In this paper, we supplied a list of unipotent representation for special rigid orbits. For other orbits, the complete list of unipotent representations will surely require unitary parabolic induction. Not only that, it seems that cohomological induction is also needed. In any case, we state Vogan's program in the following problem.

\begin{pro}[Vogan] Give a complete list of unipotent representations of $G$. Construct unitary representations using parabolic induction and cohomological induction from a unipotent representation.
\end{pro} 

Recently, Salamanca-Riba and Vogan conjectured that all irreducible unitary representation can be constructed by cohomological method from an irreducible unitary representation of a smaller group except unitarily small representations \!~\cite{sv}. Notice that all the representations in $\mc N(\mc O)$ are unitarily small and quantum inductions  keep the infinitesimal characters small. So quantum induction, combined with parabolic induction, can be used to construct unitarily small representations. Whether it is sufficient to produce all unitarily small representations, is yet to be seen. Nevertheless, quantum induction certainly raises our hope for a classification of the unitary dual of classical groups of type I, modulo the conjecture of Salamanca-Riba and Vogan. For coverings of classical groups, the situation is a little more delicate. But the same idea should apply.\\
\\
Alternatively, one can follows the route of Problem 1 and construct the unitary dual by applying different induction functors based on the wave front sets.
\begin{pro}
Arrange a process to construct irreducible unitary representations of a classical group $G$ using unitary parabolic induction,
cohomological induction (\!~\cite{k-v}) and complementary series construction, combined with
quantum induction. Can this process be arranged to produce all irreducible unitary representations?
\end{pro}
For an illustration, let $\mc O_{\mathbf D}$ be a real nilpotent coadjoint orbit of $Sp_{26}(\mb R)$, with $\bold D=$
\begin{center}
\fm \fp \fm \fp \fm \fp  \\
\fm \fp \fm \fp \fm \mbo \\
\fp \fm \fp \fm \fp  \mbo \\
\fm \fp \fm \mbo \mbo \mbo  \\
\fp \fm \fp \mbo \mbo \mbo  \\
\fp \fm \mbo \mbo \mbo \mbo \\
\fp \fm \mbo \mbo \mbo \mbo 
\end{center}
Then $\mc O_{\mathbf D-2}$ is a nilpotent orbit of $Sp_{12}(\mathbb R)$, $\mc O_{\mathbf D-4}$ is a nilpotent orbit of $Sp_4(\mathbb R)$. $\mc O_{\mathbf D-4}$ is a minimal orbit of $Sp_{4}(\mb R)$. There are essentially two representations $\pi$ of $Mp_4(\mb R)$ associated with  $\mc O_{\mathbf D-4}$---both come from the oscillator representation. To obtain representations attached to $\mc O_{\mathbf D-2}$, one applies the quantum induction $Q(4; 4,4; 12)$ to $\pi$. By our results, we obtain a unitary representation of $Mp_{12}(\mb R)$. Now notice that $\mc O_{\mathbf D}$ is contained in the induced orbit 
$${\rm Ind}_{{\f gl}_7(\mb R) {\f sp}_{12}(\mb R)}^{\f {sp}_{26}(\mb R)} \mc O_{\mathbf D-2}$$ 
One will need to apply  the cohomological induction to obtain a representation associated with $\mc O_{\mathbf D}$. Parabolic induction can also be used in this context. However, the wave front set of the resulting representations will have several components and will contain $\mc O_{\mathbf D}$ as one irreducible component. \\
\\
Many questions remain, for example, whether cohomological induction is exhaustive. For the group $U(p,q)$, Vogan conjectured that unitary cohomological induction is sufficient to produce all unitary representations with integral infinitesimal characters. In light of this conjecture,
we believe that quantum induction should be sufficient to supplement
the existing techniques to produce a classification of the unitary dual for type I classical groups.
\\
\\
Theta correspondences as formulated by Howe (\!~\cite{howe79}) originated from the theory of theta series (\!~\cite{weilsiegel}, \!~\cite{si}). The history of theta series went back to the theta functions defined by Jacobi. Special cases of theta series were studied by many people in number theory.  The upshot is that theta correspondences should produce automorphic representations of a bigger group from automorphic representations of a smaller group. This line of ideas was further explored by Rallis, Li in some quite general settings (\!~\cite{ra87}, \!~\cite{li94}) and by many other people in special cases. We are tempted to make the following conjecture.
\begin{conj}
All representations in $\mc N(\mc O)$ are automorphic.
\end{conj}
In fact, we can even say more about $\mc N(\mc O)$ and quantum induction. Recall that
$$\Pi_u(G) \supset \Pi_{auto}(G)  \supset \Pi_{rama}(G).$$
We adopt the definitions from \!~\cite{bls}. 
In principle, according to the conjectures of Ramanujan
and Selberg (\!~\cite{se}, \!~\cite{sa}), complementary series should not occur in $\Pi_{auto}(G)$ or $\Pi_{rama}(G)$. Assuming that, 
\begin{pro}
Can  $\pi \in \Pi_{auto}(G)$ be constructed by unitary quantum induction, unitary cohomological induction and unitary parabolic induction?
\end{pro}
The argument that parabolic inductions should send $\Pi_{auto}$ of a Levi subgroup to $\Pi_{auto}(G)$ was given in \!~\cite{bls} with their $H$ taking to be 
the parabolic subgroup rather than the Levi subgroup. Finally, for split classical groups, Barbasch showed that the irreducible spherical unitary representations can be constructed as a parabolically induced representation from a complementarily induced representation tensored with
a special unipotent representation (\!~\cite{barbasch2}). The Ramanujan dual is simply the intersection of the automorphic dual with the spherical dual. Thus one may ask
\begin{pro}
For $G$ a split classical group of type I, can one construct all $\pi \in \Pi_{rama}(G)$ by unitary quantum induction and unitary parabolic induction?
\end{pro}
This exhaustion question is perhaps too difficult to answer.  
For $SL(2, \mb R)$ where no quantum induction is involved, this is equivalent to the Ramanujan-Selberg conjecture. We shall refer the readers to the recent book of Arthur for more conjectures concerning unitary dual and automorphic dual (\cite{ar13}).

\chapter*{Appendix}
\addcontentsline{toc}{chapter}{\bf{Appendix}  }

\section*{Infinitesimal Characters for $Mp$}
\addcontentsline{toc}{section}{\bf{Infinitesimal Characters for $Mp$}  }
Recall that
$\mc I_{\pm}(\bold d)$ consists of an arithmetic sequence with multiplicities.
For $\bold d^t$ odd, the sequence ends with $\frac{1}{2}$; for $\bold d^t$ even, the sequence ends with $zero$ if $\bold d \neq [1]^{2n}$. In this chapter, we shall give a proof for Theorem \!~\ref{infs}. We start with a lemma.

\begin{lem*}[Triviality]\label{trivial}
\begin{enumerate}
\item If $\overline{\lambda} \prec \overline{\mu}$ and $\overline{\lambda^{\p}} \preceq \overline{\mu^{\p}}$, then
$\overline{\lambda} + \overline{\lambda^{\p}} \prec \overline{\mu}+ \overline{\mu^{\p}}.$
\item If $\overline{\lambda} \prec \overline{\mu}$ and $\overline{\lambda^{\p} }\prec \overline{\mu^{\p}}$, then
$\overline{(\lambda, \lambda^{\p})} \prec \overline{( \mu, \mu^{\p})}$.
\item If $\lambda_l \leq \mu_l$ for all $l$, then
$\overline{\lambda} \preceq \overline{\mu}$.
\item If $\lambda_l < \mu_l$ for all $l$, then $\overline{\lambda} \prec \overline{\mu}$.
\end{enumerate}
\end{lem*}
Let $G=Sp_{2n}(\mb R)$.
Then
$$\rho(G)=(n,n-1,\ldots, 1)=\mc I_{-}(2n).$$
\begin{lem*}\label{i2}
Suppose that $m \equiv r \pmod 2$. If $r \leq m$, then
$$\overline{(\mc I_{-}(m), \mc I_{+}(r))} \preceq \frac{m}{m+r} \mc I_{-}(m+r).$$
\end{lem*}
Proof: Write
$$B= \overline{(\mc I_{-}(m), \mc I_{+}(r))}, \qquad A= \frac{m}{m+r} \mc I_{-}(m+r).$$
We prove our lemma by induction on $\frac{m-r}{2}$.  If $\frac{m-r}{2}=0$, $m=r$. Then
$$A=(\frac{m}{2},\frac{m-1}{2}, \frac{m-2}{2}, \frac{m-3}{2}, \ldots, \frac{2}{2}, \frac{1}{2}),$$
$$B=(\frac{m}{2},\frac{m-2}{2}, \frac{m-2}{2}, \frac{m-4}{2}, \frac{m-4}{2}, \ldots).$$
$B$ ends with $(1,0)$ if $m$ is even and $(\frac{1}{2},\frac{1}{2})$ if $m$ is odd. Obviously, $A_k \geq B_k$.
So $B \preceq A$.
Assume that $m-2 \geq r$ and
$$\overline{(\mc I_{-}(m-2), \mc I_{+}(r))} \preceq \frac{m-2}{m-2+r} \mc I_{-}(m+r-2).$$
Notice that $\frac{m-2}{m-2+r} \leq \frac{m}{m+r}$. We have
\begin{equation*}
\begin{split}
B= \overline{(\mc I_{-}(m), \mc I_{+}(r))} & = \overline{(\frac{m}{2}, \mc I_{-}(m-2), \mc I_{+}(r))} \\ 
 & \preceq \overline{(\frac{m}{2}, \frac{m-2}{m-2+r} \mc I_{-}(m+r-2))} \\
 & \preceq \overline{(\frac{m}{m+r} \frac{m+r}{2}, \frac{m}{m+r} \mc I_{-}(m+r-2)) } \\
 & = \frac{m}{m+r} \overline{\mc I_{-}(m+r)}=A.
\end{split}
\end{equation*}
$\Box$
\begin{lem*}~\label{i3} Consider a partition of $2n$
$$\bold d^t= (\overbrace{m,m,m-2,m-2,\ldots, m-2j+2,m-2j+2}^{2j},m_0,r)$$ with 
$m \equiv m_0 \equiv r \pmod 2$ and $m-2j \geq m_0 \geq r \geq 0$.  Then
$$\overline{\mc I_{-}(\bold d)}
\preceq \frac{m}{2n} \mc I_{-}(2n).$$
\end{lem*}
Proof:  Clearly, $2n=m_0+r+(2m-2j+2)j$. If $j=0$, $\bold d^t=(m_0,r)$. We have $\mc I_{-}(\bold d) \preceq \frac{m_0}{m_0+r} \mc I_{-}(m_0+r)$ by Lemma \!~\ref{i2}. Assume $j \geq 1$. \\
\\
Use induction on $m$. When $m=2$, $B(2,0,0)=\mc I_{-}(2) \preceq I_{-}(2)$. So our lemma holds for $m=2$. Assume our lemma holds for
$m-1$. \\
\\
Suppose that $ r \geq 1 $. The case $r=0$ can be treated similarly.  Obviously,
$2n < 2m (j+1)$. So 
$$ m(2n-2j-2)= 2mn-m(2j+2) \leq 2nm-2n=2n(m-1).$$ It follows that $\frac{m}{2n} \leq
\frac{m-1}{2n-2j-2}$. Consider 
$$\bold e^t= (\overbrace{m-1,m-1,m-3,m-3,\ldots,m-2j+1,m-2j+1}^{2j},m_0-1,r-1).$$ \\
\\
Suppose that $m$ is even.
Then 
$\overline{\mc I_{-}(\bold d)}=(\overline{\mc I_{-}(\bold e)}+ \bold{\frac{1}{2}}, \overbrace{0,0,\ldots, 0}^{j+1})$.
 By induction hypothesis,
$$ \overline{\mc I_{-}(\bold e)} \preceq  \frac{m-1}{2n-2j-2} \mc I_{-}(2n-2j-2).$$
Then we have
\begin{equation*}
\begin{split}
 \overline{\mc I_{-}(\bold e)} + \bold{\frac{1}{2}} \preceq & \frac{m-1}{2n-2j-2} \mc I_{-}(2n-2j-2) +\bold{ \frac{1}{2}} \\
= & \bold {\frac{m-1}{2}}-\frac{m-1}{2n-2j-2}(0,1, \ldots n-j-2)+ \bold{\frac{1}{2}} \\
\preceq & \bold{\frac{m}{2}}-\frac{m}{2n}(0,1, \ldots, n-j-2) \\
= & \frac{m}{2n} \bold n-\frac{m}{2n}(0,1, \ldots, n-j-2) \\
= & \frac{m}{2n} (\overbrace{n,n-1, \ldots,j+2}^{n-j-1}).
\end{split}
\end{equation*}
We obtain $\overline{\mc I_{-}(\bold d)} \preceq \frac{m}{2n} \mc I_{-}(2n)$. \\
\\
Suppose $m$ is odd.
Then 
$\overline{\mc I_{-}(\bold d)}=(\overline{\mc I_{-}(\bold e)}+ \bold{\frac{1}{2}}, \overbrace{\frac{1}{2},\frac{1}{2},\ldots, \frac{1}{2}}^{j+1})$. We will again have
$$\overline{\mc I_{-}(\bold e)} + \bold{\frac{1}{2}} \preceq \frac{m}{2n} (\overbrace{n,n-1, \ldots,j+2}^{n-j-1}).$$
As $l$ changes from $n-j$ to $n$, $\frac{m}{2n} \mc I_{-}(2n)_{l} - \frac{1}{2} $ decreases from 
positive to negative. So
$$\frac{m}{2n} \sum_{i=1}^{l} \mc I_{-}(2n)_{i}- \sum_{i=1}^l \overline{\mc I_{-}(\bold d)}_i$$
increases and then decreases. It suffices to show that
$$\frac{m}{2n} \sum_{i=1}^{n} \mc I_{-}(2n)_{i}- \sum_{i=1}^n \overline{\mc I_{-}(\bold d)}_i \geq 0$$
Notice that 
$$\sum_{i} (\mc I_{-}(m-2s))_i=\frac{(m-2s+1)^2}{8}, \qquad \sum_i (\mc I_{+}(m-2s))_i=\frac{(m-2s-1)^2}{8}$$
We obtain
\begin{equation*}
\begin{split}
 & \frac{m}{2n} \sum_{i=1}^{n} \mc I_{-}(2n)_{i}- \sum_{i=1}^n \overline{\mc I_{-}(\bold d)}_i \\
= & \frac{m}{2n} \frac{(n+1)n}{2}-[ \sum_{s=0}^{j-1} \frac{(m-2s+1)^2}{8}+\frac{(m-2s-1)^2}{8}] -
\frac{(m_0+1)^2}{8}-\frac{(r-1)^2}{8} \\
= & \frac{1}{8} [ m(2n+2)- (\sum_{s=0}^{j-1} 2 (m-2s)^2+2)-(m_0+1)^2-(r-1)^2 ] \\
= & \frac{1}{8} \{ m ( m_0+1+r-1+2+  \sum_{s=0}^{j-1} 2 (m-2s))- (m_0+1)^2-(r-1)^2-2j - \sum_{s=0}^{j-1} 2 (m-2s)^2 \}  \\
=& \frac{1}{8} \{ (m_0+1)(m-m_0-1)+(r-1)(m-r+1)+ 2 (m-j)+ 2 \sum_{s=0}^{j-1} (m-2s)(m-m+2s) \} \\
\geq & 0
\end{split}
\end{equation*}
Therefore $\overline{\mc I_{-}(\bold d)} \preceq \frac{m}{2n} \mc I_{-}(2n)$. By induction, $\overline{\mc I_{-}(\bold d)} \preceq \frac{m}{2n} \mc I_{-}(2n)$ for all $m$.
$\Box$ \\
\\
Proof of Theorem \!~\ref{infs}: Fix $n$ first. Let $\mc O_{\bold d} \in \mc U(Mp_{2n}(\mb R))$. Fix the number of row $m_1=m$.
Consider 
\begin{equation*}~\label{constraint}
\bold d_0^t=(m,m,m-2,m-2, \ldots, m-2j+2, m-2j+2, m_0,r) \qquad (j \geq 0)
\end{equation*}
with 
\begin{equation*}~\label{constraint0}
0 \leq r \leq m_0 \leq m-2j-2, \qquad m \equiv m_0 \equiv r \pmod 2.
\end{equation*}
Here $j=0$ means that $\bold d_0^t=(m_0, r)$. By Lemma \!~\ref{i3},
$$\overline{\mc I_{-}(\bold d_0)} \preceq \frac{m}{2n} \mc I_{-}(2n).$$
Notice that $\mc O_{\bold d_0}$ is the minimal orbit in $\mc U(Mp_{2n}(\mb R))$ with a fixed $m$. Therefore, by Theorem \!~\ref{rev},
for any $\bold d$ with $m$ rows, $\overline{\mc I_{-}(\bold d)} \preceq \overline{\mc I_{-}(\bold d_0)}$.   We obtain for any $\mc O_{\bold d} \in \mc U(Mp_{2n}(\mb R))$,
$$\overline{\mc I_{-}(\bold d)} \preceq \frac{m_1}{2n} \rho(G)
\prec \frac{m_1+2}{2n} \rho(G).$$
$\Box$ 
\section*{Integration of Functions with Values in Hilbert Spaces }
\addcontentsline{toc}{section}{\bf{Integration of Functions with Values in Hilbert Spaces}  }
In this section, we shall give a proof of Lemma \!~\ref{fa}. By changing coordinates, it suffices to prove the following Lemma.
\begin{lem*}~\label{ivh} Let $\mc H$ be a Hilbert space and $\Phi$ be a $\mc H$-valued continuous function on $\mb R^n$. If $$\int \int_{x,y \in \mb R^n}  | (\Phi(x), \Phi(y))| \, d x \, d y < \infty,$$
 then
$\int_{x} \Phi(x) \, d x \in \mc H$ is well-defined as the unique limit of
$u_{K_i}=\int_{K_i} \Phi(x) d x $ for any increasing sequence of compact sets $K_i \rightarrow \mb R^n$. Furthermore,
$$(\int_{x \in \mb R^n} \Phi(x) \,  d x, \int_{y \in \mb R^n} \Phi(y) \, d y)=\int_{x,y \in \mb R^n}  (\Phi(x), \Phi(y)) \, d x \, d y$$
\end{lem*}
Proof: Let $\Phi(x)$ be a continuous function from $\mathbb R^n $ to $\mathcal H$ such that
$$\int_{x,y \in \mathbb R^n} |(\Phi(x), \Phi( y))| \, d x \, d y < \infty.$$
Let $K_i$ be an increasing sequence of compact sets $K_i \rightarrow \mathbb R^n$.
Then as $ i \rightarrow \infty$,
$$\int_{\mathbb R^{n} \times \mathbb R^n - K_l \times K_l} |(\Phi(x), \Phi(y))| dx dy \rightarrow 0.$$
Define $u_{K_i}= \int_{K_i} \Phi(x) d x.$
Then for any $m >i$, we have
\begin{equation*}
\begin{split}
 (u_{K_m}-u_{K_i}, u_{K_m} -u_{K_i})= &(\int_{K_{m}-K_{i}} \Phi(x) d x, \int_{K_{m}-K_{i}} \Phi(x) \, d x)  \\
 = & \int_{(K_{m}-K_{i}) \times (K_{m}-K_{i})}(\Phi(x), \Phi(y)) \, d x \, d y \\
\leq &  \int_{(\mathbb R^n -K_i) \times (\mathbb R^n - K_i)} |(\Phi(x), \Phi(y))| \, d x \, d y \\
\leq & \int_{\mb R^n \times \mb R^n - K_i \times K_i} |(\Phi(x), \Phi(y))| \, d x \, d y \rightarrow 0.
\end{split}
\end{equation*}
Therefore $\{ u_{K_i} \}$ is a Cauchy sequence and $\lim u_{K_i}$ exists. So
$u= \int \Phi(x) d x \in \mc H$ is well-defined. Furthermore,
$\| u-u_{K_l} \| \rightarrow 0$.
By a similar argument, we have
$$(u,u )=\int (\Phi(x),\Phi(y)) \, d x \, d y.$$
We can generalize Lemma \!~\ref{ivh} as follows.
\begin{thm*}~\label{ivh1} Let $H$ be a Hilbert space and $X$ a locally compact space. Let $\mu$ be a $\sigma-$finite regular Borel measure. Let $\Phi$ be a $\mc H$-valued continuous function on $X$. If $$\int \int_{x,y \in X}  | (\Phi(x), \Phi(y))| \, d \mu(x )\, d \mu (y) < \infty,$$
 then
$\int_{x} \Phi(x) \, d \mu(x) \in \mc H$ is well-defined as the unique limit of $\int_{K_i} \Phi(x) d \mu(x)$ for any increasing sequence $K_i$ with $\mu(X-\cup K_i)=0$. Furthermore,
$$(\int_{x \in X} \Phi(x) \,  d \mu(x), \int_{y \in X} \Phi(y) \, d \mu(y)=\int_{x,y \in X}  (\Phi(x), \Phi(y)) \, d \mu(x) \, d \mu(y)$$
\end{thm*}
The convergence of $\int_{x} \Phi(x) \, d \mu(x)$ is neither a strong convergence nor a weak convergence. We refer to it as {\it convergence on compact sets}.

\end{document}